\documentclass[11pt]{article}

\usepackage{amssymb,amscd,amsthm,amsmath}

\newtheorem{theorem}{Theorem}[section]
\newtheorem{lemma}[theorem]{Lemma}

\newtheorem{cor}[theorem]{Corollary}

\theoremstyle{definition}
\newtheorem{definition}[theorem]{Definition}
\newtheorem{example}[theorem]{Example}

\theoremstyle{remark}
\newtheorem{remark}[theorem]{Remark}

\theoremstyle{conjecture}
\newtheorem{conjecture}[theorem]{Conjecture}
\theoremstyle{problem}

\numberwithin{equation}{section}

\newcommand{\C}{\mathbb{C}}

\newcommand{\K}{\mathbb{K}}
\newcommand{\M}{\mathbb{M}}
\newcommand{\N}{\mathbb{N}}
\newcommand{\Q}{\mathbb{Q}}
\newcommand{\R}{\mathbb{R}}
\newcommand{\T}{\mathbb{T}}
\newcommand{\Z}{\mathbb{Z}}
\newcommand{\cB}{\mathcal{B}}

\newcommand{\cI}{\mathcal{I}}
\newcommand{\cH}{\mathcal{H}}

\newcommand{\cO}{\mathcal{O}}
\newcommand{\cT}{\mathcal{T}}

\newcommand{\cP}{\mathcal{P}}
\newcommand{\cQ}{\mathcal{Q}}
\newcommand{\cZ}{\mathcal{Z}}

\newcommand{\fo}{\mathfrak{o}}

\newcommand\Aut{\operatorname{Aut}}

\newcommand\Tr{\operatorname{Tr}}
\newcommand\id{\operatorname{id}}
\newcommand\Ad{\operatorname{Ad}}
\newcommand\Hom{\operatorname{Hom}}
\newcommand\Map{\operatorname{Map}}

\newcommand\Tor{\operatorname{Tor}}

\newcommand{\AC}{\operatorname{AC}}

\newcommand{\salpha}{\alpha^{\mathrm{s}}}
\newcommand{\sbeta}{\beta^{\mathrm{s}}}

\newcommand{\sA}{A^{\mathrm{s}}}

\newcommand{\tf}{\tilde{f}}
\newcommand{\tF}{\tilde{F}}

\newcommand{\tPhi}{\tilde{\Phi}}
\newcommand{\tPsi}{\tilde{\Psi}}
\title{Poly-$\mathbb{Z}$ group actions on Kirchberg algebras II}

\author{Masaki Izumi
\thanks{Supported in part by JSPS KAKENHI Grant Number JP15H03623}\\
Graduate School of Science \\
Kyoto University \\
Sakyo-ku, Kyoto 606-8502, Japan 
\and
Hiroki Matui 
\thanks{Supported in part by JSPS KAKENHI Grant Number JP18K03321} \\
Graduate School of Science \\
Chiba University \\
Inage-ku, Chiba 263-8522, Japan}

\begin{document} 
\maketitle
\begin{abstract} This is the second part of our serial work on the classification of poly-$\mathbb{Z}$ group actions on Kirchberg algebras. 
Based on technical results obtained in our previous work, we completely reduce the problem to the classification of continuous fields of Kirchberg algebras 
over the classifying spaces. 
As an application, we determine the number of cocycle conjugacy classes of outer $\mathbb{Z}^n$-actions on the Cuntz algebras. 
\end{abstract}
\tableofcontents
\section{Introduction} 
\bigskip
This is a continuation of our previous work \cite{IMI} on the classification of discrete amenable group actions on Kirchberg algebras.  
Kirchberg algebras are surely the most prominent class among classifiable amenable $C^*$-algebras, and it should be the first place to work on 
in order to establish classification theory of amenable group actions comparable to the case of von Neumann algebras, where satisfactory 
classification theory is known. 
The reader is referred to \cite{P00} and \cite{R2002} for the classification of Kirchberg algebras, 
and to the introduction of the first part \cite{IMI} for a brief account of the classification 
of group actions on classifiable $C^*$-algebras. 
Here we just give an incomplete list \cite{AK17JFA}, \cite{GI11}, \cite{I04Duke}, \cite{I04Adv}, \cite{IM2010}, \cite{KM08}, \cite{Ki95crelle}, \cite{Ki96}, 
\cite{Ki98JOP}, \cite{Ki98JFA}, \cite{Ma08}, \cite{MS14}, \cite{N99}, \cite{N00} of classification results on discrete amenable group 
actions on classifiable $C^*$-algebras.

Based on technical results obtained in our accompanying papers \cite{IMW}, \cite{IMI}, in this paper we propose a topological approach to 
the classification problem, which we hope will be a guiding principle for more general problems in the future. 
We recall the notion of classifying spaces in algebraic topology first.  
For any topological group $G$, there exists a universal principal $G$-bundle $EG\to BG$ 
satisfying the following property: every numerable principal $G$-bundle $\cP\to X$ is isomorphic to the pullback bundle 
$f^*EG$ of a continuous map $f:X\to BG$ so that the set of isomorphism classes of numerable principal 
$G$-bundles over $X$ is in one-to-one correspondence with the homotopy set $[X,BG]$ (see \cite[Chapter 4]{Hu}). 
The space $BG$, which is unique up to homotopy equivalence by universality, is called the classifying space of $G$. 
Since the Milnor construction of $BG$ is functorial, a continuous group homomorphism $h:G_1\to G_2$ induces a continuous map 
$Bh:BG_1\to BG_2$. 
If moreover $G_1$ and $G_2$ are discrete groups, the map 
$$\Hom(G_1,G_2)/{\textrm{conjugacy}}\ni [h]\mapsto [Bh]\in [BG_1,BG_2],$$
is a bijection, which follows from the classification of regular covering spaces over $BG_1$ (see for example \cite[Section 3.7, 3.8]{May99}). 

The main subject of this paper is the following conjecture, which is modeled after the above fact, 
on discrete amenable group actions on Kirchberg algebras.  

\begin{conjecture}\label{main conjecture s} Let $G$ be a torsion-free countable amenable discrete group, and let 
$A$ be a stable Kirchberg algebra. 
Then the map 
$$\cB:\mathcal{OA}(G,A)/{\textrm{cocycle conjugacy}} \ni [\alpha]\mapsto [B\alpha]\in [BG,B\Aut(A)]$$
is a bijection, where $\mathcal{OA}(G,A)$ denotes the set of outer actions of $G$ on $A$. 
\end{conjecture}

This is the most optimistic version of conjectures of similar kind,  and to make it more tractable 
we could assume in addition that $BG$ is homotopic to a finite CW-complex. 

We introduce necessary terminology and notation to understand the precise statement, which will be used throughout this paper. 
Let $\alpha$ and $\beta$ be actions of a discrete group $G$ on a $C^*$-algebra $A$. 
A family of unitaries $\{u_g\}_{g\in G}$ in the unitary group $U(M(A))$ of the multiplier algebra 
$M(A)$ of $A$ is said to be an $\alpha$-cocycle if it satisfies 
the 1-cocycle relation $u_{gh}=u_g\alpha_g(u_h)$ for any $g,h\in G$. 
A cocycle perturbation $\alpha^u$ of $\alpha$ is a $G$-action defined by $\alpha^u_g=\Ad u_g\circ \alpha_g$ 
with an $\alpha$-cocycle $\{u_g\}_{g\in G}$. 
We say that $\alpha$ and $\beta$ are \textit{cocycle conjugate}, 
if there exist $\theta\in \Aut(A)$ and an $\alpha$-cocycle $\{u_g\}_{g\in G}$ satisfying 
$\theta\circ \beta_g\circ \theta^{-1}=\alpha^u_g$. 
If moreover we can choose $\theta$ so that the $KK$-class $KK(\theta)$ of $\theta$ is the same as the $KK$-class 
$KK(\id_A)$ of the identify, we say that $\alpha$ and $\beta$ are \textit{$KK$-trivially cocycle conjugate}.  
We say that $\alpha$ is outer if $\alpha_g$ is outer, (i.e. not of the form $\Ad u$ with $u\in U(M(A))$), 
for any $g\in G\setminus \{e\}$. 

We can rephrase the statement of  Conjecture \ref{main conjecture s} in terms of principal bundles as follows. 
For a $G$-action $\alpha$ on $A$, we set 
$$\cP_\alpha=(EG\times \Aut(A))/G,$$
with a $G$-action on $EG\times \Aut(A)$ given by $g\cdot(x,\gamma)=(g\cdot x,\alpha_g\circ \gamma)$. 
Then $\cP_\alpha$ with the projection $[(x,\gamma)]\mapsto [x]\in BG$ is a principal $\Aut(A)$-bundle over $BG$, 
and it is identified with the pullback of the universal $\Aut(A)$-bundle $E\Aut(A)$ over $B\Aut(A)$ by 
the classifying map $B\alpha:BG\to B\Aut(A)$.  
The map $\cB$ being injective means that two outer $G$-actions $\alpha$ and $\beta$ 
are cocycle conjugate if and only if the two principal $\Aut(A)$-bundles $\cP_{\alpha}$ and $\cP_{\beta}$ are isomorphic. 
This reduces the problem of classifying group actions, which usually has difficulty of an analytic nature, to a purely topological problem. 
On the other hand $\cB$ being surjective means that every principal $\Aut(A)$-bundle over $BG$ arises from 
an outer $G$-action on $A$, and hence it shows richness of group actions. 

When $BG$ is a finite $CW$-complex, we can further rephrase the above statement in terms of locally trivial 
continuous fields of $C^*$-algebras. 
We denote by $M_\alpha$ the set of continuous sections $\Gamma(\cP_{\alpha}\times_{\Aut(A)}A)$ of the associated 
$A$-bundle $\cP_{\alpha}\times_{\Aut(A)}A$. 
Then $M_\alpha$ is a $C^*$-algebra, which we call the \textit{generalized mapping torus for $\alpha$} as it is 
the mapping torus when $G$ is the integer group $\Z$. 
Then every information of $\cP_\alpha$ is encoded in $M_\alpha$ as a $C(BG)$-algebra, 
or equivalently as a locally trivial continuous field of $A$ over $BG$.

After the first manuscript of this paper was posted on arXiv, Ralf Meyer \cite[Theorem 3.10]{Me19} obtained the following 
result. 

\begin{theorem}\label{Meyer} The map $\cB$ in Conjecture \ref{main conjecture s} is surjective for every torsion-free amenable 
second countable locally compact group $G$. 
\end{theorem}

Meyer actually showed, by using Baum-Connes conjecture techniques, 
that for $G$ as above the category $KK^G$ is equivalent to the subcategory of $KK^{BG}$ whose objects are locally trivial, 
and the equivalence restricts to the subcategories of objects whose underlying $C^*$-algebras are nuclear. 
The reader is referred to \cite{MN06} for the definition of the two categories $KK^G$ and $KK^{BG}$. 
The equivalence takes $(A,\alpha)\in KK^G$ to the $C^*$-algebra of continuous sections vanishing at infinity for the associated bundle 
$\cP_\alpha \times_{\Aut(A)}A$. 
On one hand, Kirchberg's classification result \cite[Folgerung 4.3]{Kir00} applied to a locally compact model of $BG$ 
(see \cite[Lemma 4.1]{KS03} for the existence of such a model) shows that two $KK^{BG}$-equivalent locally trivial continuous fields of 
a stable Kirchberg algebra $A$ over $BG$ are isomorphic, and hence their $KK^{BG}$-equivalence class gives rise to 
a unique isomorphism class of principal $\Aut(A)$-bundles over $BG$. 
On the other hand, Meyer \cite[Theorem 2.1]{Me19} showed that every object in $KK^G$ with a nuclear underlying $C^*$-algebra 
has a $KK^G$-equivalent object  with a stable Kirchberg algebra and an outer $G$-action. 
Thus for every principal $\Aut(A)$-bundle $\cP$ over $BG$, there exists a corresponding 
object $(A,\alpha)$ in $KK^G$ with outer $\alpha$ satisfying $\cP_\alpha\cong \cP$, which shows 
the surjectivity of $\cB$. 
In terms of equivariant KK-theory, now Conjecture \ref{main conjecture s} is rephrased that 
$KK^G$-equivalence implies cocycle conjugacy. 

It is not clear from the definition that the map $\cB$ in Conjecture \ref{main conjecture s} is well defined, 
which actually follows from Meyer's result because cocycle conjugacy always implies $KK^G$-equivalence. 
For this, we give a more elementary argument in Lemma \ref{de-stabilization} when $BG$ is a finite CW-complex.

Now we state a conjecture in the unital case, whose original version first appeared in \cite[Conjecture 1]{I10}.
For a Hilbert space $H$, we denote by $\K(H)$ the set of compact operators on $H$, and denote $\K=\K(\ell^2)$ 
for simplicity. 
The \textit{stabilization} $(\sA,\salpha)$ of a $C^*$-dynamical system $(A,\alpha)$ is defined by 
$$(A\otimes \K(\ell^2(G))\otimes \K,\alpha\otimes \Ad \rho\otimes \id_{\K}),$$
where $\rho$ is the right regular representation of $G$. 
Choosing a base point $*\in EG$, we regard $\cP_{\salpha}$ as a based space with a base point $[(*,\id)]$.

\begin{conjecture}\label{main conjecture u} Let $G$ be a torsion-free countable amenable discrete group, and let 
$A$ be a unital Kirchberg algebra. 
For two outer actions $\alpha$ and $\beta$ of $G$ on $A$, the following are equivalent: 
\begin{itemize}
\item[$(1)$] The two actions $\alpha$ and $\beta$ are $KK$-trivially cocycle conjugate. 
\item[$(2)$] There exists a base-point-preserving isomorphism between $\cP_{\salpha}$ and $\cP_{\sbeta}$.  
\end{itemize}
\end{conjecture}

In this paper, we concentrate on the case where $G$ is a poly-$\Z$ group. 

\begin{definition}\label{poly-Z} A discrete group $G$ is said to be \textit{poly-$\Z$} if there exists a subnormal series 
$$\{e\}=G_0\leq G_1\leq G_2\leq \cdots \leq G_n=G,$$
such that $G_i/G_{i-1}\cong \Z$ for any $1\leq i\leq n$. 
\end{definition}

The number $n$ in the above definition is called the Hirsch length of $G$ and denoted by $h(G)$. 
It does not depend on the choice of the subnormal series as above, and coincides with the cohomological dimension of $G$. 
A typical example of a poly-$\Z$ group is $\Z^n$, and the quotient map $\R^n\to \R^n/\Z^n=\T^n$ is a model 
of its universal bundle. 
More generally, every cocompact lattice of a simply connected solvable Lie group is poly-$\Z$. 

The following are our main results in this paper. 

\begin{theorem}\label{main theorem s} Conjecture \ref{main conjecture s} is true for every poly-$\Z$ group $G$ and every 
stable Kirchberg algebra $A$. 
\end{theorem}

\begin{theorem}\label{main theorem u} Conjecture \ref{main conjecture u} is true for every poly-$\Z$ group $G$ and every unital 
Kirchberg algebra $A$. 
\end{theorem}

We illustrate the significance of our results in the classification of group actions taking examples. 
Let $D$ be a strongly self-absorbing Kirchberg algebra (e.g. the Cuntz algebra $\cO_\infty$). 
Combining Theorem \ref{main theorem s} with Dadarlat-Pennig's generalized Dixmier-Douady theory \cite[Theorem 3.8]{DP-I}, 
we can see that the set of cocycle conjugacy classes of outer actions of a poly-$\Z$ group $G$ on $D\otimes \K$ is 
in one-to-one correspondence with the first group $E_D^1(BG)$ of their generalized cohomology theory associated with $D$. 
In particular, the cocycle conjugacy classes of outer $\Z^n$-actions on $D\otimes \K$ is in one-to-one correspondence with 
$$\Hom(\Z^n,KK(D,D)^{-1})\oplus \bigoplus_{k=1}^{[\frac{n-1}{2}]} H^{2k+1}(\T^n,K_0(D)),$$
where $[\frac{n-1}{2}]$ is the integer part of $\frac{n-1}{2}$. 
This is in a sharp contrast with the fact that there is a unique cocycle conjugacy class of outer $G$-actions on $D$ 
(see \cite[Theorem 4.4]{IMI}, and see \cite[Theorem B]{Sz1807arXiv} for a more general result). 
Theorem \ref{main theorem u} provides a better and conceptual understanding of this uniqueness result too 
because we know that $\Aut(D)$ is contractible (see \cite[Theorem 2.3]{DP15}).

It is a natural question to ask whether $\cP_{\salpha}$ in Theorem \ref{main theorem u} can be replaced by $\cP_{\alpha}$. 
Although we do not know the answer in the general case, we have a notable special case where it is affirmative. 
Using Theorem \ref{main theorem u} and Dadarlat's result \cite[Theorem 1.1]{D12} (see also \cite[Theorem 4.12]{IS}), we can show the following. 

\begin{theorem}\label{Cunzt algebra case} Let $m$ be a natural number, and let $\cO_{m+1}$ be the Cuntz algebra. 
Let $G$ be a poly-$\Z$ group whose cohomology group $H^*(BG)$ has no $m$-torsion. 
Then the map 
$$\cB:\mathcal{OA}(G,\cO_{m+1})/\textrm{cocycle conjugacy}\ni [\alpha]\mapsto [B\alpha]\in [BG,B\Aut(\cO_{m+1})]$$
is a bijection. 
In particular, the number of the cocycle conjugacy classes of outer $G$-actions on $\cO_{m+1}$ is 
$|\tilde{K}^0(BG)\otimes \Z_m|=|\tilde{H}^{\mathrm{even}}(BG;\Z_m)|$. 
\end{theorem}

For $G=\Z^n$, the cohomology group $H^*(B\Z^n)=H^*(\T^n)$ is torsion free, and we can apply the above theorem to $\Z^n$. 
Thus we see that there exist exactly $m^l$ cocycle conjugacy classes of outer $\Z^{n}$-actions on $\cO_{m+1}$, where
$$l=\sum_{k=1}^{[\frac{n}{2}]}\binom{n}{2k}.
$$
  
Now we discuss technical aspects of our paper. 
The notion of homotopy fixed points plays a crucial role throughout the paper. 

\begin{definition} Let $G$ be a topological group and let $Y$ be a $G$-space. 
A $G$-equivariant continuous map $f:EG\to Y$ is said to be a \textit{homotopy fixed point} of $Y$. 
\end{definition}

For $Y$ as above, we let $G$ act on $EG\times Y$ diagonally. 
Then $$p:EG\times _G Y:=(EG\times Y)/G\ni [(x,y)]\mapsto [x] \in BG$$ is a fiber bundle over $BG$ with fiber $Y$, called 
the Borel construction (or homotopy quotient).  
The set of continuous sections $s:BG\to EG\times_G Y$ is in one-to-one 
correspondence with the set of homotopy fixed points $f$ through the relation $s([x])=[(x,f(x))]$. 
We often identify a homotopy fixed point with the corresponding section throughout the paper. 

In various situations, our arguments in this paper can be interpreted as attempts to deduce the existence of a genuine fixed point 
from that of a homotopy fixed point, a formally weaker condition. 
For example, from the definition of the principal bundle $\cP_{\salpha}$, we see that the set of isomorphisms from $\cP_{\salpha}$ to $\cP_{\sbeta}$ 
is in one-to-one correspondence with the set of homotopy fixed points of $\Aut(\sA)$ with a $G$-action 
$g\cdot\gamma=\sbeta_g\circ \gamma\circ {\salpha_g}^{-1}$, 
and the correspondence is given as follows; if $\Phi$ is a homotopy fixed point, the corresponding isomorphism is 
$$[(x,\gamma)]\mapsto [(x,\Phi(x)\circ \gamma)]. $$
Note that this $G$-action has a genuine fixed point if and only if the two actions $\salpha$ and $\sbeta$ are conjugate in $\Aut(\sA)$, 
which actually follows from cocycle conjugacy of $\alpha$ and $\beta$.

To prove our main results by induction of the Hirsch length, we need to show a stronger statement. 
Instead of just constructing a fixed point from a given homotopy fixed point $\Phi$, 
we deform one of the actions to the other and simultaneously deform the homotopy fixed point $\Phi$ to a trivial one. 
There is a small price to pay; we need to amplify the actions by taking tensor product with an outer action $\mu$ on the Cuntz algebra $\cO_\infty$. 
But this is not a problem at all because of our absorption theorem \cite[Theorem 4.13]{IMI} (see also \cite[Corollary 3.7 ]{Sz18CMP}). 
The following is our main technical theorem. 
We denote by $\Aut(\sA)_0$ the path component of $\id$. 

\begin{theorem}[Deformation Theorem]\label{deformation}
Let $G$ be a poly-$\Z$ group, and let $\alpha$ and $\beta$ be outer $G$-actions on a 
unital Kirchberg algebra $A$. 
Assume that there exists a continuous map $\Phi:EG\rightarrow \Aut(\sA)_0$ satisfying 
$\Phi(g\cdot x)=\sbeta_g\circ \Phi(x)\circ{\salpha_g}^{-1}$ for any 
$x\in EG$ and $g\in G$. 
Then there exist a $\mu\otimes \alpha$-cocycle $\{c_g\}_{g\in G}$ in $U(\cO_\infty\otimes A)$, 
a continuous path $\{\gamma_t\}_{t\in I}$ in $\Aut(\cO_\infty\otimes A)$, and 
a continuous map $\tilde{\Phi}:EG\times [0,1]\rightarrow \Aut(\cO_\infty\otimes \sA)_0$ satisfying 
\begin{itemize}
\item $\gamma_0=\id$ and $\gamma_1\circ (\mu_g\otimes \beta_g)\circ \gamma_1^{-1}
=\Ad c_g\circ (\mu_g\otimes\alpha_g)$ for any $g\in G$, 
\item for any $x\in EG$, $t\in [0,1]$, and $g\in G$,
$$\tPhi(g\cdot x,t)=(\gamma_t\otimes \id)\circ(\mu_g\otimes \sbeta_g)\circ (\gamma_t^{-1}\otimes \id)
\circ \tPhi(x,t)\circ (\mu_g\otimes\salpha_g)^{-1},$$ 
\item $\tPhi(x,0)=\id \otimes \Phi(x)$ and $\tPhi(x,1)=\Ad W$ for any $x\in EG$,
where  
$$W=\sum_{g\in G}c_{g^{-1}}\otimes E_{g,g}\otimes 1,$$
and $\{E_{g,h}\}_{g,h\in G}$ is the canonical system of matrix units in $\K(\ell^2(G))$.  
\end{itemize}
\end{theorem}

Injectivity of $\cB$ in Theorem \ref{main theorem s} follows from Theorem \ref{main theorem u}, 
which in turn follows from Theorem \ref{deformation} (see Section 2). 
We believe that Theorem \ref{deformation} applied to the case $\beta_g=\alpha_g$ will also play an important role 
in the future for understanding topological properties of the gauge group of $\cP_{\salpha}$. 

This paper is organized as follows. 
In Section 2, we prove direct consequences of Theorem \ref{deformation} and recall necessary results for the proof of Theorem \ref{deformation}. 
In section 3, we prove two cohomology vanishing theorems, Theorem \ref{vanishing1} for 1-cocycles of group actions and Theorem \ref{cocycle actions} 
for cocycle actions, in terms of homotopy fixed points. 
They are of interest in their own right because this is the first time to precisely identify topological obstructions 
of the problems to the best of knowledge of the authors. 

Section 4 is entirely devoted to the proof of Theorem \ref{deformation}. 
Roughly speaking the main idea of the proof is as follows. 
To prove the statement by induction of the Hirsch length $n=h(G)$, we decompose $G$ as a semidirect product $G=G_{n-1}\rtimes \Z$. 
The induction hypothesis shows the existence of a section of a certain fiber bundle over $BG_{n-1}$ with fiber the loop group 
$\Omega\Aut(\sA)$ of $\Aut(\sA)$. 
Thanks to our weak homotopy equivalence type theorem in \cite{I10} (see Theorem \ref{weak homotopy} below), we can replace 
the fiber $\Omega\Aut(\sA)$ with the unitary group $U(A_\flat)$ of the continuous asymptotic centralizer $A_\flat$ of $A$, and we apply 
Theorem \ref{vanishing1}. 

For the sake of self-contained presentation, we give a proof of the surjectivity of $\cB$ in Theorem \ref{main theorem s} 
by using Theorem \ref{deformation} in section \ref{sDRT}. 
In section 6, we identify a cohomological obstruction for two actions to be $KK$-trivially cocycle conjugate obtained in 
our previous work \cite[Section 7]{IMI} with the primary obstruction for the existence of a continuous section for a certain fiber bundle. 
In section 7, we prove Theorem \ref{Cunzt algebra case}, for which Theorem \ref{cocycle actions} is used in an essential way.

\textbf{Acknowledgment.} 
We would like to thank Marius Dadarlat, Andr\'e Henriques, Daisuke Kishimoto, Ralf Meyer, Narutaka Ozawa, Ulrich Pennig, and Dai Tamaki 
for useful discussions and comments. 
Masaki Izumi would like to thank Isaac Newton Institute for Mathematical Sciences for its hospitality.


\section{Preliminaries}
In this section, we prove direct consequences of Theorem \ref{deformation}, and collect necessary facts 
for its proof. 

For a unital $C^*$-algebra $A$, we denote by $U(A)$ the unitary group of $A$. 
If $J$ is a closed ideal of $A$, we set $U(J)=\{u\in U(A);\; u-1\in J\}$. 
Recall that we denote the stabilization $A\otimes \K$ of $A$ by $\sA$. 
We often use the fact that the unitary group $U(M(\sA))$ of the stable multiplier algebra of $A$ 
is contractible in the norm topology (see \cite[Theorem 16.8]{W-O}) as well as in the strict topology (see \cite[Lemma 4.72]{RW}). 
In consequence, we have $\pi_1(\mathrm{Inn}(\sA))=\{0\}$, where $\mathrm{Inn}(\sA)$ is the group of inner automorphisms of $\sA$. 
In fact, the exact sequence 
$$0\to \T\to U(M(A^s))\to \mathrm{Inn}(\sA)\to 0,$$
implies that $\mathrm{Inn}(\sA)$ is an Eilenberg-MacLane space of type $K(\Z,2)$. 

Recall that the multiplier algebra of $\K(\ell^2(G))$ is $B(\ell^2(G))$, whose unitary group 
is the unitary group $U(\ell^2(G))$ of the Hilbert space $\ell^2(G)$.

\begin{lemma}\label{hfp for rho} Let $G$ be a discrete group with $BG$ homotopy equivalent to a CW-complex. 
Then there exists a norm continuous map $f:EG\to U(\ell^2(G))$ satisfying $f(g\cdot x)=\rho_gf(x)$ 
for every $x\in EG$ and $g\in G$. 
\end{lemma}

\begin{proof}
We consider a principal $U(\ell^2(G))$-bundle 
$EG\times_G U(\ell^2(G))\to BG$ with a $G$-action on $U(\ell^2(G))$ given by $g\cdot u=\rho_g u$. 
Since $U(\ell^2(G))$ is contractible and $BG$ is homotopy equivalent to a CW-complex, the principal bundle has a continuous section, 
and there exists a homotopy fixed point $f:EG\to U(\ell^2(G))$. 
\end{proof}

\begin{remark}
When we have two isomorphic $C^*$-algebras $A_1$ and $A_2$, we identify two groups $\Aut(A_1)$ and $\Aut(A_2)$ only through 
isomorphisms between the underlying $C^*$-algebras $A_1$ and $A_2$, and we never use abstract isomorphisms between the two 
topological groups $\Aut(A_1)$ and $\Aut(A_2)$ even when they exist. 
With this convention, it makes sense to say that a principal $\Aut(A_1)$-bundle $\cP_1$ over a space $X$  
and a principal $\Aut(A_2)$-bundle $\cP_2$ over $X$ are isomorphic. 
When $X$ is compact, this condition is equivalent to 
$$\Gamma(\cP_1\times_{\Aut(A_1)}A_1)\cong \Gamma(\cP_2\times_{\Aut(A_2)}A_2),$$
as $C(X)$-algebras. 
\end{remark}

For two $C^*$-algebras $A$ and $B$, we regard $\Aut(A)$ as subgroup of $\Aut(A\otimes B)$ via the inclusion map $\theta \mapsto \theta\otimes \id_B$, 
where $\otimes$ always means the minimal tensor product throughout the paper.  
With this convention $\cP\times_{\Aut(A)} \Aut(A\otimes B)$ is a principal $\Aut(A\otimes B)$-bundle if $\cP$ is a principal $\Aut(A)$-bundle. 

\begin{lemma}\label{tensoring} Let $A$ be a Kirchberg algebra, and let $\cP$ be a principal $\Aut(A)$-bundle over a finite CW-complex $X$. 
\begin{itemize}
\item[$(1)$] Two principal bundles $\cP$ and $\cP\times_{\Aut(A)}\Aut(A\otimes \cO_\infty)$ are isomorphic. 
\item[$(2)$] If $A$ is stable, two principal bundles $\cP$ and $\cP\times_{\Aut(A)}\Aut(A\otimes \K)$ are isomorphic. 
\end{itemize}
\end{lemma}

\begin{proof} It suffices to show the corresponding $C(X)$-algebras are isomorphic. 

(1) We denote $\cP'=\cP\times_{\Aut(A)}\Aut(A\otimes \cO_\infty)$ for simplicity. 
Then we have 
$$\cP'\times_{\Aut(A\otimes \cO_\infty)}A\otimes \cO_\infty=\cP\times _{\Aut(A)}A\otimes \cO_\infty,$$ 
and an obvious inclusion 
$$\Gamma(\cP\times_{\Aut(A)}A)\otimes \cO_\infty\subset \Gamma(\cP'\times_{\Aut(A\otimes \cO_\infty)}A\otimes \cO_\infty).$$
We claim that the two $C^*$-algebras actually coincide. 
Indeed, choosing an appropriate finite open covering $\{U_i\}_{i=1}^n$ of $X$ and local trivialization of $\cP$, we have an embedding 
$$\Gamma(\cP'\times_{\Aut(A\otimes \cO_\infty)}A\otimes \cO_\infty)\subset \bigoplus_{i=1}^nC(\overline{U_i})\otimes A\otimes \cO_\infty,$$
and the claim follows from \cite[Theorem 12.4.4]{BO} as $\cO_\infty$ is nuclear. 

Thanks to \cite[Theorem 1.5]{BK04} (see also \cite[Theorem 4.6]{HRW07}), we have  
$$\Gamma(\cP\times_{\Aut(A)}A)\otimes \cO_\infty\cong \Gamma(\cP\times_{\Aut(A)}A),$$
and so  
$$\Gamma(\cP'\times_{\Aut(A\otimes \cO_\infty)}A\otimes \cO_\infty)\cong \Gamma(\cP\times_{\Aut(A)}A),$$
as $C(X)$-algebras. 

(2) The statement follows from \cite[Proposition 3.4]{HRW07}). 
\end{proof}

\begin{lemma}\label{de-stabilization} Let $A$ be a stable Kirchberg algebra, and let $G$ be a discrete group with $BG$ homotopy equivalent to 
a finite CW complex.  
For a $G$-action $\alpha$ on $A$, the two principal bundles $\cP_\alpha$ and $\cP_{\salpha}$ are isomorphic. 
\end{lemma}

\begin{proof} Let $f:EG\to U(\ell^2(G))$ be as in Lemma \ref{hfp for rho}, and let $\Phi(x)=\id_A\otimes \Ad f(x)\otimes \id_\K$. 
Then $\Phi$ satisfies $\Phi(g\cdot x)=\salpha_g\circ \Phi(x)\circ (\alpha_g\otimes \id_{\K(\ell^2(G))\otimes \K})^{-1}$ and 
$\cP_{\salpha}$ and $\cP_{\alpha\otimes \id_{\K(\ell^2)\otimes \K}}$ are isomorphic. 
On the other hand, we have 
$$\cP_{\alpha\otimes \id_{\K(\ell^2)\otimes \K}}=\cP_\alpha\times_{\Aut(A)}\Aut(A\otimes \K(\ell^2(G))\otimes \K),$$       
and the statement follows from Lemma \ref{tensoring}. 
\end{proof}

Let $\alpha$ be an action of a discrete group $G$ on a $C^*$-algebra $A$, and let $\{u_g\}_{g\in G}$ be an $\alpha$-cocycle. 
We set 
$$W=\sum_{g\in G}u_{g^{-1}}\otimes E_{g,g}\otimes 1,$$
where convergence is in the strict topology. 
Then $W$ is a unitary in $M(\sA)$ satisfying $W\salpha_g(W^*)=u_g\otimes 1$, and we get 
$$(\alpha^u)^{\mathrm{s}}_g=\Ad W\circ \salpha_g\circ \Ad W^{-1}.$$
Thus if $\alpha$ and $\beta$ are cocycle conjugate, two principal bundles $\cP_{\salpha}$ and $\cP_{\sbeta}$ are isomorphic.

We often use the following solution of a homotopy lifting extension problem 
(see \cite[Theorem 7.16]{W78}). 

\begin{theorem}\label{lifting extension} Let $p:E\to B$ be a fibration, and let $(Y,A)$ be an NDR pair. 
Then for any continuous maps $f:I\times Y\to B$ and $g:(I\times A)\cup (0\times Y)\to E$ 
with $p\circ g=f|_{(I\times A)\cup(0\times Y)}$, there exists a continuous map 
$h:I\times Y\to E$ that is an extension of $g$ and a lifting of $f$.  
In particular, the statement holds if $p:E\to B$ is a fiber bundle over a paracompact space $B$, and 
$A$ is a subcomplex of a CW-complex $Y$. 
\end{theorem}

The following statement is probably well-known. 
Since we cannot find an appropriate reference for it, we will prove it in Appendix. 

\begin{theorem}\label{model} Let $G$ be a poly-$\Z$ group with $h(G)=n$. 
With appropriate choices of a subnormal series $\{G_i\}_{i=0}^n$ and 
$\xi_i\in G_i$ with $G_i=G_{i-1}\rtimes \langle\xi_i\rangle$  for $1\leq i\leq n$, 
there exists a free cocompact action of $G$ on $\R^n$ satisfying 
$$\xi_i\cdot(x,x_{i},y)=(\xi_i\cdot x,x_i+1,y)$$
for any $x\in \R^{i-1}$, $x_i\in \R$, and $y\in \R^{n-i}$ with convention $\R^0=\emptyset$. 
In particular, we can adopt $\R^n$ with the $G$-action as above for $EG$ so that $EG_i$ is identified 
with $$\{(x,0_{n-i})\in EG;\; x\in \R^i\}.$$
\end{theorem}

In what follows, we assume that $G$ is a poly-$\Z$ group unless otherwise stated. 
We denote $n=h(G)$, and fix $\{G_i\}_{i=0}^n$, $\{\xi_i\}_{i=1}^n$, 
and a $G$-action on $\R^n$ satisfying the statement of the above theorem. 
Since we will often use induction arguments on the Hirsch length, we denote $N=G_{n-1}$ and $\xi=\xi_{n}$ for simplicity. 
We often use variables $z=(x,y)\in \R^n=EG$ with $x\in \R^{n-1}=EN$ and $y\in \R$.  
For $z\in EG$, we denote by $[z]$ its class in $BG$. 
We identify $BN=EN/N$ with the submanifold $\{[(x,0)]\in BG;\;x\in EN\}$ of $BG$. 
We choose $0_n$ as a base point of $EG$. 

Although the statements of our main results do not depend on the model of $EG$, we need to fix a particular model for their proofs 
because we do not have a topological interpretation of Theorem \ref{vanishing2}, and we do not know if the statement holds independent 
of the model of $EG$. 

Assume that we are given two $G$-actions $\alpha$ and $\beta$ on a $C^*$-algebra $A$. 
Recall that there is a one-to-one correspondence between the set of isomorphisms from $\cP_{\salpha}$ to $\cP_{\sbeta}$ and 
the set of homotopy fixed points of $\Aut(\sA)$ with a $G$-action $g\cdot\gamma=\sbeta_g\circ \gamma\circ {\salpha_g}^{-1}$. 
In particular, an isomorphism corresponding to a homotopy fixed point $\Phi:EG\to \Aut(\sA)$ preserves the base points if and only if 
$\Phi(0_n)=\id_{\sA}$. 
This is possible only if the image of $\Phi$ is in $\Aut(\sA)_0$, and in consequence $\sbeta_g\circ {\salpha_g}^{-1}\in \Aut(\sA)_0$. 
Note that if $A$ is a Kirchberg algebra, the last condition is equivalent to  $KK(\alpha_g)=KK(\beta_g)$.

\begin{lemma}\label{reduction to hfp} Let $G$ be a poly-$\Z$ group, and let $\alpha$ and $\beta$ be $G$-actions on a unital Kirchberg algebra $A$. 
Then the following two conditions are equivalent. 
\begin{itemize}
\item[$(1)$] There exists  a base point preserving isomorphism between the two principal $\Aut(\sA)$-bundles 
$\cP_{\salpha}$ and $\cP_{\sbeta}$.  
\item[$(2)$] $KK(\alpha_g)=KK(\beta_g)$ for each $g\in G$, 
and the $G$-action $g\cdot \gamma=\sbeta_g\circ \gamma\circ {\salpha_g}^{-1}$ on $\Aut(\sA)_0$ has a homotopy fixed point. 
\end{itemize}
If $\alpha$ and $\beta$ are $KK$-trivially cocycle conjugate, the above equivalent conditions hold. 
\end{lemma}

\begin{proof} We have already seen that (1) implies (2). 
On the other hand, if $\Phi$ is a homotopy fixed point for the $G$-action on $\Aut(\sA)_0$ in (2), 
we can continuously deform $\Phi$ so that it satisfies 
$\Phi(0_n)=\id$ thanks to Theorem \ref{lifting extension}. 
Therefore (2) implies (1). 

Assume that $\alpha$ and $\beta$ are $KK$-trivially cocycle conjugate now. 
Then there exist $\gamma\in \Aut(A)$ with $KK(\gamma)=KK(\id_A)$ and an $\alpha$-cocycle $\{u_g\}_{g\in G}$ in $U(A)$ satisfying 
$\gamma\circ \beta_g\circ \gamma^{-1}=\alpha^u_g$. 
Let $W$ be a unitary in $M(\sA)$ as before. 
Then
$$\Ad W\circ\salpha_g\circ \Ad W^{-1}=(\gamma\otimes \id)\circ \sbeta_g\circ(\gamma\otimes \id)^{-1}.$$    
Note that we have $(\gamma\otimes \id)^{-1}\circ \Ad W\in \Aut(\sA)_0$, and $\Phi(x)=(\gamma\otimes \id)^{-1}\circ \Ad W$ 
is a desired (homotopy) fixed point in (2).  
\end{proof}

\begin{proof}[Proof of Theorem \ref{main theorem u} assuming Theorem \ref{deformation}] 
Implication from (1) to (2) was already shown in Lemma \ref{reduction to hfp}. 
On the other hand,  Lemma \ref{reduction to hfp} and Theorem \ref{deformation} imply that 
the KK-trivial cocycle conjugacy of $\mu\otimes \alpha$ and $\mu\otimes \beta$ follows from (2), 
and (1) follows too thanks to \cite[Theorem 4.13]{IMI}. 
\end{proof}

\begin{proof}[Proof of the injectivity of $\cB$ in Theorem \ref{main theorem s} assuming Theorem \ref{main theorem u}] 
Assume that $\alpha$ and $\beta$ are outer $G$-actions of $G$ on a stable Kirchberg algebra $A$ satisfying 
$\cP_{\alpha}\cong \cP_{\beta}$. 
We may assume $A=B\otimes \K$ with $B$ a unital Kirchberg algebra in the Cuntz standard form, i.e. $[1_B]_0=0$ in $K_0(B)$. 
We choose a system of matrix units $\{E_{i,j}\}_{i,j\in \N}$ in $\K$. 
Then thanks to \cite[Theorem 4.17]{IMI}, we have a family of partial isometries $\{u_g\}_{g\in G}$ in $A$ satisfying 
$u_g^*u_g=\alpha_g(1\otimes E_{11})$, $u_gu_g^*=1\otimes E_{11}$, and $u_g\alpha_g(u_h)=u_{gh}$ for any $g,h\in G$.
Let 
$$V_g=\sum_{i=1}^\infty (1\otimes E_{i,1})u_g\alpha_g(1\otimes E_{1,i}),$$
where convergence is in the strict topology. 
Then $\{V_g\}_{g\in G}$ is an $\alpha$-cocycle satisfying $\alpha_g^V(1\otimes E_{ij})=1\otimes E_{ij}$. 
Therefore letting $\alpha'$ be the restriction of $\alpha^V$ to $(1\otimes E_{11})A(1\otimes E_{11})$ and identifying $(1\otimes E_{11})A(1\otimes E_{11})$ with $B$, 
we get $\alpha^V=\alpha'\otimes \id_\K$. 
Thus to prove that $\alpha$ and $\beta$ are cocycle conjugate, it suffices to assume that 
$\alpha=\alpha'\otimes \id_\K$ and $\beta=\beta'\otimes \id_\K$ from the beginning.  

Since $\cP_{\salpha}\cong \cP_{\sbeta}$, there exists a continuous map $\Phi:EG\to \Aut(\sA)$ satisfying 
$\Phi(g\cdot x)=\sbeta_g\circ \Phi(x)\circ {\salpha_g}^{-1}$ for any $x\in EG$ and $g\in G$. 
We choose $\theta\in \Aut(B)$ with $KK(\theta)=KK(\Phi(0_n))$. 
Replacing $\alpha'_g$ with $\theta\circ \alpha'_g\circ \theta^{-1}$ and $\Phi(x)$ with $\Phi(x)\circ (\theta^{-1}\otimes \id)$, 
we may assume that the image of $\Phi$ is in $\Aut(\sA)_0$. 
Now Lemma \ref{reduction to hfp} and Theorem \ref{main theorem u} imply that $\alpha'$ and $\beta'$ are cocycle conjugate, 
and so are $\alpha$ and $\beta$. 
\end{proof}

\begin{remark}\label{perturbation} The above proof shows that for any outer $G$-action $\alpha$ on a stable Kirchberg algebra, 
there exist a factorization $A=B\otimes \K$ with $B$ a unital Kirchberg algebra in the Cuntz standard form 
and an outer $G$-action $\alpha'$ on $B$ such that $\alpha'\otimes \id$ is a cocycle perturbation of $\alpha$. 
\end{remark}

Surjectivity of $\cB$ in Theorem \ref{main theorem s} will be shown in Section \ref{sDRT}.

In the rest of this section, we collect necessary ingredients for the proof of Theorem \ref{deformation}. 
For a group $G$ and $g,h\in G$, we denote $g^h=h^{-1}gh$. 
For a group automorphism $\theta\in \Aut(G)$, we use the notation $g^{\theta}$ for $g\in G$ acted by $\theta$.

Our standard reference for algebraic topology is \cite{W78}. 
We use notation $I=[0,1]$, $D^n=\{x\in \R^n;\;|x|\leq 1\}$, and $S^n=\partial D^{n+1}$. 
For two topological spaces $X$ and $Y$, we denote by $[X,Y]$ the set of homotopy equivalence classes of continuous maps from $X$ to $Y$, 
and by $[X,Y]_0$ its based counterpart if $X$ and $Y$ are based spaces. 
For a topological group $\Gamma$, we denote by $\Gamma_0$ the path connected component of the neutral element $e\in \Gamma$. 
We always regard $\Gamma$ as a based space with a base point $e$.  
We denote by $\Omega\Gamma$ the based loop space of $\Gamma$.

\begin{lemma}\label{NDR} Let the notation be as above.  
\begin{itemize}
\item[$(1)$] $(BG,BN)$ is an NDR pair.
\item[$(2)$] $(BG, BN\cup [\{0_{n-1}\}\times I)])$ is an NDR pair. 
\end{itemize}
\end{lemma}

\begin{proof} Let $T_\xi$ be a homeomorphism of $BN=EN/N$ induced by $\xi$. 
Then $BG$ is identified with the mapping torus 
$$\frac{BN\times I}{(x,0)\sim (T_\xi(x),1)},$$
and we denote by $q:BN\times I\to BG$ the quotient map. 
From the definition of an NDR pair, we see that if $X$ is a subset of $BG$ including $q(BN\times \{0\})$ such that 
$(BN\times I,q^{-1}(X))$ is an NDR pair, so is $(BG,X)$. 
Thus (1) holds as $(BN\times I, BN\times \partial I)$ is an NDR pair. 
From \cite[I.5.2]{W78} we see that 
$$(BN\times I, (BN\times \partial I)\cup (\{[0_{n-1}]\}\times I))$$
is an NDR pair too, which implies (2).  
\end{proof}

Thanks to \cite{AJ76}, there exists a triangulation of $BG$, and we can choose one with 
$[0_{n}]$ being a vertex of it. 
Then $EG$ has a free $G$-CW-complex structure coming from this triangulation, and we can choose a finite subcomplex 
$K$ of $EG$ satisfying $0_n\in K$ and $G\cdot K=EG$. 
By construction, if $e$ is an open cell in $EG$ of dimension $k$, the characteristic map $f:D^k\to \overline{e}$ is a 
homeomorphism, and  $[\overline{e}]\subset BG$ is a simplex.

We set 
$$A^\flat=C^b([0,1),A)/C_0([0,1),A),$$
and we often identify an element in $A^\flat$ with one of its representative if there is no possibility of confusion. 
Identifying $A$ with the subalgebra of $A^\flat$ consisting of constant functions, we set 
$A_\flat=A^\flat\cap A'$. 

We recall a weak homotopy equivalence type result from \cite[Corollary 2.9]{IMW}. 

\begin{theorem}\label{weak homotopy} Let $A$ be a unital Kirchberg algebra, and let $X$ be a finite CW-complex 
with a base point. 
Then there exists an isomorphism 
$$\Pi_{A,X}:[X,U(A_\flat)]_0\to [X,\Omega\Aut(\sA)]_0,$$
that is natural in $X$. 
\end{theorem}

We denote $\Pi_{A,i}=\Pi_{A,S^i}$ for simplicity. 
The map $\Pi_{A,X}$ is constructed as follows. 
Since $X$ is compact, for a given continuous map $u:X\to U(A_\flat)_0$, we can choose its representative 
$u(x)(y)$ that is a continuous map from $X\times [0,1)$ to $U(A)_0$. 
For this representative, we can choose a continuous map $U:X\times [0,1)\to U(M(\sA))$ with $U(x,y)=u(x)(y)\otimes 1$ 
for any $x\in X$, $y\in [1/2,1)$ and $U(x,0)=1$ for any $x\in X$. 
Then $\Pi_{A,X}([u])$ is given by the homotopy class of $\Ad U(x,y)$ regarded as a map from $X$ to $\Omega\Aut(A\otimes \K)$.  

Let $G$ be a discrete group. 
A cocycle action $(\alpha,u)$ on $A$ is a pair of a map $\alpha:G\to \Aut(A)$ and a map $u:G\times G\to U(A)$ satisfying 
$$\alpha_g\circ \alpha_h=\Ad u_{g,h}\circ \alpha_{gh},\quad \forall g,h\in G,$$
$$u_{g,h}u_{gh,k}=\alpha_g(u_{h,k})u_{g,hk},\quad \forall g,h,k\in G.$$
Two cocycle actions $(\alpha,u)$ and $(\alpha',u')$ are equivalent if there exists a map $v:G\to U(A)$ satisfying 
$$\alpha'_g=\Ad v_g\circ \alpha_g,\quad \forall g\in G,$$
$$u'_{g,h}=v_g\alpha_g(v_h)u_{g,h}v_{gh}^{-1},\quad \forall g,h\in G.$$

Let $A^\infty=\ell^\infty(\N,A)/c_0(\N,A)$ and let $A_\infty=A^\infty\cap A'$. 
Then $(\alpha,u)$ induces cocycle actions on $A^\flat$ and $A^\infty$, and actions on $A_\flat$ and $A_\infty$. 
We abuse the notation and we still denote these induced (cocycle) actions by $(\alpha,u)$ and $\alpha$ for simplicity. 

We recall the following definition from \cite[Definition 5.1]{IM2010}. 

\begin{definition}
Let $G$ be a countable discrete group and 
let $\mu:G\curvearrowright\mathcal{O}_\infty$ be an action. 
Let $(\alpha,u)$ be a cocycle action of $G$ 
on a unital (not-necessarily separable) $C^*$-algebra $A$. 
We say that $(\alpha,u)$ admits 
an approximately central embedding of $(\mathcal{O}_\infty,\mu)$ 
if for any separable subset $S\subset A$ there exists 
a unital homomorphism $\phi:\mathcal{O}_\infty\to A^\infty\cap S'$ 
such that $\phi\circ\mu_g=\alpha_g\circ\phi$ for any $g\in G$. 
We denote by $\AC(\mathcal{O}_\infty,\mu)$ 
the class of such cocycle actions. 
\end{definition}

For a poly-$\Z$ group $G$, the class $\AC(\mathcal{O}_\infty,\mu)$ does not depend 
on the choice of $\mu$ (see \cite[Lemma 5.2]{IMI}) as long as it is outer, and  
we choose and fix an outer action $\mu^G$ on $\cO_{\infty}$. 
We often suppress $G$ in $\mu^G$ if there is no possibility of confusion.

\begin{lemma}[{ \cite[Lemma 5.3]{IMI}}]\label{AC}
Let $G$ be a poly-$\Z$ group, and let $(\alpha,u):G\curvearrowright A$ be an outer cocycle action 
on a unital Kirchberg algebra $A$.  
Then, 
\begin{enumerate}
\item[$(1)$] $(\alpha,u):G\curvearrowright A$ is in $\AC(\mathcal{O}_\infty,\mu^G)$. 
\item[$(2)$] $(\alpha,u):G\curvearrowright A^\infty$ and 
$\alpha:G\curvearrowright A_\infty$ 
are in $\AC(\mathcal{O}_\infty,\mu^G)$. 
\item[$(3)$] $(\alpha,u):G\curvearrowright A^\flat$ and 
$\alpha:G\curvearrowright A_\flat$ 
are in $\AC(\mathcal{O}_\infty,\mu^G)$. 
\end{enumerate}
\end{lemma}

We recall the following classification theorem, which is a special case of 
\cite[Theorem 6.4 and Theorem 6.5]{IMI}. 

\begin{theorem}\label{KKtrivcc}
Let $G$ be a poly-$\Z$ group and 
let $\alpha$ and $\beta$ be outer actions of $G$ on a unital Kirchberg algebra $A$. 
Suppose that 
there exists a family $(x_g)_g$ of unitaries in $C^b([0,1),A)$ such that 
$$
\lim_{s\to 1}(\Ad x_g(s)\circ\alpha_g)(a)=\beta_g(a),
\quad\forall g\in G,\ \forall a\in A, $$
$$
\lim_{s\to 1}x_g(s)\alpha_g(x_h(s))x_{gh}(s)^*
=1, \quad\forall g,h\in G. $$
Then 
there exist a continuous map $w:[0,1)\to U(\mathcal{O}_\infty\otimes A)$, 
$\gamma\in\Aut(\mathcal{O}_\infty\otimes A)$ and 
a $\mu^G\otimes \alpha$-cocycle $(c_g)_g$in $U(\mathcal{O}_\infty \otimes A)$ satisfying 
$$\gamma\circ (\mu^G_g\otimes \beta_g)\circ \gamma^{-1}=\Ad c_g\circ (\mu^G_g\otimes \alpha_g),\quad \forall g\in G,$$
$$
\lim_{s\to 1}\Ad w(s)(a)=\gamma(a),
\quad\forall a\in \mathcal{O}_\infty\otimes A, $$
$$
\lim_{s\to 1}w(s)(1\otimes x_g(s))(\mu^G_g\otimes \alpha_g)(w(s))^*=c_g
\quad\forall g\in G.$$
Moreover, we may assume $w(0)=1$. 
\end{theorem}

\section{Cohomology vanishing via homotopy fixed points}
\subsection{First cohomology vanishing}
Let $\alpha$ be an action of a discrete group $G$ on a unital $C^*$-algebra $A$, and let $u=\{u_g\}_{g\in G}$ be an $\alpha$-cocycle 
in $U(A)_0$. 
We introduce a $G$-action on $U(A)_0$ by $g\cdot v=u_g\alpha_g(v)$. 
Note that this action has a fixed point if and only if the cocycle $\{u_g\}_{g\in G}$ is a coboundary 
within $U(A)_0$, i.e. there exists $v\in U(A)_0$ satisfying $u_g=v\alpha_g(v^{-1})$ for any $g\in G$. 
We would like to establish a similar criterion for $\{u_g\}_{g\in G}$ to be asymptotic coboundary in terms 
of a homotopy fixed point. 
If $F$ is a homotopy fixed point for this action, we simply say that $F$ is a homotopy fixed point for $u$. 
We denote by $E^u$ the Borel construction $EG\times_G U(A)_0\to BG$. 

\begin{remark} We can identify the class $\kappa^2(\alpha,u)\in H^2(G,K_0(A))$ introduced in \cite[Section 7]{IMI} with 
the primary obstruction in $H^2(BG,\pi_1(U(A)_0))$ for the existence of a continuous section of the fiber bundle $p:E^u\to BG$ 
in a similar way as in the proof of Theorem \ref{obstruction}. 
\end{remark}

\begin{lemma}\label{section} Let the notation be as above. 
Assume that the cocycle $\{u_g\}_{g\in G}$ is an asymptotic coboundary in $U(A)_0$, that is, there 
exists a continuous map $v:[0,1)\to U(A)_0$ satisfying 
$$\lim_{t\to 1}v(t)\alpha_g(v(t)^*)=u_g,\quad  \forall g\in G.$$ 
Then $p:E^u\to BG$ has a continuous section, and hence there exists a homotopy fixed point for $u$.  
\end{lemma}

\begin{proof} Let 
$$u(t)_g=\left\{
\begin{array}{ll}
v(t)\alpha_g(v(t)^*) , &\quad t\in [0,1) \\
u_g , &\quad t=1
\end{array}
\right..
$$
All we need is that $\{u(t)_g\}_{t\in [0,1],\; g\in G}$ is a continuous family of $\alpha$-cocycles that is 
a coboundary at $t=0$ and $u_g$ at $t=1$. 
We introduce a homotopy of fiber bundles 
$$p:E=(EG\times I\times U(A)_0)/G\to BG\times I$$ 
by $g\cdot(x,t,v)=(g\cdot x,t,u(t)_g\alpha_g(v))$. 
Since $BG\times \{0\}$ has a continuous section $s([x])=[(x,0,v(0))]$, the statement follows from the homotopy 
lifting property of the fiber bundle $p:E\to BG\times I$. 
\end{proof}

From now, we assume that $G$ is a poly-$\Z$ group, and hence $n=h(G)$ and $EG=\R^n$ from our convention. 
Embedding $G$ into $EG$ via $G\cdot 0_n\subset EG$, we can think of 
a homotopy fixed point $F:EG\to U(A)_0$ with $F(0_n)=1$ as a continuous extension 
of the cocycle $u_g$. 

The following theorem is a strengthened converse of Lemma \ref{section}. 

\begin{theorem} \label{vanishing1}
Let $\alpha$ be an action of a poly-$\Z$ group $G$ on a unital $C^*$-algebra $A$ in $\AC(\cO_\infty,\mu^G)$, 
and let $J$ be a closed ideal of $A$ that is globally preserved by $\alpha$.  
Assume that $u=\{u_g\}_{g\in G}\subset U(J)_0$ is an $\alpha$-cocycle, and $F:EG\rightarrow U(J)_0$ 
is a homotopy fixed point for $u$. 
Then there exist a continuous family $\{v(t)\}_{t\in [0,1)}$ of unitaries in $U(J)_0$ and 
a continuous map $\tF:EG\times I\rightarrow U(J)_0$ satisfying  
\begin{itemize}
\item $v(0)=1$ and $\lim_{t\to 1}v(t)\alpha_g(v(t)^*)=u_g$ for any $g\in G$,  
\item $\tF(gx,t)=v(t)^*u_g\alpha_g(v(t)\tF(x,t))$ for any $x\in EG$, $g\in G$, and $t\in I$, 
\item $\tF(x,0)=F(x)$ and $\tF(x,1)=1$ for any $x\in EG$. 
\end{itemize}
Moreover, if $A=J=B_\flat$ with a unital Kirchberg algebra $B$ and $\alpha$ comes from an outer action of $G$ on $B$, 
we can choose $v(t)$ so that it 
continuously extends to $[0,1]$. 
\end{theorem}

Since $v(t)$ in the statement does not allow the limit at $t=1$ in general, proof by induction of the 
Hirsch length of $G$ would be considerably complicated, and we take an alternative approach. 
We first prove a 2-cocycle version, for which induction works well.   

Let $(\alpha,u)$ be a cocycle action of $G$ on a unital $C^*$-algebra $A$ with $u(g,h)\in U(A)_0$. 
Assume that $u$ is a coboundary within $U(A)_0$, that is, there exist $w_g\in U(A)_0$ for $g\in G$ satisfying 
$u(g,h)=\alpha_g(w_h^{-1})w_g^{-1}w_{gh}$. 
We can extend $w_g$ to a continuous map $w:EG\to U(A)_0$ in the sense that we have $w(g\cdot 0_n)=w_g$, 
and set 
$$v(g,x)=\alpha_g(w(x)^{-1})w(g\cdot 0_n)^{-1}w(g\cdot x)$$ 
for $g\in G$ and $x\in EG$. 
Then $v:G\times EG\to U(A)_0$ is a continuous map satisfying $v(g,h\cdot 0_n)=u(g,h)$ and 
$$u(g,h)v(gh,x)=\alpha_g(v(h,x))v(g,h\cdot x).$$ 
The existence of the extension $v$ of $u$ as above is indeed a sufficient condition for $u$ to be 
a coboundary within $U(A)_0$ in a good situation. 

\begin{theorem}\label{vanishing2}
Let $(\alpha,u)$ be a cocycle action of a poly-$\Z$ group $G$ on a unital $C^*$-algebra $A$ in $\AC(\cO_\infty,\mu^G)$,  
and let $J$ be a closed ideal of $A$ with $u(g_1,g_2)\in U(J)_0$ for any $g_1,g_2\in G$ 
that is globally preserved by $\alpha$. 
Then for any continuous map $v:G\times EG\to U(J)_0$ satisfying $u(g_1,g_2)=v(g_1,g_2\cdot 0_n)$ and 
$$u(g_1,g_2)v(g_1g_2,z)=\alpha_{g_1}(v(g_2,z))v(g_1,g_2\cdot z),$$
for any $g_1,g_2\in G$ and $z\in EG$, 
there exists a continuous map $w:EG\to U(J)_0$ satisfying 
$$w(g\cdot 0_n)\alpha_g(w(z))v(g,z)w(g\cdot z)^{-1}=1,$$
for any $g\in G$, $z\in EG$. 
\end{theorem}

To prove Theorem \ref{vanishing2}, we need the following lemma.  

\begin{lemma}\label{small section} Let $G$ be a poly-$\Z$ group, let $S\subset G$ be a finite symmetric generating set, 
and let $K\subset EG$ be a compact subset. 
Then for any $0<\varepsilon<2$, there exists $\delta>0$ satisfying the following (1) and (2). 
If moreover $K$ is a finite subcomplex of $EG$ satisfying $0_n\in K$ and $G\cdot K=EG$, 
we can choose $\delta =\frac{\varepsilon}{m(n+1)}$, where $m$ is a natural number satisfying 
$\{g\in G;\;  (g\cdot K)\cap K\neq \emptyset \}\subset S^m$. 

\begin{itemize}
\item[$(1)$] 
Let $\alpha$ be an action of $G$ on a unital $C^*$-algebra $A$, and 
let $u=\{u_g\}_{g\in G}$ and $u'=\{u'_g\}_{g\in G}$ be $\alpha$-cocycles in $U(A)_0$ satisfying 
$\|u_g-u'_g\|<\delta$ for any $g\in S$. 
Then for any homotopy fixed point $F:EG\to U(A)_0$ for $u$, there exists a homotopy fixed point 
$F':EG\to U(A)_0$ for $u'$ satisfying $F(0_n)=F'(0_n)$ and 
$\|F(x)-F'(x)\|<\varepsilon$ for any $x\in K$. 
\item[$(2)$] 
Let $\alpha$ be an action of $G$ on a unital $C^*$-algebra $A$, and  
let $u=\{u_g\}_{g\in G}$ and $u'=\{u'_g\}_{g\in G}$ be $\id\otimes \alpha$-cocycles in $U(C_0([0,1))\otimes A)_0$ satisfying 
$\|u_g-u'_g\|<\delta$ for all $g\in S$ and $u_g(0)=u'_g(0)$ for all $g\in G$. 
Then for any homotopy fixed point $F:EG\to U(C_0([0,1))\otimes A)_0$ for $u$, there exists a homotopy fixed point 
$F':EG\to U(A)_0$ for $u'$ satisfying $F(0_n)=F'(0_n)$, $F(x)(0)=F'(x)(0)$ for all $x\in EG$, and 
$\|F(x)-F'(x)\|<\varepsilon$ for all $x\in K$. 
\end{itemize}
\end{lemma}

\begin{proof} (1) We fix a free $G$-CW-complex structure of $EG$ with $0_n$ in the 0-skeleton, 
and choose a finite subcomplex $K_0$ of $EG$ satisfying $0_n\in K_0$, $K\subset K_0$, and $G\cdot K_0=EG$. 
Let 
$$S_1=\{g\in G;\; (g\cdot K_0)\cap K_0\neq \emptyset \}.$$ 
Then $S_1$ is a finite subset of $G$, and there exists $m\in \N$ satisfying $S_1\subset S^m$. 
We set $\delta=\frac{\varepsilon}{m(n+1)}$. 
Assume that $u_g$, $u'_g$ and $F$ satisfy the condition of (1). 
Then we have $\|u_g-u'_g\|<\frac{\varepsilon}{n+1}$ for all $g\in S_1$. 

We claim that if $f:D^k\to U(A)_0$ and $f_0':\partial D^k\to U(A)_0$ 
are continuous maps satisfying 
$$l:=\max_{x\in \partial D^k}\|f(x)-f_0'(x)\|<2,$$
there exists a continuous extension $f':D^k\to U(A)_0$ of $f'_0$ satisfying $\|f(x)-f'(x)\|\leq l$ for any $x\in D^k$. 
Indeed, it follows from the following fact: 
if $f''_0:\partial D^k\to U(A)_0$ is a continuous map satisfying 
$$l:=\max_{x\in \partial D^k}\|f''_0(x)-1\|<2,$$
there exists a continuous extension $f'':D^k\to U(A)_0$ of $f''_0$ satisfying $\|f''(x)-1\|\leq l$ for any $x\in D^k$. 
Indeed, we can set 
$$f''(x)=e^{|x|\log f''_0(\frac{x}{|x|})}$$
for $x\in  D^k\setminus \{0_n\}$ and set $f''(0)=1$. 
Applying this fact to $f''_0(x)=f(x)^{-1}f_0'(x)$ for $x\in \partial D^k$,  
and setting $f'(x)=f(x)f''(x)$ for $x\in D^k$, we get the claim. 

By a cell-by-cell argument using the above claim, we can construct a continuous map $F'_0:K_0\to U(A)_0$ satisfying 
$F'_0(0_n)=F(0_n)$, $\|F(x)-F'_0(x)\|<\frac{(k+1)\varepsilon}{n+1}$ for any $x$ in the $k$-skeleton $K_0^k$ of $K_0$, 
and $F'_0(g\cdot x)=u'_g\alpha_g(F'_0(x))$ 
for any $x\in K_0$ and $g\in S_1$ with $g\cdot x\in K_0$. 
Indeed, we can show this claim by induction of $k$.  
For $x=g\cdot 0_n\in K_0^0$, we set $F_0'(x)=u_g'\alpha_g(F(0_n))$, which satisfies 
$\|F(x)-F_0'(x)\|<\frac{\varepsilon}{n+1}$. 
Assume that we have $F_0'$ defined on $K_0^k$ satisfying the above condition. 
Let $f:D^{k+1}\to K_0$ be the characteristic map of a $k+1$-cell of $K_0$. 
Then we can extend $F_0'$ to $K_0^k\cup f(D^{k+1})$ keeping the estimate $\|F(x)-F_0'(x)\|<\frac{(k+1)\varepsilon}{n+1}$. 
If $g\in S_1$ satisfies $g\cdot f(D^{k+1})\cap K_0\neq \emptyset$, we can further extend $F_0'$ by setting 
$F_0'(g\cdot x)=u'_g\alpha_g(F_0'(x))$ for $x\in f(D^{k+1})$ with estimate 
$$\|u_g\alpha_g(F(x))-u_g'\alpha_g(F_0'(x))\|<\|u_g-u'_g\|+\|F(x)-F_0'(x)\|<\frac{(k+2)\varepsilon}{n+1}.$$
Repeating the same argument, we can extend $F_0'$ to $K_0^{k+1}$. 
Thus the claim is shown by induction. 

The map $F'_0$ extends to a homotopy fixed point $F'$ for $u'$ with the desired properties. 

(2) follows from essentially the same argument as above using $D^k\times I$ instead of $D^k$. 
Indeed, note that $D^k\times I\cong D^{k+1}$ and $\partial (D^k\times I)=(\partial D^k\times I)\cup (D_k\times \{0,1\})$. 
Thus the boundary condition $F'(x)(y)=F(x)(y)$ for $y\in \{0,1\}$ can be handled by $f'_0$ in (1).     
\end{proof}

We recall the fact that every poly-$\Z$ group satisfies the following $H^2$-stability \cite[Theorem 5.13]{IMI}. 

\begin{theorem}\label{H2-stability} For every poly-$\Z$ group $G$, there exists a finite generating subset $S$ of $G$ 
such that for every $\varepsilon>0$ there exists $\delta>0$ satisfying the following property. 
Suppose $(\alpha, u)$ is a cocycle action in $\mathrm{AC}(\cO_\infty,\mu^G)$ and $J\subset A$ 
is a globally $\alpha$-invariant closed ideal with $u(g,h)\in U(J)_0$ for all $g,h\in G$. 
If $\|u(g,h)-1\|<\delta$ for all $g,h\in S$, then there exists a family of unitaries $\{s_g\}_g$ 
in $U(J)_0$ such that 
$$s_g\alpha_g(s_h)u(g,h)s_{gh}^*=1,\quad \forall g,h\in G,$$
and $\|s_g-1\|<\varepsilon$ for all $g\in S$. 
\end{theorem}

\begin{remark}\label{H2-stability'}
Assume that $G$ is a poly-$\Z$ and $G=G_{n-1}\rtimes \langle{\xi}\rangle$ is the semi-direct product decomposition as before. 
In the proof of \cite[Theorem 5.13]{IMI} (or actually \cite[Lemma 5.11]{IMI}), the finite generating set $S$ in the statement 
of Theorem \ref{H2-stability} is constructed in the following form: $S=T\cup \{\xi\}$ with some symmetric finite generating set $T$ of $G_{n-1}$. 
Moreover, among the condition 
$$\|u(g,h)-1\|<\delta,\quad \forall g,h\in S,$$ 
in the statement, we can see that the condition $\|u(\xi,\xi)-1\|<\delta$ can be omitted. 
\end{remark}

\begin{proof}[Proof of Theorem \ref{vanishing2}] 
We show the statement by induction of the Hirsch length of $G$. 
The statement is trivially true when the Hirsch length is 0. 
Assume that the statement is true for $G_{n-1}=N$. Recall that $G$ is generated by $N$ and $\xi=\xi_n$. 

Step I. By induction hypothesis, there exists a continuous map $w_0:EN\to U(J)_0$ satisfying 
$$w_0(h\cdot 0_{n-1})\alpha_h(w_0(x))v(h,(x,0))w_0(h\cdot x)^{-1}=1$$
for any $h\in N$ and $x\in EN$. 
Note that we automatically have $w_0(0_{n-1})=1$. 
We can extend $w_0$ continuously to $w:EG\to U(J)_0$ with $w(x,0)=w_0(x)$ and 
$$w(x,1)=w(\xi\cdot 0_{n-1},1)\alpha_\xi(w(\xi^{-1}\cdot x,0))v(\xi,(\xi^{-1}\cdot x,0)),$$
where we can choose $w(\xi\cdot 0_{n-1},1)\in U(J)_0$ arbitrarily. 
Then we have 
$$w(\xi\cdot 0)\alpha_\xi(w(x,0))v(\xi,(x,0))w(\xi\cdot(x,0))^{-1}=1.$$ 
Thus by replacing $\alpha_g, u(g_1,g_2), v(g,z)$ with $\Ad w(g\cdot 0_n)\circ \alpha_g$, 
$$w(g_1\cdot0_n)\alpha_{g_1}(w(g_2\cdot 0_n))u(g_1,g_2)w(g_1g_2\cdot 0_n)^{-1},$$    
$$w(g\cdot 0_n)\alpha_g(w(z))v(g,z)w(g\cdot z)^{-1},$$ 
respectively, we may and do assume $u(h,k)=u(\xi,h)=v(h,(x,0))=v(\xi,(x,0))=1$ for any $h,k\in N$ and $x\in EN$ from the beginning.

Step II. Note that 
$$u(\xi,h^\xi)v(\xi h^\xi,(x,0))=\alpha_\xi(v(h^\xi,(x,0)))v(\xi,h^\xi\cdot(x,0))$$
implies $v(h\xi,(x,0))=1$ for any $h\in N$ and $x\in EN$, and 
$$u(h,\xi)v(h\xi,(x,0))=\alpha_h(v(\xi,(x,0)))v(h,\xi\cdot(x,0))$$
implies $u(h,\xi)=v(h,\xi\cdot (x,0))$ for any $x\in EN$. 
Therefore $v(h,\xi\cdot (x,0))$ does not depend on $x\in EN$, 
and the family $\{u(h,\xi)^{-1}\}_{h\in N}$ forms an $\alpha|_N$-cocycle. 

For $y\in I$ and $h,k\in N$, let $\beta_h^{(y)}=\Ad v(h,(0_{n-1},y))^{-1}\circ \alpha_h$ and, let  
\begin{align*}
\lefteqn{u'(h,k)(y)} \\
 &=v(h,(0_{n-1},y))^{-1}\alpha_h(v(k,(0_{n-1},y))^{-1})v(hk,(0_{n-1},y))\\
 &=v(h,(0_{n-1},y))^{-1}v(h,(k\cdot 0_{n-1},y)),
\end{align*}
 for $h,k\in N$. 
Since $v(h,\xi\cdot(x,0))=v(h,(\xi\cdot x,1))$, which does not depend on $x$, we get $u'(h,k)\in U(SJ)$. 
Then $(\beta,u')$ gives a cocycle action of $N$ on $C(I)\otimes A$, 
where $\beta_h(a)(y)=\beta_h^{(y)}(a(y))$. 
We have $(\beta,u')\in \AC(\cO_\infty,\mu^N)$.
Let $v'(h,x)(y)=v(h,(0_{n-1},y))^{-1}v(h,(x,y))$. 
Then $v': N\times EN\to U(SJ)$ is a continuous map satisfying $v'(h,0_{n-1})=1$ and 
$$u'(h,k)v'(hk,x)=\beta_h(v'(k,x))v'(h,k\cdot x),$$
and hence $v'(h,x)\in U(SJ)_0$. 
Since $u'(h,k)=v'(h,k\cdot 0_{n-1})$ we also have $u'(h,k)\in U(SJ)_0$.  
Thus the induction hypothesis implies that there exists a continuous map $w':EN\to U(SJ)_0$ satisfying 
$$w'(h\cdot 0_{n-1})\beta_h(w'(x))v'(h,x)w'(h\cdot x)^{-1}=1,$$
and so \begin{align*}
\lefteqn{1=w'(h\cdot 0_{n-1})(y)\beta^{(y)}_h(w'(x)(y))v(h,(0_{n-1},y))^{-1}v(h,(x,y))w'(h\cdot x)(y)^{-1}} \\
 &=w'(h\cdot 0_{n-1})(y)v(h,(0_{n-1},y))^{-1}\alpha_h(w'(x)(y)) v(h,(x,y))w'(h\cdot x)(y)^{-1}, \\
\end{align*}
for any $h\in H$, $x\in EN$, and $y\in I$.  
This shows 
$$w'(h\cdot 0_{n-1})(0)\alpha_h(w'(x)(y))v(h,(x,y))w'(h\cdot x)(y)^{-1}$$
does not depend on $x\in EN$. 
We also have 
$$w'(\xi\cdot 0_{n-1})(1)\alpha_\xi(w'(x)(0))v(\xi,(x,0))w'(\xi\cdot x)(1)^{-1}=1\cdot 1\cdot 1\cdot 1=1.$$
We define $w_1:EN\times I\to U(J)_0$ by  $w_1(x,y)=w'(x)(y)$. 
Then we have $w_1(x,0)=w_1(x,1)=1$, 
$$w_1(h\cdot 0_n)\alpha_h(w_1(x,y))v(h,(x,y))w_1(h\cdot (x,y))^{-1},$$
does not depend on $x\in EN$ for $y\in I$, and it is 1 at $y=0$. 
We also have 
$$w_1(\xi\cdot 0_n)\alpha_\xi(w_1(x,0))v(\xi,(x,0))w_1(\xi\cdot (x,0))^{-1}=1,$$
for any $x\in EN$. 
Thus from the beginning, we may and do assume that $v(h,(x,0))=v(\xi,(x,0))=1$ for any $h\in N$, $x\in EN$, 
and $v(h,(x,y))$ does not depend on $x\in EN$ for any $h\in N$ and $y\in I$. 

Step III. We choose a finite symmetric generating set $S\subset N$ so that $H^2$-stability holds for $S\cup \{\xi\}$ 
(see Theorem \ref{H2-stability} and Remark \ref{H2-stability'}). 
We fix a free $N$-CW-complex structure of $EN$ with $0_{n-1}$ in the 0-skeleton, 
and choose a finite subcomplex $K_0$ of $EN$ satisfying $0_{n-1}\in K_0$ and $N\cdot K_0=EN$. 
We choose a finite subcomplex $K$ of $EN$ satisfying $K_0\cup (\xi^{-1}\cdot K_0)\subset K$, and let 
$S_1=\{h\in N;\; (h\cdot K)\cap K\neq \emptyset\}$. 
Taking $m\in \N$ satisfying $S_1\subset S^m$, we set $\varepsilon=\frac{1}{2mn}$. 
For this $\varepsilon$, Theorem \ref{H2-stability} and the Remark \ref{H2-stability'} imply that there exists $\delta>0$ satisfying the following property: 
if $u(h,k)=u(\xi,h)=1$ and $\|u(h,\xi)-1\|<\delta$ for any $h,k\in S$, there exist $s_g\in U(J)_0$ for $g\in G$ satisfying 
$s_{g_1}\alpha_{g_1}(s_{g_2})u(g_1,g_2)s_{g_1g_2}^{-1}=1$ for any $g_1,g_2\in G$, and $\|s_g-1\|<\varepsilon$ for any $g\in S\cup \{\xi\}$. 
We may assume $\delta\leq \varepsilon$. 

Recall that $v(h,(x,y))$ for $h\in N$ and $y\in I$ does not depend on $x\in EN$. 
For any $h,k\in N$,  
$$u(h,k)v(hk,(x,y))=\alpha_h(v(k,(x,y)))v(h,(k\cdot x,y))$$
implies that $\{v(h,(0_{n-1},y))^{-1}\}_{h\in N,\; y\in I}$ is a continuous family 
of $\alpha|_N$-cocycles that is 1 at $y=0$. 
Thus \cite[Lemma 5.6]{IMI} implies that there exists a continuous map $r:[0,1]\to U(J)_0$ 
satisfying $r(0)=1$ and
$$\sup_{y\in I}\|v(h,(x,y))^{-1}-r(y)^{-1}\alpha_h(r(y))\|<\delta,$$
for any $h\in S$ and $x\in EN$.  
Defining a continuous map $w_2: EN\times I\to U(J)_0$ by $w_2(x,y)=r(y)$, we get 
$$\sup_{y\in I}\|w_2(h\cdot 0_n)\alpha_h(w_2(x,y))v(h,(x,y))w_2(h\cdot(x,y))^{-1}-1\|<\delta,$$
for any $h\in S$ and $x\in EN$. 
We also have 
$$w_2(\xi\cdot 0_n)\alpha_\xi(w_2(x,0))v(\xi,(x,0))w_2(\xi\cdot (x,0))^{-1}
=r(1)r(1)^{-1}=1.$$
Thus from the beginning, we may and do further assume $\|v(h,(x,y))-1\|<\delta$ for any $h\in S$, $x\in EN$, $y\in I$. 

Now since $u(h,k)=1$, $u(\xi,h)=1$, $\|u(h,\xi)-1\|<\delta$ for any $h,k\in S$,  
we get $s_g\in U(J)_0$ for $g\in G$ as in Theorem \ref{H2-stability}. 
In particular, the family $\{s_h\}_{h\in N}$ forms a $\alpha|_N$-cocycle, $\|s_h-1\|<\varepsilon$ for any $h\in S$ 
and $\|s_\xi-1\|<\varepsilon$.  
We also have $s_\xi\alpha_\xi(s_h)s_{\xi h}^{-1}=1$ and $s_h\alpha_h(s_\xi)u(h,\xi)s_{h\xi}^{-1}=1$ for any $h\in N$.  

We claim that there exists a continuous map $w_3:EN\times I\to U(J)_0$ satisfying 
$$w_3(h\cdot 0_{n-1},0)\alpha_h(w_3(x,y))v(h,(x,y))w_3(h\cdot(x,y))^{-1}=1,$$
$$w_3(\xi\cdot 0_n)\alpha_\xi(w_3(x,0))v(\xi,(x,0))w_3(\xi\cdot (x,0))^{-1}=1,$$
$w_3(h\cdot 0_{n-1},0)=s_h$, and $w_3(\xi\cdot 0_n)=s_\xi$ for any $h\in N$, $x\in EN$, and $y\in I$.  
The first two equations can be rephrased as 
\begin{equation}\label{E1}
w_3(h\cdot x,y)=s_h\alpha_h(w_3(x,y))v(h,(0_{n-1},y)),
\end{equation}
\begin{equation}\label{E2}
w_3(x,1)=s_\xi\alpha_\xi(w_3(\xi^{-1}\cdot x,0)).
\end{equation}
Note that $\{s_h\}_{h\in N}$ and $\{v(h,(0_{n-1},y))^{-1}\}_{h\in N}$ are $\alpha|_N$-cocycles satisfying 
$\|s_h-1\|<\varepsilon$ and $\|v(h,(0_{n-1},y))^{-1}-1\|<\varepsilon$ for all $h\in S$. 
We also have $v(h,0_n)=1$ for all $h\in N$. 

We first apply Lemma \ref{small section},(1) to $N$ (in place of $G$), $S$, $K$, $mn\varepsilon$ (in place of $\varepsilon$), 
$1$ (in place of $u_g$), $s_g$ (in place of $u'_g$), and $1$ (in place of $F$). 
This is possible as we have $S_1\subset S^m$. 
Then we get a homotopy fixed point $f:EN \to U(J)_0$ for $\{s_h\}_{h\in N}$ satisfying $f(0_{n-1})=1$, $f(h\cdot x)=s_hf(x)$, 
and $\|f(x)-1\|<mn\varepsilon$ for all $x\in K$ and $h\in N$. 
We set $w_3(x,0)=f(x)$ and $w_3(x,1)=s_\xi\alpha_\xi(f(\xi^{-1}\cdot x))$ for $x\in EN$. 
Then 
$$w_3(h\cdot 0_{n-1},0)=f(h\cdot0_{n-1})=s_hf(0_{n-1})=s_h,$$ 
$$w_3(\xi\cdot 0_n)=s_\xi\alpha_\xi(f(0_{n-1}))=s_\xi,$$
and Eq.(\ref{E1}) with $y=0$ and Eq.(\ref{E2}) are satisfied. 
Since 
\begin{align*}
\lefteqn{w_3(h\cdot x,1)=s_\xi\alpha_\xi(f(\xi^{-1}h\cdot x))=s_\xi\alpha_\xi(f(h^{\xi}\xi^{-1}\cdot x))} \\
 &= s_\xi\alpha_\xi(s_{h^\xi}\alpha_{h^\xi}(f(\xi^{-1}\cdot x)))= s_\xi\alpha_\xi(s_{h^\xi})\alpha_{\xi h^\xi}(f(\xi^{-1}\cdot x))\\
 &=s_\xi\alpha_\xi(s_{h^\xi})u(h,\xi)^{-1}\alpha_h(\alpha_\xi(f(\xi^{-1}\cdot x))u(h,\xi) \\
 &=s_{\xi h^\xi}u(h,\xi)^{-1}\alpha_h(s_\xi^{-1})\alpha_h(w_3(x,1))v(h,(\xi\cdot 0_{n-1},1))\\
 &=s_h\alpha_h(w_3(x,1))v(h,(0_{n-1},1)),
\end{align*}
Eq.(\ref{E1}) with $y=1$ is also satisfied. 
We have $\|w_3(x,0)-1\|<mn\varepsilon$ for all $x\in K$. 
Since $\xi^{-1}\cdot K_0\cup K_0\subset K$, we have $\|w_3(x,1)-1\|<(mn+1)\varepsilon<2mn\varepsilon$ for all $x\in K_0$. 
We would like to extend $w_3$ to the whole $EN\times I$ to get a continuous map satisfying Eq.(\ref{E1}). 
For this purpose, since $N\cdot K_0=EN$, it suffices to continously extend $w_3$ to $K_0\times I$ so that 
Eq.(\ref{E1}) holds whenever $x$ and $h\cdot x$ are in $K_0$. 
We can do it by a cell-by-cell argument. 

Let $K_0^{k}$ be the $k$ skeleton of $K_0$. 
As in the proof of Lemma \ref{small section} using $D^k\times I$ instead of $D^k$, we can continuously extend $w_3$ 
to $K_0\times I$ by induction of $k$ so that Eq.(\ref{E1}) holds whenever $x,h\cdot x\in K_0$,  
and $\|w_3(x,y)-1\|<2mn\varepsilon+2(k+1)m\varepsilon $ for all $(x,y)\in K_0^k\times I$.  
Indeed, since $w_3(0_n)=1$ and $\|w_3(0_{n-1},1)-1\|<2mn\varepsilon<2$, we can continuously extend $w_3$ to $\{0_{n-1}\}\times I$ 
so that $\|w_3(0_{n-1},y)-1\|<2mn\varepsilon$ for all $y\in I$. 
For $h\cdot 0_{n-1}\in K_0^0$ with $h\neq e$, we set $w_3(h\cdot 0_{n-1},y)=s_h\alpha_h(w_3(0_{n-1},y))v(h,(0_{n-1},y))$. 
Then we get $\|w_3(x,y)-1\|<4mn\varepsilon$ for all $(x,y)\in K_0^0\times I$. 
By induction of $k$, we can extend $w_3$ to $K_0\times I$ as desired with the estimate 
$\|w_3(x,y)-1\|<4mn\varepsilon=2$ for all $(x,y)\in K_0\times I$. 
Now we can extend $w_3$ to the whole $EN\times I$ by the relation Eq.(\ref{E1}), and the claim is shown.

Now extending $w_3$ continuously to the whole $EG$ and replacing $v(g,z)$ with 
$$w_3(g\cdot 0_n)\alpha_g(w_3(z))v(g,z)w_3(g\cdot z)^{-1},$$ 
we may and do assume that $v(h,(x,y))=v(\xi,(x,0))=1$ for any $h\in N$, $x\in EN$, and $y\in I$. 

Step IV. The rest of the proof is routine. 
First we set $w(x,y)=1$ for $x\in EN$ and $y\in I$. 
Secondly, using $w(\xi\cdot 0_n)=w(0_n)=1$, we define $w$ on $\xi^m\cdot 0_n=(\xi^m\cdot 0_{n-1},m)$ for $m\in \Z$ by 
$$w(\xi^m\cdot 0_n)=w(\xi^{m-1}\cdot 0_n)u(\xi^{m-1},\xi).$$ 
The relation 
$$u(\xi^{m-1},\xi)v(\xi^m,(\xi^{-m}\cdot x,0))
 =\alpha_{\xi^{m-1}}(v(\xi,(\xi^{-m}\cdot x,0)))v(\xi^{m-1},(\xi^{1-m}\cdot x,1)),$$
implies 
$$u(\xi^{m-1},\xi)v(\xi^m,(\xi^{-m}\cdot x,0))=v(\xi^{m-1},(\xi^{1-m}\cdot x,1)),$$
and so 
$$w(\xi^m\cdot 0_n)v(\xi^m,(\xi^{-m}\cdot x,0))=w(\xi^{m-1}\cdot 0_n)v(\xi^{m-1},(\xi^{1-m}\cdot x,1)).$$
This allows us to set 
$$w(x,m+y)=w(\xi^m\cdot 0_n)v(\xi^m,(\xi^{-m}\cdot x,y)),$$
for any $x\in EN$, $y\in I$ and $m\in \Z$, and we get a well-defined continuous map $w:EG\to U(J)_0$. 
Let 
$$v_1(g,z)=w(g\cdot 0_n)\alpha_g(w(z))v(g,z)w(g\cdot z)^{-1}.$$
Then by construction, we have $v_1(h,(x,y))=v_1(\xi^m,(x,y))=1$ for any $m\in \Z$, $h\in N$, $x\in EN$, and $y\in I$. 
We claim that $v_1(\xi^m,z)=1$ holds for any $m\in \Z$ and $z\in EG$. 
Indeed, letting $u_1(g_1,g_2)=v_1(g_1,g_2\cdot 0_n)$, we have $u_1(\xi^m,\xi)=1$ for any $m\in \Z$. 
This together with the 2-cocycle relation of $u_1$ implies $u_1(\xi^m,\xi^l)=1$ for any $m,l\in \Z$. 
Now the claim follows from the relation
$$u_1(\xi^m,\xi^l)v_1(\xi^{m+l},(x,y))=\alpha_{\xi^m}(v_1(\xi^l,(x,y)))v_1(\xi^m,\xi^l\cdot(x,y)).$$

For $x\in EN$, $m\in \Z$, and $y\in I$, 
\begin{align*}
\lefteqn{v_1(h,(x,m+y))} \\
 &=w(h\cdot 0_n)\alpha_h(w(x,m+y))v(h,(x,m+y))w(h\cdot x,m+y)^{-1} \\
 &=\alpha_h(w(\xi^m\cdot 0_n)v(\xi^m,(\xi^{-m}\cdot x,y)))v(h,(x,m+y))\\
 &\quad\times v(\xi^m,(\xi^{-m}h\cdot x,y))^{-1}w(\xi^m\cdot 0_n)^{-1} \\
 &=\alpha_h(w(\xi^m\cdot 0_n)) u(h,\xi^m)v(h\xi^m,(\xi^{-m}\cdot x,y))v(\xi^m,(\xi^{-m}h\cdot x,y))^{-1}w(\xi^m\cdot 0_n)^{-1}\\
 &=\alpha_h(w(\xi^m\cdot 0_n)) u(h,\xi^m)\\
 &\quad \times v(\xi^mh^{\xi^m},(\xi^{-m}\cdot x,y))v(\xi^m,(h^{\xi^m}\xi^{-m}\cdot x,y))^{-1}w(\xi^m\cdot 0_n)^{-1}\\
 &=\alpha_h(w(\xi^m\cdot 0_n)) u(h,\xi^m)u(\xi^m,h^{\xi^m})^{-1}\alpha_{\xi^m}(v(h^{\xi^m},(\xi^{-m}\cdot x,y))) w(\xi^m\cdot 0_n)^{-1}\\
 &=\alpha_h(w(\xi^m\cdot 0_n)) u(h,\xi^m)u(\xi^m,h^{\xi^m})^{-1} w(\xi^m\cdot 0_n)^{-1},
\end{align*}
and to show $v_1(h,(x,m+y))=1$, it suffices to show that 
$$\alpha_h(w(\xi^m\cdot 0_n)) u(h,\xi^m)u(\xi^m,h^{\xi^m})^{-1} w(\xi^m\cdot 0_n)^{-1}=1$$
by induction of $m$. 
Note that it is the case for $m=0$ and $m=1$. 
By construction, we have \begin{align*}
\lefteqn{\alpha_h(w(\xi^m\cdot 0_n)) u(h,\xi^m)u(\xi^m,h^{\xi^m})^{-1} w(\xi^m\cdot 0_n)^{-1}} \\
 &=\alpha_h(w(\xi^{m-1}\cdot 0_n)u(\xi^{m-1},\xi)) u(h,\xi^m)u(\xi^m,h^{\xi^m})^{-1}u(\xi^{m-1},\xi)^{-1} w(\xi^{m-1}\cdot 0_n)^{-1} \\
 &=\alpha_h(w(\xi^{m-1}\cdot 0_n))u(h,\xi^{m-1})u(h\xi^{m-1},\xi)\\
 &\quad \times u(\xi^{m-1},\xi h^{\xi^m})^{-1}\alpha_{\xi^{m-1}}(u(\xi,h^{\xi^m})^{-1})w(\xi^{m-1}\cdot 0_n)^{-1} \\
 &=\alpha_h(w(\xi^{m-1}\cdot 0_n))u(h,\xi^{m-1})u(\xi^{m-1}h^{\xi^{m-1}},\xi) 
 u(\xi^{m-1},h^{\xi^{m-1}}\xi)^{-1}w(\xi^{m-1}\cdot 0_n)^{-1} \\
  &=\alpha_h(w(\xi^{m-1}\cdot 0_n))u(h,\xi^{m-1})
  u(\xi^{m-1},h^{\xi^{m-1}})^{-1}\alpha_{\xi^{m-1}}(u(h^{\xi^{m-1}},\xi))w(\xi^{m-1}\cdot 0_n)^{-1} \\
  &=\alpha_h(w(\xi^{m-1}\cdot 0_n))u(h,\xi^{m-1})
  u(\xi^{m-1},h^{\xi^{m-1}})^{-1}w(\xi^{m-1}\cdot 0_n)^{-1}. 
\end{align*}
Thus we get $v_1(h,z)=1$ for any $h\in N$ and $z\in EG$. 

Now the relation 
$$u_1(h,\xi^m)v_1(h\xi^m,z)=\alpha_h(v_1(\xi^m,z))v_1(h,\xi^m\cdot z)$$
for $h\in N$, $m\in \Z$, and $z\in BG$ implies $v_1(g,z)=1$ for any $g\in G$ and $z\in EG$. 
\end{proof}

\begin{proof}[Proof of Theorem \ref{vanishing1}] Since $(BG,\{[0_n]\})$ is an NDR pair, we may and do assume $F(0_n)=1$ 
thanks to Theorem \ref{lifting extension}. 
For $g\in G$, $x\in EG$ and $t\in I$, we set $\tilde{\alpha}^{(t)}_g=\Ad F(tg\cdot 0_n)\circ \alpha_g$, and 
$$v(g,x)(t)=F(t g\cdot 0_n)\alpha_g(F(tx))F(tg\cdot x)^{-1}.$$
Then we have $v(g,x)(0)=v(g,x)(1)=1$, and $v:G\times EG\to U(SJ)$ is a continuous map. 
We also have $v(g,0_n)=1$ and the image of $v$ is in $U(SJ)_0$. 
Let $u_1(g,h)=v(g,h\cdot 0_n)$. 
Then by construction $(\tilde{\alpha},u_1)$ is a cocycle action of $G$ on $C(I)\otimes A$ in $\AC(\cO_\infty,\mu^G)$. 

For any $g,h\in G$, $x\in EG$, and $t\in I$, we have 
\begin{align*}
\lefteqn{u_1(g,h)(t)v(gh,x)(t)} \\
 &=F(tg\cdot 0_n)\alpha_g(F(th\cdot 0_n))F(tgh\cdot 0_n)^{-1}F(tgh\cdot 0_n)\alpha_{gh}(F(tx))F(tgh\cdot x)^{-1} \\
 &=F(tg\cdot 0_n)\alpha_g(F(th\cdot 0_n))\alpha_{gh}(F(tx))F(tgh\cdot x)^{-1},
\end{align*}
\begin{align*}
\lefteqn{\tilde{\alpha}^{(t)}_g(v(h,x)(t))v(g,h\cdot x)(t)} \\
 &=F(tg\cdot 0_n)\alpha_g(F(th\cdot 0_n)\alpha_h(F(tx))F(th\cdot x)^{-1})F(tg\cdot 0_n)^{-1}\\
 &\quad \times F(tg\cdot 0_n)\alpha_g(F(th\cdot x))F(tgh\cdot x)^{-1} \\
 &=F(tg\cdot 0_n)\alpha_g(F(th\cdot 0_n))\alpha_{gh}(F(tx))F(tgh\cdot x)^{-1},
\end{align*}
and so $u_1(g,h)v(gh,x)=\tilde{\alpha}_g(v(h,x))v(g,h\cdot x)$ holds. 
Thus thanks to Theorem \ref{vanishing2}, there exists a continuous map $w:EG\to U(SJ)_0$ satisfying 
$$w(g\cdot 0_n)\tilde{\alpha}_g(w(x))v(g,x)w(g\cdot x)^{-1}=1,$$ which means 
$$w(g\cdot 0_n)(t)F(tg\cdot 0_n)\alpha_g(w(x)(t)F(tx))=w(g\cdot x)(t)F(tg\cdot x).$$
In particular, $\{w(g\cdot 0_n)(t)F(tg\cdot 0_n)\}_{g\in G,\; t\in I}$ is a continuous family of $\alpha$-cocycles 
giving $1$ at $t=0$ and $u_g$ at $t=1$. 
We set 
$$\tilde{u}_g(t)=w(g\cdot 0_n)(1-t)F((1-t)g\cdot 0_n),$$
$$\tilde{F}_0(x)(t)=w(x)(1-t)F((1-t)x).$$ 
Then $\{\tilde{u}_g\}_{g\in G}$ is an $\id_{C(I)}\otimes \alpha$-cocycle in 
$U(C_0([0,1))\otimes J)_0$ with $\tilde{u}_g(0)=u_g$, and 
$\tilde{F}_0:EG\to U(C_0([0,1))\otimes J)_0$ is a homotopy fixed point for $\tilde{u}$ with $\tilde{F}_0(x)(0)=F(x)$. 

We choose a finite symmetric generating set $S\subset G$ and 
a finite subcomplex $K\subset EG$ satisfying $G\cdot K=EG$. 
Let $\delta>0$ be the constant derived from Lemma \ref{small section} with these $S$ and $K$ and $\varepsilon=1$. 
Since the cocycle $\{\tilde{u}_g\}_{g\in G}$ is homotopic to 1 within the set of $\id_{C(I)}\otimes \alpha$-cocycles 
in $U(C_0([0,1))\otimes J)_0$, it is an asymptotic coboundary in $U(C_0([0,1))\otimes J)_0$ thanks to 
Lemma 5.6 and Lemma 5.13 in \cite{IMI}. 
Thus there exists a continuous family $\{r_0(s)\}_{s\in [0,1)}$ of unitaries 
in $U(C_0([0,1))\otimes J)_0$ satisfying 
$$\lim_{s\to 1-0}\|r_0(s)(\id\otimes \alpha_g)(r_0(s)^{-1})-\tilde{u}_g\|=0.$$ 
We may further assume $\|r_0(s)(\id\otimes \alpha_g)(r_0(s)^{-1})-\tilde{u}_g\|<\varepsilon/2$ 
for any $s\in [0,1)$ and $g\in S$. 
Writing $r_0(s)(t)=r_0(s,t)$, this means that $r_0:[0,1)\times I\to U(J)_0$ is a continuous map 
with $r(s,1)=1$ for any $s\in [0,1)$, 
$$\lim _{s\to 1-0}\sup_{t\in I}\|r_0(s,t)\alpha_g(r_0(s,t)^{-1})-\tilde{u}_g(t)\|=0,$$
for any $g\in G$, and 
$$\sup_{t\in I}\|r_0(s,t)\alpha_g(r_0(s,t)^{-1})-\tilde{u}_g(t)\|<\varepsilon/2, $$
for any $s\in [0,1)$ and $g\in S$. 
Choosing an orientation preserving homeomorphism $f$ of $I$ satisfying $\lim_{t\to 1-0}r_0(f(t),t)=1$, 
we set $r(t)=r_0(t,0)r_0(f(t),t)^{-1}$, which is in $U(J)_0$. 
Then 
\begin{align*}
\lefteqn{\|r(t)^{-1}u_g\alpha_g(r(t))-\tilde{u}_g(t)\|} \\
 &=\|r_0(f(t),t)r_0(t,0)^{-1}\tilde{u}_g(0)\alpha_g(r_0(t,0)r_0(f(t),t)^{-1})-\tilde{u}_g(t)\| \\
 &\leq\|r_0(t,0)^{-1}\tilde{u}_g(0)\alpha_g(r_0(t,0))-1\|
 +\|r_0(f(t),t)\alpha_g(r_0(f(t),t)^{-1})-\tilde{u}_g(t)\|<\varepsilon,
\end{align*}
for any $g\in S$, $t\in [0,1)$, and 
$$\lim_{t\to 1-0}r(t)^{-1}u_g\alpha_g(r(t))=1,$$
for any $g\in G$. 

Let 
$$\tilde{u}'_g(t)=\left\{
\begin{array}{ll}
r(t)^{-1}u_g\alpha_g(r(t)) , &\quad t\in [0,1) \\
1 , &\quad t=1
\end{array}
\right..$$
Then $\{\tilde{u}'_g\}_{g\in G}$ is an $\id_{C(I)}\otimes \alpha$-cocycle in $U(C_0([0,1))\otimes J)_0$ 
satisfying $\tilde{u}_g(0)=\tilde{u}'_g(0)=u_g$ for any $g\in G$, and 
$\|\tilde{u}_g-\tilde{u}'_g\|<\varepsilon$ for any $g\in S$. 
By Lemma \ref{small section},(2), we can get a homotopy fixed point $\tilde{F}:EG\to U(C_0([0,1))\otimes J)_0$ for $\tilde{u}'$ satisfying 
$\tilde{F}(x)(0)=\tilde{F}_0(x)(0)=F(x)$ for any $x\in EG$. 
Writing $\tilde{F}(x,t)=\tilde{F}(x)(t)$, we get the desired $\tilde{F}$. 

It remains to show the last statement. 
Assume that $A=B_\flat$ with a unital Kirchberg algebra $B$ and $\alpha$ comes from an outer action of $G$ on $B$. 
Then thanks to \cite[Lemma 6.2]{IMI}, we can modify $r(t)$ for $t$ sufficiently close to 1 so that it 
continuously extends to $[0,1]$. 
\end{proof}

\subsection{Second cohomology vanishing}
Let $(\alpha,u)$ be a cocycle action of a poly-$\Z$ group $G$ on a unital $C^*$-algebra $A$ in $\AC(\cO_\infty,\mu^G)$.  
Let $\{E_{g,h}\}_{g,h}\subset \K(\ell^2(G))$ be the canonical system of matrix units, and let 
$\{\rho_g\}_{g\in G}$ be the right regular representation of $G$. 
For simplicity we denote $f=E_{e,e}\in \K(\ell^2(G))$. 

We set 
$$W_g=(\sum_{s\in G}\alpha_s^{-1}(u(s,g)^{-1})\otimes E_{s,s})(1_A\otimes \rho_g)\in U(M(A\otimes \K(\ell^2(G)))).$$
Then 
$$W_g(\alpha_g\otimes \id)(W_h)(u(g,h)\otimes 1)W_{gh}^{-1}=1,$$
and we get a genuine action of $G$ on $A\otimes \K(\ell^2(G))$ given by 
$$\hat{\hat{\alpha}}_g=\Ad W_g\circ (\alpha_g\otimes\id).$$ 
In fact $\hat{\hat{\alpha}}$ is conjugate to the second dual action of $G$ on 
$$(A\rtimes_{(\alpha,u)}G)\rtimes_{\hat{\alpha}}\hat{G}\cong A\otimes \K(\ell^2(G)),$$
and its inner conjugacy class depends only on the equivalence class of $(\alpha,u)$.  

Note that $K_0(A)$ is a $G$-module through $\alpha$, and we denote by $K_0(A)^G$ the set of $G$-fixed points in $K_0(A)$.  
For $a\in K_0(A)^G$, we choose a non-zero projection $p\in A$ with $[p]=a$, and 
denote by $P^{a}$ the connected component of $p\otimes f$ in the projections in $A\otimes \K(\ell^2(G))$. 
Then $P^a$ does not depend on the choice of $p$, and it is a $G$-space through $\hat{\hat{\alpha}}$. 
We denote by $E^a$ the Borel construction $EG\times_G P^a$, which is a fiber bundle over $BG$.

\begin{theorem}\label{cocycle actions}
Let $(\alpha,u)$ be a cocycle action of a poly-$\Z$ group $G$ on a unital $C^*$-algebra $A$ in $\AC(\cO_\infty,\mu^G)$. 
Then $(\alpha,u)$ is equivalent to a genuine action if and only if the $G$-action $\hat{\hat{\alpha}}$ on $P^{[1_A]}$ 
has a homotopy fixed point. 
\end{theorem}

The following two lemmas hold in a much more general situation. 

\begin{lemma}\label{existence} 
If $(\alpha,u)$ is equivalent to an action, the $G$-action $\hat{\hat{\alpha}}$ on $P^{[1_A]}$ 
has a homotopy fixed point. 
\end{lemma}

\begin{proof} Assume that $(\alpha,u)$ is equivalent to an action. 
Since the inner conjugacy class of $\hat{\hat{\alpha}}$ depends only on the equivalence class of $(\alpha,u)$, 
we may and do assume that $u(g,h)=1$, and hence $W_g=1_A\otimes \rho_g$. 
Thanks to Lemma \ref{hfp for rho}, there exists a continuous map $U:EG\to U(M(\K(\ell^2(G))))$ satisfying 
$U(g\cdot x)=\rho_gU(x)$. 
Let $p(x)=1_A\otimes U(x)fU(x)^{-1}$. 
Then $p:EG\to P^{[1_A]}$ is a continuous map satisfying 
$$p(g\cdot x)=(1\otimes \rho_gU(x)fU(x)^{-1}\rho_g^{-1})=(\alpha_g\otimes \Ad \rho_g)(p(x)),$$
and $p$ is a homotopy fixed point. 
\end{proof}

\begin{lemma}\label{extension} 
If there exists a homotopy fixed point $p:EG\to P^{[1_A]}$, there exists a continuous map $v:G\times EG\to U(A)$ satisfying 
$$u(g,h)v(gh,x)=\alpha_g(v(h,x))v(g,h\cdot x),$$
for any $g,h\in G$ and $x\in EG$. 
\end{lemma}

\begin{proof} Since the map
$$U((A\otimes \K(\ell^2(G))))\ni U\mapsto U(1_A\otimes f)U^{-1}\in P^{[1_A]}$$
is a fiber bundle (see \cite[Lemma 2.8, Corollary 2.9]{DP-I}) and $EG$ is contractible, 
the homotopy lifting property implies that there exists a continuous map $U:EG\to U((A\otimes \K(\ell^2(G))))$ 
satisfying $p(x)=U(x)(1_A\otimes f)U(x)^{-1}$. 
Since $p(g\cdot x)=\hat{\hat{\alpha}}_g(p(x))$, we have 
$$U(g\cdot x)(1_A\otimes f)U(g\cdot x)^{-1}=\hat{\hat{\alpha}}_g(U(x))W_g(1_A\otimes f)W_g^{-1}\hat{\hat{\alpha}}_g(U(x)^{-1}),$$
and so there exists a continuous map 
$$V:G\times EG\to U(M(A\otimes (1-f)\K(\ell^2(G))(1-f))+A\otimes f))$$
satisfying $U(g\cdot x)=\hat{\hat{\alpha}}_g(U(x))W_gV(g,x)$. 

On one hand, we have 
$U(gh\cdot x)=\hat{\hat{\alpha}}_{gh}(U(x))W_{gh}V(gh,x)$, 
and on the other hand, 
\begin{align*}
\lefteqn{U(gh\cdot x)=\hat{\hat{\alpha}}_g(U(h\cdot x))W_gV(g,h\cdot x)} \\
 &=\hat{\hat{\alpha}}_g(\hat{\hat{\alpha}}_h(U(x))W_hV(h,x))W_gV(g,h\cdot x) \\
 &=\hat{\hat{\alpha}}_{gh}(U(x))W_g(\alpha_g\otimes \id)(W_hV(h,x))V(g,h\cdot x)\\
 &=\hat{\hat{\alpha}}_{gh}(U(x))W_{gh}(u(g,h)^{-1}\otimes 1)
 (\alpha_g\otimes \id)(V(h,x))V(g,h\cdot x).
\end{align*}
Thus we get 
$$u(g,h)\otimes 1=(\alpha_g\otimes \id)(V(h,x))V(g,h\cdot x)V(gh,x)^{-1}.$$ 
Now $v:G\times EG\to U(A)$ given by $V(g,x)(1_A\otimes f)=v(g,x)\otimes f$ satisfies the desired property. 
\end{proof}

\begin{proof}[Proof of Theorem \ref{cocycle actions}] It remains to show that the existence of $v$ as in the previous lemma implies that 
$(\alpha,u)$ is equivalent to an action. 
Let $\alpha'_g=\Ad v(g,0_n)^{-1}\circ \alpha_g$ and let $$u'(g,h)=v(g,0_n)^{-1}\alpha_g(v(h,0_n)^{-1})u(g,h)v(gh,0_n).$$ 
Then $(\alpha,u)$ is equivalent to $(\alpha',u')$ and 
\begin{align*}
\lefteqn{u'(g,h)=v(g,0_n)^{-1}\alpha_g(v(h,0_n)^{-1}v(h,x))v(g,h\cdot x)v(gh,x)^{-1}v(gh,0_n) } \\
 &=\alpha'_g(v(h,0_n)^{-1}v(h,x))v(g,0_n)^{-1}v(g,h\cdot x)v(gh,x)^{-1}v(gh,0_n)\\
 &=\alpha'_g(v'(h,x))v'(g,h\cdot x)v'(gh,x)^{-1},
\end{align*}
where $v'(g,x)=v(g,0_n)^{-1}v(g,x)$. 
Thus we have $v'(g,0_n)=1$, $u'(g,h)=v'(g,h\cdot 0_n)$, and 
$$u'(g,h)v'(gh,x)=\alpha'_g(v'(h,x))v'(g,h\cdot x).$$
Now Theorem \ref{vanishing2} implies the statement. 
\end{proof}

\begin{remark}\label{promap} 
Since the generalized mapping torus $M_{\hat{\hat{\alpha}}}=\Gamma(\cP_{\hat{\hat{\alpha}}}\times _G \sA)$ is identified with 
$$\{f\in C^b(EG,\sA);\; f(g\cdot x)=\hat{\hat{\alpha}}_g(f(x)),\;\forall g\in G,\;\forall x\in EG\},$$
the condition in Theorem \ref{cocycle actions} is equivalent to the following $K$-theoretical condition: 
there exists a projection $q\in M_{\hat{\hat{\alpha}}}$ satisfying $K_0(\mathrm{ev}_{x_0})([q])=[1_A]$ in $K_0(A)$, 
where $x_0\in BG$ is an arbitrary point and $\mathrm{ev}_{x_0}:M_{\hat{\hat{\alpha}}}\to \sA$ is the evaluation map at $x_0$. 
\end{remark}

Let $\alpha: G\to \Aut(A)$ be an action of $G$ on a unital $C^*$-algebra $A$ in $\AC(\cO_\infty,\mu^G)$,  
and let $p\in A$ be a non-zero projection with $[p]\in K_0(A)^G$. 
For each $g\in G$, we choose a partial isometry $u_g\in A$ satisfying $u_g^*u_g=\alpha_g(p)$ and 
$u_gu_g^*=p$. 
We set $w_{g,h}=u_g\alpha_g(u_h)u_{gh}^*\in U(pAp)$, and for $b\in pAp$ we set 
$\alpha^p_g(b)=u_g\alpha_g(b)u_g^*$. 
Then $(\alpha^p,w)$ is a cocycle action of $G$ on $pAp$, whose equivalence class does not depend on 
the choice of $\{u_g\}_{g\in G}$. 
Recall that we have $\hat{\hat{\alpha}}_g=\alpha_g\otimes \Ad\rho_g$ as $\alpha$ is a genuine action.  

\begin{theorem} Let $\alpha:G\to \Aut(A)$ be an action of a poly-$\Z$ group $G$ on a unital $C^*$-algebra $A$ in $\AC(\cO_\infty,\mu^G)$, 
and let $(\alpha^p,w)$ be the cocycle action obtained as above. 
Then $(\alpha^p,w)$ is equivalent to a genuine action if and only if the $G$-action $\hat{\hat{\alpha}}$ 
on $P^{[p]}$ has a homotopy fixed point. 
\end{theorem}

\begin{proof} Assume first that $(\alpha^p,w)$ is equivalent to a genuine action. 
Then we may assume that partial isometries $u_g\in A$ satisfy $u_gu_g^*=p$, $u_g^*u_g\in \alpha_g(p)$, and 
$u_g\alpha_g(u_h)=u_{gh}$. 
Since $\alpha^p$ is a genuine action with this choice of $\{u_g\}_{g\in G}$, 
the second dual action $\widehat{\widehat{\alpha^p}}$ of $\alpha^p$ on $pAp\otimes \K(\ell^2(G))$ is given by 
$$\widehat{\widehat{\alpha^p}}_g(X)=(u_g\otimes 1)\hat{\hat{\alpha}}_g(X)(u_g^*\otimes 1)=R\hat{\hat{\alpha}}_g(R^*XR)R^*,$$
where 
$$R=\sum_{g\in G}u_{g^{-1}}\otimes E_{g,g}\in M(A\otimes \K(\ell^2(G))),$$
is a partial isometry with $RR^*=p\otimes 1$, and 
$$R^*R=\sum_{g\in G}\hat{\hat{\alpha}}_g(p\otimes f),$$
and convergence is in the strict topology. 
Thanks to Lemma \ref{existence}, there exists a continuous map 
$q_0:EG\to P^{[p]}$ satisfying $q_0(g\cdot x)=\widehat{\widehat{\alpha^p}}_g(q_0(x))$ for any $x\in EG$ and $g\in G$. 
Letting $q(x)=R^*q_0(x)R$, we get a homotopy fixed point $q:EG\to P^{[p]}$ for $\hat{\hat{\alpha}}$. 

Assume on the other hand that there exists a homotopy fixed point $q:EG\to P^{[p]}$ for $\hat{\hat{\alpha}}$. 
Then as in the proof of Lemma \ref{extension}, there exists a norm continuous map $U:EG\to U(M(A \otimes \K(\ell^2(G))))$ 
satisfying $q(x)=U(x)(p\otimes f) U(x)^{-1}$. 
As before there exists a continuous map $S:G\times EG\to U(M(A\otimes \K(\ell^2(G))))$ satisfying 
$$U(g\cdot x)=(1\otimes \rho_g)(\alpha_g\otimes \id)(U(x))S(g,x),$$
and $S(g,x)(p\otimes f)S(g,x)^*=\alpha_g(p)\otimes f$. 
Moreover, 
\begin{align*}
\lefteqn{U(gh\cdot x)=(1\otimes \rho_g)(\alpha_g\otimes \id)(U(h\cdot x))S(g,h\cdot x)} \\
 &=(1\otimes \rho_g)(\alpha_g\otimes \id)((1\otimes \rho_h)(\alpha_h\otimes \id)(U(x))S(h,x))S(g,h\cdot x)\\
 &=(1\otimes \rho_{gh})(\alpha_{gh}\otimes \id)(U(x))(\alpha_g\otimes \id)(S(h,x))S(g,h\cdot x),
\end{align*}
implies $S(gh,x)=(\alpha_g\otimes \id)(S(h,x))S(g,h\cdot x)$. 
We define $s(g,x)\in \alpha_g(p)Ap$ by 
$S(g,x)(p\otimes f)=s(g,x)\otimes f$. 
Then each $s(g,x)$ is a partial isometry with $s(g,x)^*s(g,x)=p$, $s(g,x)s(g,x)^*=\alpha_g(p)$, and 
we have $s(gh,x)=\alpha_g(s(h,x))s(g,h\cdot x)$. 
We may assume $u_g=s(g,0)^*$, and 
\begin{align*}
\lefteqn{w_{g,h}=s(g,0)^*\alpha_g(s(h,0)^*)s(gh,0)} \\
 &=s(g,0)^*\alpha_g(s(h,0)^*)\alpha_g(s(h,0))s(g,h\cdot 0)
=s(g,0)^*s(g,h\cdot 0). 
\end{align*}
Let $v(g,x)=s(g,0)^*s(g,x)\in U(pAp)_0$. 
Then we have $w_{g,h}=v(g,h\cdot 0)$, and 
$$w_{g,h}v(gh,x)=\alpha^p_g(v(h,x))v(g,h\cdot x)$$
for any $g,h\in G$ and $x\in EG$. 
Thus Theorem \ref{vanishing2} implies that $(\alpha^p,w)$ is equivalent to an action. 
\end{proof}


\section{Proof of Theorem \ref{deformation}}

Throughout this section, we freely use the notation in Introduction. 
We assume that $G$ is a poly-$\Z$ group, and for $G$ we use the convention in Section 2.  
We denote by $\{E_{g,h}\}_{g,h\in G}$ the canonical system of matrix units in $\K(\ell^2(G))$. 
We assume that $A$ is a unital Kirchberg algebra.

\begin{definition}\label{isobundle} Let $\alpha$ and $\beta$ be two $G$-actions on $A$ 
satisfying $KK(\alpha_g)=KK(\beta_g)$ for any $g\in G$. 
We introduce a fiber bundle $\cI_{\salpha,\sbeta}^0$ over $BG$ with fiber $\Aut(\sA)_0$ by the Borel construction 
$p: EG\times_G \Aut(\sA)\to BG$  with a $G$-action on $\Aut(\sA)_0$ given by   
$g\cdot\gamma=\sbeta_g\circ \gamma\circ {\salpha_g}^{-1}$. 
\end{definition}

The homotopy fixed point $\Phi:EG\to \Aut(\sA)_0$ in Theorem \ref{deformation} is identified with the continuous 
section $[x]\mapsto [(x,\Phi(x))]$ of $p$. 

\begin{remark} Thanks to our choice of $EG$, the whole information of the homotopy fixed point $\Phi$ is encoded 
in its restriction to $EN\times I$, and we often regard $\Phi$ as a map from $EN\times I$ to $\Aut(\sA)$. 
\end{remark}

The statement of Theorem \ref{deformation} is trivially true if the Hirsch length of $G$ is 0. 
As induction hypothesis, we assume that it is true for $N$. 

\subsection{Step I}
\begin{lemma}\label{homotopy} Theorem \ref{deformation} is true if it is the case for $\alpha$, $\beta$, and $\Phi$ with the following properties. 
There exists an $\alpha|_N$-cocycle $\{a_h\}_{h\in N}$ with $\beta_h=\Ad a_h\circ \alpha_h$ for any $h\in N$, and 
$\Phi(x,0)=\Ad Z$ for any $x\in EN$, where 
$$Z=\sum_{h\in N,\;n\in \Z}a_{h^{-1}}\otimes E_{\xi^nh,\xi^nh}\otimes 1.$$ 
\end{lemma}

\begin{proof}
Let $\alpha$ and $\beta$ be outer actions of $G$ on a 
unital Kirchberg algebra $A$. 
Assume that there exists a continuous map $\Phi:EG\rightarrow \Aut(\sA)_0$ satisfying 
$\Phi(g\cdot z)=\sbeta_g\circ \Phi(z)\circ{\salpha_g}^{-1}$ for any 
$z\in EG$ and $g\in G$. 

For $x\in EN$, let $\Psi(x)=\Phi(x,0)$. 
Then by the induction hypothesis, if we replace $(A,\alpha,\beta)$ with 
$(\cO_\infty\otimes A,\mu\otimes \alpha,\mu\otimes \beta)$, 
we may assume that there exist an $\alpha|_N$-cocycle $\{a_h\}_{h\in N}$ in $U(A)$, 
a continuous path $\{\gamma_t\}_{t\in I}$ 
in $\Aut(A)$, a continuous map $\tilde{\Psi}:EN\times I\rightarrow \Aut(\sA)_0$, 
satisfying 
\begin{itemize}
\item $\gamma_0=\id$ and $\gamma_1\circ \beta_h\circ \gamma_1^{-1}=\Ad a_h\circ \alpha_h$ for any $h\in N$, 
\item 
$$\tPsi(h\cdot x,t)=(\gamma_t\otimes \id)\circ \sbeta_h\circ (\gamma_t^{-1}\otimes \id)
\circ \tPsi(x,t)\circ {\salpha_h}^{-1}$$ 
for any $x\in EN$, $t\in I$, and $h\in N$, 
\item $\tPsi(x,0)=\Psi(x)$ and $\tPsi(x,1)=\Ad Z$ for any $x\in EN$, where 
$$Z=\sum_{h\in N,\;n\in \Z}a_{h^{-1}}\otimes E_{\xi^nh,\xi^nh}\otimes 1.$$ 
\end{itemize}

We claim that there exists a continuous map $\tPhi:EG\times I\to \Aut(\sA)_0$ satisfying
\begin{itemize}
\item $\tPhi(g\cdot z,t)=(\gamma_t\otimes \id)\circ \sbeta_g\circ (\gamma_t^{-1}\otimes \id)
\circ \tPhi(z,t)\circ {\salpha_g}^{-1}$ 
for any $z\in EG$, $t\in I$, and $g\in G$, 
\item $\tPhi(z,0)=\Phi(z)$ and $\tPhi(x,0,t)=\tPsi(x,t)$ for any $z\in EG$, $x\in EN$ and $t\in I$. 
\end{itemize} 

Since $\cI^0_{\salpha,{\gamma}_t\circ \sbeta\circ \gamma_t^{-1}}\to BG$ is a homotopy 
of fiber bundles, the continuous section of $\cI^0_{\salpha,\sbeta}\to BG$ given by $\Phi$ 
extends to a continuous family of sections of $\cI^0_{\salpha,\gamma_t\circ \sbeta\circ \gamma_t^{-1}}\to BG$ 
for $t\in [0,1]$, and we get a continuous map 
$\tPhi_1 :EG\times I \to \Aut(\sA)_0$ satisfying $\tPhi_1(z,0)=\Phi(z)$ for all $z\in EG$ and 
$$\tPhi_1(g\cdot z,t)=(\gamma_t\otimes \id)\circ \sbeta_g\circ (\gamma_t^{-1}\otimes \id)
\circ \tPhi_1(z,t)\circ {\salpha_g}^{-1}$$ 
for any $z\in EG$, $t\in I$, and $g\in G$. 

Let $\omega(x,t)=\tPhi_1(x,0,t)^{-1}\circ \tPsi(x,t)$. 
Then $\omega:EN\times I\to \Aut(\sA)_0$ is a continuous map satisfying $\omega(x,0)=\id$, and 
$$\omega(h\cdot x,t)=\salpha_h\circ \omega(x,t)\circ {\salpha_h}^{-1},$$
for any $h\in N$. 
Thus we have a homotopy of partial sections $r_t:BN\to \cI^0_{\salpha,\salpha}$ for $t\in I$ with 
$r_0([x])=[(x,\id)]$.  
On the other hand, we have a section $s_0:BG\ni [z]\mapsto [(z,\id)]\in \cI^0_{\salpha,\salpha}$. 
Since $(BG,BN)$ is an NDR pair, Theorem \ref{lifting extension} implies that 
the homotopy $r_t$ extends to a homotopy of sections 
$s_t:BG\to \cI^0_{\salpha,\salpha}$ for $t\in I$. 
Thus we have a continuous map $\omega_1:EG\times I\to \Aut(\sA)_0$ satisfying 
$\omega_1(x,0,t)=\omega(x,t)$ for any $x\in EN$ and $t\in I$, $\omega_1(z,0)=\id$ for any $z\in EG$, and 
$$\omega_1(g\cdot z,t)=\salpha_g\circ \omega(x,t)\circ {\salpha_g}^{-1},$$
for any $g\in G$. 
Now $\tPhi(z,t)=\tPhi_1(z,t)\omega_1(z,t)$ has the desired property. 

By replacing $\beta$ with $\gamma_1\circ \beta\circ\gamma_1^{-1}$, and $\Phi$ with $\tPhi(\cdot,1)$, 
we get the statement.  
\end{proof}
\subsection{Step II}
The above lemma shows that from the beginning, we may and do assume that 
$\Phi(x,0)=\Ad Z$ and $\beta_h=\Ad a_h\circ \alpha_h$ for any $x\in BN$ and $h\in N$. 
Since $a_h\otimes 1=Z\salpha_h(Z^{-1})$, we also have 
$$\sbeta_h=\Ad Z\circ \salpha_h\circ \Ad Z^*.$$

Since $\pi_0(\mathrm{Inn}(M(\sA)))=\{0\}$ and $\Phi(0,0)=\Ad Z$, 
$\Phi(0,1)=\sbeta_\xi\circ \Ad Z\circ {\salpha_\xi}^{-1}$,
the path $\{\Phi(0,y)\}_{y\in I}$ determines a homotopy class of a path from $\id$ to 
$\sbeta_\xi\circ {\salpha_\xi}^{-1}$ in $\Aut(\sA)_0$.  
On the other hand, there exists a unitary $r\in U(A^\flat)$ satisfying 
$$\beta_\xi=\lim_{y\to 1-0}\Ad r(y)\circ \alpha_\xi,$$ 
which also determines a homotopy class of a path from $\id$ to $\sbeta_\xi\circ {\salpha_\xi}^{-1}$. 
Since 
$$\Pi_{A,0}:\pi_0(U(A_\flat))\to \pi_0(\Omega\Aut(\sA))=\pi_1(\Aut(\sA)_0)$$ 
is an isomorphism, we can arrange $r$ so that the two homotopy classes coincide.  

Now $\Ad a_h\circ \alpha_h|_{A_\flat}=\beta_h|_{A_\flat}$  and $\Ad r\circ \alpha_\xi|_{A_\flat}$ 
give rise to a cocycle action $(\gamma,\omega)$ of $G$ on $A_\flat$. 
In order to obtain the statement appealing to Theorem \ref{KKtrivcc}, 
our primary concern is whether it is equivalent to a genuine action or not. 
By construction, the 2-cocycle $\omega$ is trivial on $N\times N$. 
Thus $\omega$ is determined by a $\beta|_N$-1-cocycle in $U(A_\flat)$, given by 
$$u_h=r\alpha_\xi(a_{h^\xi})\alpha_h(r^{-1})a_h^{-1}= r\alpha_\xi(a_{h^\xi})a_h^{-1}\beta_h(r^{-1})\in U(A_\flat)$$
as we have
$$\Ad r\circ \alpha_\xi \circ \Ad a_{h^\xi}\circ \alpha_{h^\xi}\circ (\Ad r\circ \alpha_\xi)^{-1}=
\Ad (r\alpha_\xi(a_{h^\xi})\alpha_h(r^{-1})a_h^{-1})\circ \Ad a_h\circ \alpha_h. $$  
The cocycle action of $G$ is equivalent to an action if and only if the $\beta|_N$-cocycle 
$\{u_h\}_{h\in H}$ is a coboundary in $U(A_\flat)$.  

\begin{lemma} The unitary $u_h$ belongs to $U(A_\flat)_0$.
\end{lemma}

\begin{proof} 
To show $u_h\in U(A_\flat)_0$, we compute the homotopy class of the loop in $\Aut(A^s)_0$ determined by 
$\Ad (u_h(t)\otimes 1)$ (note that it is a continuous path from $\Ad (u_h(0)\otimes 1)$ to $\id$). 
Using $a_h\otimes 1=Z\salpha_{h}(Z^{-1})$, 
We get
\begin{align*}
\lefteqn{u_h(y)\otimes 1=(r(y)\otimes 1)\salpha_\xi(a_{h^\xi}\otimes 1)\salpha_h(r(y)^{-1}\otimes 1) (a_h^{-1}\otimes 1)} \\
 &=(r(y)\otimes 1)\salpha_\xi(Z)\salpha_h(\salpha_\xi(Z^{-1})(r(y)^{-1}\otimes 1))(a_h^{-1}\otimes 1), 
\end{align*}
and 
\begin{align*}
\lefteqn{\Ad (u_h(y)\otimes 1)} \\
 &=\Ad((r(y)\otimes 1)\salpha_\xi(Z))\circ \Ad \salpha_h(\salpha_\xi(Z^{-1})(r(y)^{-1}\otimes 1)) 
 \circ \Ad (a_h^{-1}\otimes 1)\\
 &=\Ad((r(y)\otimes 1)\salpha_\xi(Z))\circ \salpha_h\circ \Ad((r(y)\otimes 1)\salpha_\xi(Z))^{-1}
 \circ {\sbeta_h}^{-1}. 
\end{align*}
Note that we have 
$$\lim_{y\to1-0}\Ad ((r(y)\otimes 1)\salpha_\xi(Z))=\sbeta_\xi\circ {\salpha_\xi}^{-1}\circ \Ad \salpha_\xi(Z)
=\sbeta_\xi \circ \Ad Z\circ {\salpha_\xi}^{-1}=\Phi(0,1).$$
To compute the homotopy class of the loop in $\Aut(\sA)_0$ given by $\Ad(u_h(y)\otimes 1)$, 
we can replace 
$\Ad ((r(y)\otimes 1)\salpha_\xi(Z))$ in the above with $\Phi(0,y)$, $y\in I$, and we get 
$$\Phi(0,y)\circ \salpha_h\circ \Phi(0,y)^{-1}\circ {\sbeta_h}^{-1}=\Phi(0,y)\circ  \Phi(h\cdot 0,y)^{-1}.$$
Since $\Phi(x,0)=\Ad Z$ and $\Phi(x,1)=\sbeta_\xi\circ \Ad Z\circ {\salpha_\xi}^{-1}$, we have 
$$(\Phi(0,y)\circ \Phi(x,y)^{-1})_{y\in I}\in \Omega\Aut(\sA),
\quad \forall x\in EN,$$ 
and 
$\{(\Phi(0,y)\circ \Phi(th\cdot 0,y)^{-1})_{y\in I}\}_{t\in I}$ 
gives a homotopy from the trivial loop to the loop $(\Phi(0,y)\circ \Phi(h\cdot 0,y)^{-1})_{y\in I}$. 
\end{proof}

We take a norm continuous map $R:[0,1)\to U(M(\sA))$ satisfying $R(y)=r(y)\otimes 1$ for any $y\in[1/2,1)$ and 
$R(0)=Z\salpha_\xi(Z^{-1})$. 
Then $\Ad (R(0)\salpha_\xi(Z))=\Ad Z=\Phi(0_n)$, and 
$$\lim_{y\to 1-0}\Ad (R(y)\salpha_\xi(Z))=\Phi(0_{n-1},1).$$
We have seen in the proof of the previous lemma that the two paths $\{\Phi(0_{n-1},y)\}_{y\in I}$ 
and $\{\Ad (R(y)\salpha_\xi(Z))\}_{y\in [0,1)}$, continuously extended to $[0,1]$ by setting $\Phi(0,1)$ at $y=1$, 
are homotopic. 
Thus thanks to Theorem \ref{lifting extension} and Lemma \ref{NDR}, we may assume 
$$\Phi(0_{n-1},y)=\Ad (R(y)\salpha_\xi(Z))$$ for any $y\in [0,1)$ keeping the condition $\Phi(x,0)=\Ad Z$ for any $x\in EN$. 
Indeed, recall that $\Phi$ is identified with a continuous section $s$ of the fiber bundle $\cI^0_{\salpha,\sbeta}\to BG$. 
Since $(BG, BN\cup[\{0_{n-1}\}\times I])$ is an NDR pair, any homotopy of sections on $BN\cup[\{0_{n-1}\}\times I]$ starting 
from the restriction of $s$ to $BN\cup[\{0_{n-1}\}\times I]$ extends to $BG$.

Since $a_h\otimes 1=\sbeta_h(Z^{-1})Z$, we have 
\begin{align*}
\lefteqn{u_h(y)\otimes 1} \\
 &=(r(y)\otimes 1)\salpha_\xi(Z)\salpha_h(\salpha_\xi(Z^{-1})(r(y)^{-1}\otimes 1))(a_h^{-1}\otimes 1) \\
 &=(r(y)\otimes 1)\salpha_\xi(Z)(a_h^{-1}\otimes 1)\sbeta_h(\salpha_\xi(Z^{-1})(r(y)^{-1}\otimes 1))\\
 &=(r(y)\otimes 1)\salpha_\xi(Z)Z^{-1}\sbeta_h((r(y)\otimes 1)\salpha_\xi(Z)Z^{-1})^{-1}.
 \end{align*}
We set 
$$U_h(y)=R(y)\salpha_\xi(Z)Z^{-1}\sbeta_h(R(y)\salpha_\xi(Z)Z^{-1})^{-1}.$$
Then $\{U_h(y)\}_{h\in N,\;y\in [0,1)}$ is a continuous family of $\sbeta$-cocycles 
with $U_g(0)=1$ and $U_h(y)=u_h(y)\otimes 1$ for any $h\in N$ and $y\in [1/2,1)$. 
By definition, 
\begin{align*}
\lefteqn{\Ad U_h(y)=\Ad(R(y)\salpha_\xi(Z)Z^{-1}\sbeta_h(R(y)\salpha_\xi(Z)Z^{-1})^{-1})} \\
 &=\Phi(0_{n-1},y)\circ \Ad Z^{-1}\circ \sbeta_h\circ \Ad Z\circ \Phi(0_{n-1},y)^{-1}\circ {\sbeta_h}^{-1} \\
 &=\Phi(0_{n-1},y)\circ \salpha_h\circ \Phi(0_{n-1},y)^{-1}\circ {\sbeta_h}^{-1}\\
 &=\Phi(0_{n-1},y)\circ \Phi(h\cdot 0_{n-1},y)^{-1}. 
\end{align*}

Note that 
$$(\Phi(0_{n-1},y)\circ \salpha_h\circ \Phi(0_{n-1},y)^{-1})_{y\in I}
=(\Phi(0_{n-1},y)\circ \Phi(h\cdot 0_{n-1},y)^{-1}\circ \sbeta_h)_{y\in I}$$  
is a continuous family of $N$-actions with $\sbeta|N$ at the two end points, 
and we can introduce an $N$-action on $\Omega\Aut(\sA)_0$ by 
$$(h\cdot\theta)(y)= \Phi(0_{n-1},y)\circ \Phi(h\cdot 0_{n-1},y)^{-1}\circ \sbeta_h\circ\theta(y)\circ {\sbeta_h}^{-1}.$$
For $y\in [0,1)$, it is $\Ad U_h(y)\circ \sbeta_h\circ\theta(y)\circ {\sbeta_h}^{-1}$. 

\begin{lemma} Let $\Theta(x)(y)=\Phi(0_{n-1},y)\circ \Phi(x,y)^{-1}$ for $x\in EN$ and $y\in I$. 
Then $\Theta:EN\to \Omega\Aut(\sA)$ is a homotopy fixed point for the $N$-action on $\Omega\Aut(\sA)_0$ 
satisfying $\Theta(0_{n-1})=\id$.   
\end{lemma}

\begin{proof} For $y=0$, we have 
$$\Theta(x)(0)=\Phi(0_{n-1},0)\circ \Phi(x,0)^{-1}=\Ad Z\circ \Ad Z^{-1}=\id.$$
For $y=1$, 
$$\Theta(x)(1)=\Phi(0_{n-1},1)\circ \Phi(x,1)^{-1}=\id,$$
and $\Theta(x)\in \Omega\Aut(\sA)$. 

For $h\in N$, $x\in EN$, and $y\in [0,1)$,  
\begin{align*}
\lefteqn{\Theta(h\cdot x)(y)=\Phi(0_{n-1},y)\circ \Phi(h\cdot x,y)^{-1}} \\
 &=\sbeta_h\circ \Phi(h^{-1}\cdot 0_{n-1},y)\circ \Phi(x,y)^{-1}\circ {\sbeta_h}^{-1}\\
 &=\sbeta_h\circ \Phi(h^{-1}\cdot 0_{n-1},y)\circ \Phi(0_{n-1},y)^{-1}\circ \Theta(x)(y)\circ {\sbeta_h}^{-1}\\
 &=\Phi(0_{n-1},y)\circ \Phi(h\cdot 0_{n-1},y)^{-1}\circ \sbeta_h\circ \Theta(x)(y)\circ {\sbeta_h}^{-1},
\end{align*}
and $\Theta$ is a homotopy fixed point. 
\end{proof}

\begin{lemma}\label{hfp} There exist a homotopy fixed point $f:EN\to U(A_\flat)_0$ for $u$ with $f(0_{n-1})=1$, a continuous map 
$F:EN\times [0,1)\to U(M(\sA))$, and a continuous map $\tilde{\Theta}:EN\times I\to \Omega\Aut(\sA)_0$ 
satisfying 
\begin{itemize} 
\item $F(h\cdot x,y)=U_h(y)\sbeta_h(F(x,y))$ for any $h\in N$, $x\in EN$, $y\in [0,1)$. 

\item $F(x,0)=F(0_{n-1},y)=1$ for any $x\in EN$, $y\in [0,1)$, 
and $F(x,y)=\hat{f}(x,y)\otimes 1$ for any $x\in EN$, $y\in [1/2,1)$, 
where $\hat{f}(x,\cdot)$ represents $f(x)$.   

\item $\{\tilde{\Theta}(\cdot,t)\}_{t\in I}$ is a continuous family of homotopy fixed points 
of the $N$-action on $\Omega\Aut(\sA)_0$. 

\item $\tilde{\Theta}(x,0)=\Theta(x)$ for any $x\in EN$, and $\tilde{\Theta}(x,1)(y)=\Ad F(x,y)$ 
for any $x\in EN$, $y\in [0,1)$.
\end{itemize}  
\end{lemma}

\begin{proof} Recall that $K$ is a finite subcomplex of $EN$ satisfying $N\cdot K=EN$. 
Let $\{e_i\}_{i=0}^m$ be the set of cells in $K$. 
We may assume that $e_0=\{0_{n-1}\}$, and $\dim e_i\leq \dim e_j$ whenever $i<j$. 
We set $K_k=\cup_{i=0}^k e_i$, which is a subcomplex of $K$. 
We inductively construct continuous maps $f_i:K_i\to U(A_\flat)_0$, $F_i:K_i\times [0,1)\to U(M(\sA))$, 
$\tilde{\Theta}_i:EN\times I\to \Omega\Aut(\sA)_0$ for $i=0,1,\ldots,m$ satisfying the following conditions: 
\begin{itemize}
\item $f_0(0_{n-1})=1$, $F_0(0_{n-1},y)=1$ for any $y\in [0,1)$, and $\tilde{\Theta}_0(x,t)=\Theta(x)$ for any $x\in EN$. 

\item $f_i(h\cdot x)=u_h\beta_h(f_i(x))$ for any $h\in N$ and $x\in K_i$ with $h\cdot x\in K_i$. 
The map $f_{i+1}$ is an extension of $f_i$. 

\item $F_i(h\cdot x,y)=U_h(y)\sbeta_h(F_i(x,y))$ for any $h\in N$, $x\in K_i$ with $h\cdot x\in K_i$, 
and $y\in [0,1)$. 
The map $F_{i+1}$ is an extension of $F_i$.

\item  $F_i(x,0)=1$ for any $x\in K_i$, and $F_i(x,y)=\hat{f_i}(x,y)\otimes 1$ for any $x\in K_i$, $y\in [1/2,1)$, 
where $\hat{f}_i(x,\cdot)$ represents $f_i(x)$. 

\item $\{\tilde{\Theta}_i(\cdot,t)\}_{t\in I}$ is a continuous family of homotopy fixed points 
of the $N$-action on $\Omega\Aut(\sA)_0$. 

\item $\tilde{\Theta}_i(\cdot,1)=\tilde{\Theta}_{i+1}(\cdot,0)$ for $i=0,1,\ldots,k-1$.  

\item $\Ad F_i(x,y)=\tilde{\Theta}_i(x,1)(y)$ for any $x\in K_i$, $y\in [0,1)$, and $i=0,1,\ldots,m$. 
\end{itemize}

Since $e_0=\{0_{n-1}\}$ and $\Theta(0_{n-1})=\id$, we can construct $f_0$, $F_0$, and 
$\tilde{\Theta}_0$ as desired. 
Assume we have constructed $f_i$, $F_i$, and $\tilde{\Theta}_i$ for $0\leq i\leq k$. 
If there exists $i<k+1$ and $h\in N$ satisfying $h\cdot e_i=e_{k+1}$, 
we put $f_{k+1}(x)=u_{h}\beta_h(f_i(h^{-1}\cdot x))$, $F_{k+1}(x,y)=U_h(y)\sbeta_h(F_i(h^{-1}\cdot x,y))$ and 
$\tilde{\Theta}_{k+1}(x,t)=\tilde{\Theta}_k(x,1)$ for $x\in e_{k+1}$. 

Assume that there is no $h\in N\setminus \{e\}$ with $(h\cdot K_k)\cap e_{k+1}\neq \emptyset$. 
We denote $l=\dim e_{k+1}$ for simplicity. 
Since $\tilde{\Theta}_k(\cdot,1)$ restricted to $\overline{e_{k+1}}$ is a continuous extension of $\Ad F_k$ restricted to 
$\partial e_{k+1}$, Theorem \ref{weak homotopy} implies that $f_k|_{\partial e_{k+1}}$ is homotopic to a constant map, and 
we can continuously extend $f_k$ to $f'_{k+1}:K_{k+1}\to U(A_\flat)_0$. 
We choose a continuous map $\hat{f}'_{k+1}:K_{k+1}\times [0,1)\to U(A)_0$ representing $f'_{k+1}$. 
We claim that we can choose $f'_{k+1}$ and $\hat{f}'_{k+1}$ so that $\hat{f}'_{k+1}$ is an extension of $\hat{f}_k$. 
Indeed, since 
$$\lim_{y\to 1-0}\sup_{x\in K_k}\|\hat{f}_k(x,y)-\hat{f}_{k+1}(x,y)\|=0,$$
we may assume $\|\hat{f}_k(x,y)-\hat{f}_{k+1}(x,y)\|< 1/2$ for any $x\in K_k$ and $y\in [0,1)$. 
Let $\varphi:D^l\to \overline{e_k}$ be the characteristic map of $e_k$. 
Recall that since our $N$-CW-complex structure of $EN$ comes from a triangulation of $BN$, it is a homeomorphism. 
We set $\psi(x)=\varphi(\frac{\varphi^{-1}(x)}{|\varphi^{-1}(x)|})$ for $x\in e_{k+1}\setminus \{\varphi(0_l)\}$. 
Replacing $\hat{f}'_{k+1}$ with 
$$\left\{
\begin{array}{ll}
\hat{f}'_{k+1}(\varphi(0_{l}),y) , &\quad x=\varphi(0_{l}) \\
 e^{|\varphi^{-1}(x)|\log (\hat{f}'_{k+1}(\psi(x),y)
\hat{f}'_{k+1}(\psi(x),y)^{-1})}
\hat{f}'_{k+1}(x,y), &\quad x\in e_{k+1}\setminus\{\varphi(0_{l})\} \\
\hat{f}_k(x,y) , &\quad x\in \partial e_k
\end{array}
\right.,
$$ 
we get the claim. 

Using the fact that $U(M(\sA))$ is contractible, we can choose a continuous map $F'_{k+1}:K_{k+1}\times [0,1)$, which is an 
extension of $F_k$, satisfying $F'_k(x,y)=\hat{f}'_{k+1}(x,y)\otimes 1$ for any $x\in K_{k+1}\times [1/2,1)$ and $F'_{k+1}(x,0)=1$ 
for any $x\in K_{k+1}$. 
Now $\Ad F'_{k+1}(x,\cdot)^{-1}\circ \tilde{\Theta}_k(x,1)$ gives a continuous map from $K_{k+1}$ to $\Omega\Aut(\sA)$ whose 
restriction to $K_k$ is a constant map $\id$. 
Since $K_{k+1}/K_k=\overline{e_{k+1}}/\partial e_{k+1}\cong S^l$, we may regard this as an element in 
$\Map(S^l,\Omega\Aut(\sA))_0$, giving a homotopy class $a\in \pi_l(\Omega\Aut(\sA))$. 
Since $\Pi_{A,l}:\pi_l(U(A)_\flat)\to \pi_l(\Omega\Aut(\sA))$ is an isomorphism, 
there exists a continuous map $f''_{k+1}:K_{k+1}\to U(A_\flat)_0$ with $f''_{k+1}(x)=1$ for $x\in K_k$ such that 
if we regard $f''_{k+1}$ as an element of $\Map(S^l,U(A_\flat))_0$, we have $\Pi_{A,l}([f''_{k+1}])=a$. 
Setting $f_{k+1}(x)=f'_{k+1}(x)f''_{k+1}(x)$, and repeating the same argument as above, we get desired maps 
$\hat{f}_{k+1}$ and $F_{k+1}$ so that $\Ad F_{k+1}(x,\cdot)^{-1}\circ \tilde{\Theta}_k(x,1)$ regarded as an element 
of $\Map(S^l,\Omega\Aut(\sA))_0$ is null-homotopic. 
Thus we can make a desired homotopy $\tilde{\Theta}_{k+1}$ on $K_{k+1}$. 
Since $q_N(K_{k+1})$ is a subcomplex of $BN$, we can apply Theorem \ref{lifting extension} to get $\tilde{\Theta}_{k+1}$ on the whole $EN$. 

Now it is easy to show that $f_m$, $\hat{f}_m$ and $F_m$ extend to $EN$ with desired properties, and we get $f$, $\hat{f}$ and $F$. 
Concatenating the homotopies $\tilde{\Theta}_k$, $k=0,1,\cdots,m$, we get $\tilde{\Theta}$. 
\end{proof}

Thanks to Theorem \ref{vanishing1}, the cocycle $\{u_h\}_{h\in N}$ is a coboundary within $U(A_\flat)_0$, 
and there exists $v\in U(A_\flat)_0$ satisfying $v^{-1}u_h\beta_h(v)=1$. 
By replacing $r$ with $v^{-1}r$, we get $r\alpha_\xi(a_{h^\xi})\alpha_h(r^{-1})a_h^{-1}=1$ in $U(A^\flat)$. 

Let $x_h(t)=a_h$ for $h\in N$ and let $x_\xi(t)=r(t)$ for $t\in [0,1)$. 
We can define $x_g(t)$ for the other $g\in G$ in an appropriate way so that 
$\{x_g(t)\}_{g\in G,\;t\in [0,1)}$ forms a $\alpha$-cocycle in $U(A^\flat)$ satisfying 
$$\lim_{t\to 1-0}\Ad x_g(t)\circ \alpha_g=\beta_g.$$ 
Theorem \ref{KKtrivcc} shows that if we replace $(A,\alpha,\beta)$ with $(\cO_\infty\otimes A,\mu\otimes \alpha,\mu\otimes \beta)$, 
there exist a continuous map $w :[0,1)\to U(A)$, $\gamma\in \Aut(A)$, and 
a $\alpha$-cocycle $\{b_g\}_{g\in G}$ in $U(A)$ satisfying 
$$\gamma\circ \beta_g\circ \gamma^{-1}=\Ad b_g\circ \alpha_g,\quad \forall g\in G,$$
$$\lim_{t\to1-0}\Ad w(t)=\gamma,$$
$$\lim_{t\to1-0}w(t)x_g(t)\alpha_g(w(t)^*)=b_g,\quad \forall g\in G.$$
We may and do assume $w(0)=1$. 
We set 
$$\gamma_t=\left\{
\begin{array}{ll}
\Ad w(t) , &\quad t\in [0,1)  \\
\gamma , &\quad t=1
\end{array}
\right.,
$$
$$b_g(t)=\left\{
\begin{array}{ll}
w(t)x_g(t)\alpha_g(w(t)^*) , &\quad t\in [0,1) \\
b_g , &\quad t=1
\end{array}
\right..
$$
Note that $\{b_h(t)\}_{h\in N}$ is a continuous family of $\alpha|_N$-cocycle, and we set 
$$Z_t=\sum_{h\in N,\; n\in \Z}b_{h^{-1}}(t)\otimes E_{\xi^nh,\xi^nh}\otimes 1.$$
Then $\{Z_t\}_{t\in I}$ is a strictly continuous family of unitaries in 
$U(M(A\otimes \K(\ell^2(G)))\otimes 1)\subset U(M(\sA))$ satisfying $Z_t\salpha_h(Z_t^*)=b_h(t)\otimes 1$ 
for any $h\in N$ and $t\in I$. 
We set 
$$W=\sum_{g\in G}b_{g^{-1}}\otimes E_{g,g}\otimes 1,$$
which is a unitary in $U(M(A\otimes\K(\ell^2(G)))\otimes 1)\subset U(M(\sA))$ satisfying 
$W\salpha_g(W^*)=b_g\otimes 1$ for $g\in G$. 

Let $\tPsi(x,t)=\Ad Z_t$ for $(x,t)\in EN\times I$. 
Then $\tPsi:EN\times I\to \Aut(\sA)_0$ is a continuous map satisfying $\tPsi(x,0)=\Ad Z=\Phi(x,0)$ and 
$$\tPsi(h\cdot x,t)=(\gamma_t\otimes \id)\circ \sbeta_h\circ (\gamma_t^{-1}\otimes \id)
\circ \tPsi(x,t)\circ {\salpha_h}^{-1},$$
for any $x\in EN$, $t\in I$, and $h\in N$. 
As in Lemma \ref{homotopy}, we can continuously extend $\tPsi$ to $\tPhi:BG\times I\to \Aut(\sA)_0$ satisfying 
$\tPhi(x,0,t)=\tPsi(x,t)=\Ad Z_t$ for any $(x,t)\in EN\times I$, and 
$$\tPhi(g\cdot z,t)=(\gamma_t\otimes \id)\circ \sbeta_g\circ (\gamma_t^{-1}\otimes \id)
\circ \tPhi(z,t)\circ {\salpha_g}^{-1},$$  
for any $(z,t)\in EG\times I$ and $g\in G$. 

Note that the two end points of the path $\{\tPhi(0_{n-1},y,1)\}_{y\in I}$ are in $\mathrm{Inn}(\sA)$, and 
it determines a homotopy type of a loop in $\Aut(\sA)_0$. 
 
\begin{lemma} The homotopy type of the loop determined by $\{\tPhi(0_{n-1},y,1)\}_{y\in I}$ is trivial.  
\end{lemma}

\begin{proof} Since $\pi_1(\mathrm{Inn}(\sA))$ is trivial and $\tPhi(0,0,t)=\Ad Z_t$, 
the homotopy class we would like to determine is 
given by the homotopy class of the concatenation of 
$\{\tPhi(0,y,0)\}_{y\in I}=\{\Phi(0,y)\}_{y\in I}$ and $\{\tPhi(0,1,t)\}_{t\in I}$. 
For the latter, we have    
\begin{align*}
\lefteqn{\tPhi(0,1,t)
=(\gamma_t\otimes \id)\circ \sbeta_\xi\circ (\gamma_t^{-1}\otimes \id)\circ \tPhi(0,0,t)\circ {\salpha_\xi}^{-1}} \\
 &=\Ad((\gamma_t\otimes \id)\circ \sbeta_\xi\circ (\gamma_t^{-1}\otimes \id)(Z_t))\circ 
 (\gamma_t\otimes \id)\circ \sbeta_\xi\circ (\gamma_t^{-1}\otimes \id)\circ {\salpha_\xi}^{-1}.
\end{align*}
Since $\{\Ad((\gamma_t\otimes \id)\circ \sbeta_\xi\circ (\gamma_t^{-1}\otimes \id)(Z_t))\}_{t\in I}$ 
is a path in $\mathrm{Inn}(\sA)$, we can ignore it, and the remaining part is 
\begin{align*}
\lefteqn{(\gamma_t\otimes \id)\circ \sbeta_\xi\circ (\gamma_t^{-1}\otimes \id)\circ {\salpha_\xi}^{-1}} \\
 &=(\gamma_t\otimes \id)\circ \sbeta_\xi\circ  {\salpha_\xi}^{-1}\circ \Ad (\alpha_\xi(w(t)^*)\otimes 1) \\
 &=(\gamma_t\otimes \id)\circ \sbeta_\xi\circ  {\salpha_\xi}^{-1}\circ 
 \Ad (r(t)^*w(t)^*b_\xi(t)\otimes 1)) \\
 &=(\gamma_t\otimes \id)\circ \sbeta_\xi\circ  {\salpha_\xi}^{-1}\circ \Ad (r(t))^*\otimes 1)\circ 
 (\gamma_t^{-1}\otimes \id)\circ \Ad(b_\xi(t)\otimes 1). 
\end{align*}
To compute the homotopy class of this path, we can ignore the inner part $\Ad(b_\xi(t)\otimes 1)$ again. 
Moreover, note that we can replace $\sbeta_\xi\circ {\salpha_\xi}^{-1}$ with $\Phi(0,1)$, and 
$\Ad (r(t))^*\otimes 1)$ with $\Phi(0,t)^{-1}$ for the same reason. 
Thus we end up with the path 
$$(\gamma_t\otimes \id)\circ \Phi(0,1)\circ \Phi(0,t)^{-1}\circ (\gamma_t^{-1}\otimes \id).$$
This path is homotopic to $\{\Phi(0,1)\circ \Phi(0,t)^{-1}\}_{t\in I}$ via the homotopy 
$$\{\{(\gamma_{(1-s)t}\otimes \id)\circ \Phi(0,1)\circ \Phi(0,t)^{-1}\circ (\gamma_{(1-s)t}^{-1}\otimes \id)\}_{t\in I}\}_{s\in I}.$$
Thus the concatenation of 
$\{\Phi(0,y)\}_{y\in I}$ and $\{\tPhi(0,1,t)\}_{t\in I}$ is 
homotopic to the trivial path. 
\end{proof}

We are left with the case where $\beta_g=\Ad b_g\circ \alpha_g$ with an $\alpha$-cocycle $\{b_g\}_{g\in G}$,  
$\Phi(x,0)=\Ad Z_1$, and the homotopy type of the loop in $\Aut(\sA)$ determined by $\{\Phi(0_{n-1},y)\}_{y\in I}$ 
is trivial. 
Since $b_h\otimes 1=Z_1\salpha_h(Z_1^*)=W\salpha_h(W^*)$ for any $h\in N$, we have 
$$Z_1^*W\in U(M(A\otimes \K(\ell^2(G)))^{(\alpha\otimes \Ad \rho)|_N}\otimes 1).$$ 
Since $U(M(M(A\otimes \K(\ell^2(G)))^{(\alpha\otimes \Ad \rho)|_N})\otimes \K)$ is contractible, 
there is a continuous path from $1$ to $Z_1^*W$ in $U(M(\sA)^{\salpha|_N})$. 
Thus applying Theorem \ref{lifting extension} and Lemma \ref{NDR}, we may further assume 
$\Phi(x,0)=\Phi(0_{n-1},y)=\Ad W$ for any $x\in EN$, $y\in I$. 
Since $b_g\otimes 1=W\salpha_g(W^{-1})$, if we set $\Psi(z)=\Ad W^{-1}\circ \Phi(z)$, 
we have $\Psi(g\cdot z)=\salpha_g\circ \Psi(z)\circ {\salpha_g}^{-1}$.

\begin{lemma} Theorem \ref{deformation} is true if it is the case for $\beta=\alpha$ and 
$\Psi$ in place of $\Phi$ satisfying $\Psi(x,0)=\Psi(0_{n-1},y)=\id$ for any $x\in EN$ and $y\in I$. 
\end{lemma}

\begin{proof} Let $\alpha$, $\beta$, $\{b_g\}_{g\in G}$, $W$, $\Phi$, and $\Psi$ be as above. 
Assume that there exist a continuous path $\{\gamma_t\}_{t\in I}$ in $\Aut(A)_0$ with $\gamma_0=\id$ 
and an $\alpha$-cocycle $\{c_g\}_{g\in G}$ satisfying $\gamma_1\circ\alpha_g\circ\gamma_1^{-1}=\Ad c_g\circ \alpha_g$. 
Assume that there exists a continuous map $\tPsi:EG\times I\to \Aut(\sA)$ satisfying 
\begin{itemize}
\item $\tPsi(g\cdot z,t)=(\gamma_t\otimes \id)\circ \salpha_g\circ(\gamma_t\otimes \id)^{-1}\circ \tPsi(z,t)
\circ {\salpha_g}^{-1}$ for any $g\in G$, $z\in EG$, and $t\in I$.
\item $\tPsi(z,0)=\Psi(z)$ and $\tPsi(z,1)=\Ad W'$ for any $z\in EG$, where 
$$W'=\sum_{g\in G}c_{g^{-1}}\otimes E_{g,g}\otimes 1.$$
\end{itemize}
Then $\tPhi(z,t)=\Ad(\gamma_t\otimes\id)(W)\circ \tPsi(z,t)$ has the desired properties 
with the same $\gamma_t$, $\gamma_1(b_g)c_g$ in place of $c_g$, and $\gamma_1(W)W'$ in place of $W$ 
in Theorem \ref{deformation}. 
\end{proof}

\subsection{Step III}
Assume that $\alpha=\beta$ and $\Phi(x,0)=\Phi(0_{n-1},y)=\id$ for any $x\in EN$ and $y\in I$.  

Let $\Theta(x)(y)=\Phi(x,y)$ for $x\in EN$ and $y\in I$. 
Then $\Theta:EN\to \Omega\Aut(\sA)_0$ is a homotopy fixed point for the $N$-action on $\Omega\Aut(\sA)_0$ 
given by $(h\cdot \theta)(y)=\salpha_h\circ \theta(y)\circ{\salpha_h}^{-1}$ with $\Theta(0_{n-1})=\id$. 
We apply Lemma \ref{hfp} to $\Theta$ with $u_h=1$ and $U_h(y)=1$, and obtain 
continuous maps $f:EN\to U(A_\flat)_0$, 
$F:EN\times [0,1)\to U(M(\sA))$, and $\tilde{\Theta}:EN\times I\to \Omega\Aut(\sA)_0$ 
satisfying 
\begin{itemize} 
\item $f(h\cdot x)=\alpha_h(f(x))$ for any $x\in EN$ and $h\in N$, and $f(0_{n-1})=1$. 

\item $F(h\cdot x,y)=\salpha_h(F(x,y))$ for any $h\in N$, $x\in EN$, $y\in [0,1)$. 

\item $F(x,0)=F(0_{n-1},y)=1$ for any $x\in EN$, $y\in [0,1)$, 
and $F(x,y)=\hat{f}(x,y)\otimes 1$ for any $x\in EN$, $y\in [1/2,1)$, 
where $\hat{f}(x,\cdot)$ represents $f(x)$.   

\item $\{\tilde{\Theta}(\cdot,t)\}_{t\in I}$ is a continuous family of homotopy fixed points 
of the $N$-action on $\Omega\Aut(\sA)_0$. 

\item $\tilde{\Theta}(x,0)=\Theta(x)$ for any $x\in EN$, and $\tilde{\Theta}(x,1)(y)=\Ad F(x,y)$ 
for any $x\in EN$, $y\in [0,1)$.
\end{itemize}

Let $\tilde{\Phi}(x,y,t)=\tilde{\Theta}(x,t)(y)$ for any $x\in EN$, $y\in I$, and $t\in I$. 
Then $\tilde{\Phi}$ extends to a homotopy of the homotopy fixed point for the $G$-action on $\Aut(\sA)_0$ 
with $\tilde{\Phi}(x,y,0)=\Phi(x,y)$ and $\tPhi(x,y,1)=\Ad F(x,y)$ for $y\in [0,1)$. 
Thus we may assume $\Phi(x,y)=\Ad F(x,y)$ for any $x\in EN$ and $y\in [0,1)$ from the beginning.  

\begin{lemma} There exist a continuous path of unitaries $\{v_t\}_{t\in I}$ in $U(A_\flat)_0$ with $v_0=1$ and 
$\alpha_h(v_1)=v_1$ for any $h\in N$, 
a continuous map $v:[0,1)\times I\to U(A)_0$ representing $v_t$, 
and continuous maps $\tilde{f}:EN\times I\to U(A_\flat)_0$, 
$\tilde{F}:EN\times [0,1)\times I\to U(M(\sA))$ satisfying 
\begin{itemize}
\item $v(y,0)=v(0,t)=1$ for any $y\in [0,1)$ and $t\in I$.
\item $\tilde{f}(h\cdot x,t)=v_t^{-1}\alpha_h(v_t\tilde{f}(x,t))$ for any $h\in N$, $x\in EN$, $t\in I$. 
\item $\tilde{f}(x,0)=f(x)$ and $\tilde{f}(x,1)=1$ for any $x\in EN$. 
\item $\tilde{F}(h\cdot x,y,t)=(v(y,t)^{-1}\otimes 1)\salpha_h((v(y,t)\otimes 1) \tilde{F}(x,y,t))$ for any 
$h\in N$, $x\in EN$, $y\in [0,1)$, $t\in I$. 
\item $\tilde{F}(x,y,t)=\hat{\tf}(x,y,t)\otimes 1$ for any  $x\in EN$, $t\in I$, and $y\in [1/2,1)$ where 
$\hat{\tf}(x,y,t)$ represents $\tilde{f}(x,t)$.  
\item $\tilde{F}(0_{n-1},y,t)=\tilde{F}(x,0,t)=1$ for any $x\in EN$, $y\in [0,1)$, $t\in I$. 
\item $\tilde{F}(x,y,0)=F(x,y)$ for any $x\in EN$ and $y\in [0,1)$. 
\end{itemize}
\end{lemma}

\begin{proof} Applying Theorem \ref{vanishing1} to $\alpha$ and $f$ with $u_h=1$, we get 
$\tilde{f}$ and $v_t$. 
Since $v_0=1$, we can get $v:[0,1)\times I\to U(A)_0$ as above. 
Using the cell decomposition $K=\cup_{i=1}^m e_i$ as in the proof of Lemma \ref{hfp}, 
we can inductively construct the map $\tilde{F}$ and $\hat{\tilde{f}}$ on $K\times [0,1)\times I$ 
so that the desired conditions hold whenever $x$ and $h\cdot x$ belong to $K$. 
Now we can extend $\tilde{F}$ and $\hat{\tilde{f}}$ to the whole $EN\times [0,1)\times I$ by the relations 
$$\tilde{F}(h\cdot x,s,t)=(v(y,t)^{-1}\otimes 1)\salpha_h((v(y,t)\otimes 1) \tilde{F}(x,y,t)),$$
$$\hat{\tilde{f}}(h\cdot x,s,t)=v(y,t)^{-1}\alpha_h(v(y,t)\hat{\tilde{f}}(x,y,t)).$$
\end{proof}

\begin{lemma} Theorem \ref{deformation} is true if it is the case for $\beta=\alpha$ with $\Phi$ satisfying 
the following properties. 
There exist a continuous map $r:[0,1)\to U(A)_0$ with $r(0)=1$, satisfying
\begin{itemize}
\item $r(y)$ represents an element in $U(A_\flat)_0$ and  
$$\lim_{y\to 1-0} \|r(y)-\alpha_h(r(y))\|=0,$$
for any $h\in N$.
\item $\Phi(x,y)=\Ad (r(y)\otimes 1)\circ \Ad S(y)$ for any $x\in EN$ and $y\in [0,1)$, 
where $S:I\to U(M(\sA))$ is a strictly continuous map with $S(0)=S(1)=1$ given by 
$$S(y)=\sum_{h\in N,\;n\in \Z}r(y)^{-1}\alpha_{h^{-1}}(r(y))\otimes E_{\xi^nh,\xi^nh}\otimes 1, $$
for $y\in[0,1)$. 
\end{itemize} 
\end{lemma}

\begin{proof}
Recall that we are in the situation with  $\Phi(x,y)=\Ad F(x,y)$. 
The above lemma shows that we can replace $\Phi$ with 
$$\Phi(x,y)=\Ad (v(y,1)\otimes 1)\circ \Ad \tilde{F}(x,y,1).$$
We denote $r(y)=v(y,1)$. 
Then $r(y)$ and $\tilde{F}(x,y,1)$ satisfy 
$$\lim_{y\to 1-0}\|r(y)^{-1}\alpha_h(r(y))-1\|=0,$$
for any $h\in N$, and 
$$\lim_{y\to 1-0}\|\tilde{F}(x,y,1)-1\|=0,$$
uniformly on any compact subset of $EN$ as $\hat{\tilde{f}}(x,\cdot,1)$ represents 1 in $U(A_\flat)$. 
Thus we regard  $y\mapsto r(y)^{-1}\alpha_h(r(y))$ as a continuous map defined on $I$, 
and $(x,y)\mapsto \tilde{F}(x,y,1)$ as a continuous map defined on $EN\times I$.

Note that $S$ defined from $r(y)$ as in the statement is a strictly continuous loop in $U(M(\sA))$ satisfying 
$S(0)=1$, and we have $\salpha_h((r(y)\otimes 1)S(y))=(r(y)\otimes 1)S(y)$ for any $h\in N$ and $y\in [0,1)$. 
Then
$\salpha_h(S(y)^{-1})=S(y)^{-1}(r(y)^{-1}\alpha_h(r(y))\otimes 1)$ for any $y\in [0,1)$ and $h\in N$. 
Let $\tilde{F}'(x,y)=S(y)^{-1}\tilde{F}(x,y,1)$ for $x\in EN$ and $y\in I$. 
Then 
$$\Phi(x,y)=\Ad (r(y)\otimes 1)\circ \Ad S(y)\circ \Ad \tilde{F}'(x,y).$$ 
Note that we have $\tilde{F}'(h\cdot x,y)=\salpha_h(\tilde{F}'(x,y))$. 

We regard $\tilde{F'}$ as a continuous map from $EN$ to $\Omega U(M(\sA))$, where $U(M(\sA))$ is equipped with 
the strict topology, which is a homotopy fixed point of the $N$-action on $\Omega U(M(\sA))$ given by $\salpha$. 
Since $U(M(\sA))$ is contractible, so is $\Omega U(M(\sA))$. 
Thus a cell-by-cell argument as in the proof of Lemma \ref{hfp} shows that 
$\tilde{F}'$ is homotopic to the constant map 1 within the homotopy fixed points. 
This shows that we may further assume 
$$\Phi(x,y)=\Ad ((r(y)\otimes 1))\circ \Ad S(y),$$
for any $x\in EN$ and $y\in [0,1)$. 
\end{proof}

\begin{proof}[Proof of Theorem \ref{deformation}] 
Assume that we are given $\alpha$, $r$, and $\Phi$ as in the above lemma. 
Throughout the proof, the automorphism $\Ad r(t)\in \Aut(A)$ at $t=1$ is understood as $\id$. 

We set $x_h(t)=1$ for $h\in N$ and $x_\xi(t)=r(t)$ for $t\in [0,1)$. 
We can define $x_g(t)$ for the other $g\in G$ in an appropriate way so that 
$\{x_g(t)\}_{g\in G,\;t\in [0,1)}$ forms an $\alpha$-cocycle in $U(A^\flat)$ satisfying 
$$\lim_{t\to 1-0}\Ad x_g(t)\circ \alpha_g=\alpha_g.$$ 
Then Theorem \ref{KKtrivcc} shows that if we replace $(A,\alpha)$ with $(\cO_\infty\otimes A,\mu\otimes \alpha)$, 
there exist a continuous map $w :[0,1)\to U(A)$, $\gamma\in \Aut(A)$, and 
a $G$-cocycle $\{c_g\}_{g\in G}$ in $U(A)$ satisfying 
$$\gamma\circ \alpha_g\circ \gamma^{-1}=\Ad c_g\circ \alpha_g,\quad \forall g\in G,$$
$$\lim_{t\to1-0}\Ad w(t)=\gamma,$$
$$\lim_{t\to1-0}w(t)x_g(t)\alpha_g(w(t)^*)=c_g,\quad \forall g\in G.$$
We may and do assume $w(0)=1$. 
As before we set
$$Z=\sum_{h\in N,\;n\in \Z}c_{h^{-1}}\otimes E_{\xi^nh,\xi^nh}\otimes 1,$$
$$W=\sum_{g\in G}c_{g^{-1}}\otimes E_{g,g}\otimes 1.$$

Let 
$$\nu_t=\left\{
\begin{array}{ll}
\Ad w(t) , &\quad t\in [0,1)  \\
\gamma , &\quad t=1
\end{array}
\right.,
$$
and let $\gamma_t=\nu_{1-t}^{-1}\circ \gamma$. 
Then $\{\gamma_t\}_{t\in I}$ is a continuous path in $\Aut(A)$ from $\id$ to $\gamma$. 
Let 
$$b_g(t)=\left\{
\begin{array}{ll}
w(t)x_g(t)\alpha_g(w(t)^{-1}) , &\quad t\in [0,1) \\
c_g , &\quad t=1
\end{array}
\right.,
$$
for $g\in G$, and let 
$$c_h(t)=\left\{
\begin{array}{ll}
1, &\quad t=0 \\
w(1-t)^{-1}c_h\alpha_h(w(1-t)), &\quad t\in(0,1]
\end{array}
\right.,
$$
for $h\in N$. 
Then $\{b_h(t)\}_{h\in N,\; t\in I}$ and  $\{c_h(t)\}_{h\in N,\;t\in I}$ are 
continuous families of $\alpha|_N$-cocycles satisfying 
$\nu_t\circ \alpha_h\circ \nu_t^{-1}=\Ad b_h(t)\circ \alpha_h$ and  
$\gamma_t\circ \alpha_h\circ \gamma_t^{-1}=\Ad c_h(t)\circ \alpha_h$ respectively.

For $h\in N$, and $(y,t)\in I^2$, we set 
$$u_h(y,t)=\left\{
\begin{array}{ll}
r((1-t)y)^{-1}c_h(t)\alpha_h(r((1-t)y)) , &\quad (y,t)\neq (1,0)   \\
1 , &\quad (y,t)=(1,0)
\end{array}
\right..
$$
Then $\{u_h(y,t)\}_{h\in N,\;(y,t)\in I^2}$ is a continuous family of $\alpha|N$-cocycles. 
We set 
$$\tilde{S}(y,t)=\sum_{h\in N,\;n\in \Z}u_{h^{-1}}(y,t)\otimes E_{\xi^nh,\xi^nh}\otimes 1,$$
which gives a strictly continuous map $\tilde{S}:I^2\to U(M(\sA))$ satisfying $\tilde{S}(y,0)=S(y)$, 
$\tilde{S}(y,1)=Z$, and 
$$(c_h(t)\otimes 1)\salpha_h((r((1-t)y)\otimes 1)\tilde{S}(y,t))=(r((1-t)y)\otimes 1)\tilde{S}(y,t),$$
or equivalently, 
\begin{equation}\label{hfu}
(u_h(y,t)\otimes 1)\salpha_h(\tilde{S}(y,t))=\tilde{S}(y,t),
\end{equation}
for any $h\in N$, $(y,t)\in I^2$. 
Thus we have 
$$(\gamma_t\otimes \id)\circ \salpha_h\circ (\gamma_t^{-1}\otimes \id)\circ \Ad((r((1-t)y)\otimes 1)\tilde{S}(y,t))\circ {\salpha_h}^{-1}
=\Ad((r((1-t)y)\otimes 1)\tilde{S}(y,t)),$$
for any $h\in N$, $(y,t)\in I^2$.

We choose a strictly continuous map  
$$T:I^2\to U(M(M(A\otimes \K(\ell^2(G)))^{(\alpha\otimes \Ad \rho)|_N}\otimes \K))\subset U(M(\sA)^{\salpha|N}),$$ 
satisfying $T(y,0)=1$, $T(y,1)=Z^{-1}W$, 
and define a continuous map $\tPhi:EN\times I^2\to \Aut(\sA)$ by 
$$\tPhi(x,y,t)=\Ad (r((1-t)y)\otimes 1)\circ \Ad \tilde{S}(y,t)\circ \Ad T(y,t).$$
We deduce the precise condition for $T$ so that $\tPhi$ fulfills the required properties 
in Theorem \ref{deformation}.  
We have $\tPhi(x,y,0)=\Phi(x,y)$, and 
$$\tPhi(x,y,1)=\Ad (\tilde{S}(y,1)T(y,1))=\Ad(ZZ^{-1}W)=\Ad W.$$
For $h\in N$, 
\begin{align*}
\lefteqn{(\gamma_t\otimes \id)\circ \salpha_h\circ (\gamma_t^{-1}\otimes \id)\circ 
\tPhi(x,y,t)\circ{\salpha_h}^{-1}} \\
 &=(\gamma_t\otimes \id)\circ \salpha_h\circ (\gamma_t^{-1}\otimes \id)\circ 
\Ad (r((1-t)y)\otimes 1)\tilde{S}(y,t))\circ \Ad T(y,t)\circ{\salpha_h}^{-1} \\
 &=(\gamma_t\otimes \id)\circ \salpha_h\circ (\gamma_t^{-1}\otimes \id)\circ 
\Ad (r((1-t)y)\otimes 1)\tilde{S}(y,t))\circ{\salpha_h}^{-1}\circ \Ad T(y,t) \\
 &=\Ad (r((1-t)y)\otimes )\tilde{S}(y,t))\circ \Ad T(y,t)\\
 &=\tPhi(h\cdot x,y,t). 
\end{align*}
It remains to verify the condition 
$$(\gamma_t\otimes \id)\circ \salpha_\xi\circ (\gamma_t^{-1}\otimes \id)\circ 
\tPhi(x,0,t)\circ{\salpha_\xi}^{-1}=\tPhi(\xi\cdot x,1,t).$$
Since 
$$\gamma_t\circ \alpha_\xi\circ \gamma_t^{-1}=\nu_{1-t}^{-1}\circ \Ad c_\xi\circ \alpha_\xi\circ \nu_{1-t},$$
the left-hand side is 
$$(\nu_{1-t}^{-1}\circ \Ad c_\xi\circ \alpha_\xi\circ \nu_{1-t}\circ \alpha_\xi^{-1}\otimes \id)\circ
\Ad \salpha_\xi(\tilde{S}(0,t)T(0,t)).$$
On the other hand, since $x_\xi(t)=r(t)$, 
we have $r(t)=w(t)^{-1}b_\xi(t)\alpha_\xi(w(t))$, and  
$$\Ad r(y)=\nu_{y}^{-1}\circ \Ad b_\xi(y)\circ \alpha_\xi \circ \nu_y\circ \alpha_\xi^{-1}.$$
Thus the right-hand side is 
$$((\nu_{1-t}^{-1}\circ \Ad b_\xi((1-t))\circ \alpha_\xi \circ \nu_{1-t}\circ \alpha_\xi^{-1})\otimes \id)
\circ \Ad (\tilde{S}(1,t)T(1,t)).$$
This means that we are done if we can find $T$ satisfying the following additional condition: 
\begin{equation}\label{T3}
T(1,t)=\tilde{S}(1,t)^{-1}(\alpha_\xi\circ\nu_{1-t}^{-1}\circ \alpha_\xi^{-1}(b_\xi(1-t)^{-1}c_\xi)\otimes 1))\salpha_\xi (\tilde{S}(0,t)T(0,t)).
\end{equation}
We first check that Eq.(\ref{T3}) is consistent with the condition $T(y,0)=1$ and $T(y,1)=Z^{-1}W$. 
The both sides of Eq.(\ref{T3}) are 1 at $t=0$. 
At $t=1$, the left-hand side is $Z^{-1}W$, and the right-hand side is  
$$\tilde{S}(1,1)^{-1}(c_\xi\otimes 1)\salpha_\xi (\tilde{S}(0,1)T(0,1))
=Z^{-1}W\salpha_\xi (W^{-1}ZZ^{-1}W)=Z^{-1}W,$$
which shows that Eq.(\ref{T3}) is consistent with the other two conditions. 
Since the group $U(M(M(A\otimes \K(\ell^2(G)))^{(\alpha\otimes \Ad \rho)|_N}\otimes \K))$ is contractible, we can obtain $T$ satisfying Eq.(\ref{T3}) if 
$$\tilde{S}(1,t)^{-1}(
\alpha_\xi\circ\nu_{1-t}^{-1}\circ \alpha_\xi^{-1}(b_\xi(1-t)^{-1}c_\xi)
\otimes 1))\salpha_\xi (\tilde{S}(0,t)) $$ 
belongs to $M(A\otimes \K(\ell^2(G)))^{(\alpha\otimes \Ad \rho)|_N}\otimes \C$. 
We assume $0<t<1$ now. 
In this case we have 
\begin{align*}
\lefteqn{\alpha_\xi\circ\nu_{1-t}^{-1}\circ \alpha_\xi^{-1}(b_\xi(1-t)^{-1}c_\xi)} \\
 &=\alpha_\xi(w(1-t)^{-1}\alpha_\xi^{-1}(\alpha_\xi(w(1-t))r(1-t)^{-1}w(1-t)^{-1}c_\xi)w(1-t)) \\
 &=r(1-t)^{-1}w(1-t)^{-1}c_\xi\alpha_\xi(w(1-t)), 
\end{align*}
which we denote by $\varphi(t)$. 
Thanks to Eq.(\ref{hfu}), all we have to verify is 
$$u_h(1,t)\alpha_h(\varphi(t))\alpha_\xi(u_{h^\xi}(0,t)^{-1})=\varphi(t).$$
Indeed, \begin{align*}
\lefteqn{u_h(1,t)\alpha_h(\varphi(t))\alpha_\xi(u_{h^\xi}(0,t)^{-1})} \\
 &=u_h(1,t)\alpha_h(r(1-t)^{-1}w(1-t)^{-1}c_\xi)\alpha_\xi(\alpha_{h^\xi}(w(1-t))u_{h^\xi}(0,t)^{-1})\\
 &=r(1-t)^{-1}c_h(t)\alpha_h(w(1-t)^{-1}c_\xi)\alpha_\xi(\alpha_{h^\xi}(w(1-t))c_{h^\xi}(t)^{-1})\\
 &=r(1-t)^{-1}w(1-t)^{-1}c_h\alpha_h(c_\xi)\alpha_\xi(c_{h^\xi}^{-1}w(1-t))\\
 &=r(1-t)^{-1}w(1-t)^{-1}c_{\xi h^\xi}\alpha_\xi(c_{h^\xi}^{-1})\alpha_\xi(w(1-t))  \\
 &=r(1-t)^{-1}w(1-t)^{-1}c_\xi\alpha_\xi(w(1-t))\\
 &=\varphi(t),
\end{align*}
which finishes the proof. 
\end{proof}

\section{Dynamical realization theorem}\label{sDRT} 
Although the surjectivity of $\cB$ in Conjecture \ref{main conjecture s} holds in full generality thanks to 
Ralf Meyer's result \cite[Theorem 3.10]{Me19}, we give a proof in the case of poly-$\Z$ groups without using equivariant KK-theory 
for the sake of self-contained presentation.  

Recall that $BG=EG/G$ is identified with the mapping torus 
$$\frac{BN\times I}{(x,0)\sim (T_\xi(x),1)},$$
where $T_\xi$ is a homeomorphism of $BN$ induced by $\xi$. 
We denote by $q$ the quotient map $q:BN\times I\to BG$, and identify 
$BN$ with the submanifold $q(BN\times \{0\})$ of $BG$. 

Let $\Gamma$ be a topological group, and let $p:\cQ\to BN$ be a principal $\Gamma$-bundle over $BN$. 
Let $T_\xi^*\cQ$ be the pullback bundle, that is, 
$$T_\xi^*\cQ=\{(x,\eta)\in BN\times \cQ;\; T_\xi(x)=p(\eta) \}.$$
Let $\hat{T}_\xi:T_\xi^*\cQ\to \cQ$ be the projection onto the second component, which is a bundle map covering $T_\xi$. 
Assume that $f:\cQ\to {T_\xi}^*\cQ$ is an isomorphism as principal $\Gamma$-bundles over $BN$. 
Then we can extend $\cQ$ to a principal $\Gamma$-bundle $\cQ_{f}$ over $BG$ by setting 
$$\cQ_f=\frac{\cQ\times I}{(\eta,0)\sim (\hat{T}_\xi \circ f(\eta),1)}.$$  

\begin{lemma}\label{mappingtorus} Let $\Gamma$ be a topological group, and let $\cQ$ be a principal $\Gamma$-bundle over $BN$. 
\begin{itemize}
\item[$(1)$] The principal bundle $\cQ$ extends to a principal $\Gamma$-bundle over $BG$ if and only if $\cQ$ is isomorphic to 
$T_\xi^*\cQ$. 
If such an extension exists, every extension is isomorphic to  $\cQ_f$ for some isomorphism $f:\cQ\to T_\xi^*\cQ$. 
\item[$(2)$] Let $f_0:\cQ\to {T_\xi}^*\cQ$ be an isomorphism, and set $l_\xi=\hat{T}_\xi\circ f_0$. 
Then the isomorphism classes of extensions of $\cQ$ to principal $\Gamma$-bundles over $BG$ are parameterized by 
$$H^1(\Z,\pi_0(\Aut(\cQ)))=\pi_0(\Aut \cQ)/\sim,$$
where $[\varphi]\sim [\varphi']$ if and only if there exists $[\psi]$ satisfying 
$[\varphi']= [\psi^\xi][\varphi][\psi^{-1}]$ with  
$\psi^\xi=l_\xi^{-1}\circ \psi\circ l_\xi$.   
\end{itemize} 
\end{lemma}

\begin{proof} (1) Let $p:\cP\to BG$ be a principal $\Gamma$-bundle extending $\cQ$. 
Then $q^*\cP$ is a principal $\Gamma$-bundle over $BN\times I$ whose restriction to $BN\times \{0\}$ is $\cQ$. 
By definition, we have 
$$q^*\cP=\{(x,t,\zeta)\in BN\times I\times \cP;\; q(x,t)=p(\zeta)\},$$
and $\cP$ is identified with the quotient of $q^*\cP$ by the relation $(x,0,\zeta)\sim (T_\xi(x),1,\zeta)$. 
From \cite[Theorem 4.9.9]{Hu}, there exists an isomorphism $\Sigma:\cQ\times I\to q^*\cP$ of the form 
$\Sigma(\eta,t)=(p(\eta),t,\sigma_t(\eta))$ with $\sigma_t:\cQ\to \cP$ such that $\sigma_0:\cQ\to \cP$ is the inclusion map, and 
$p(\sigma_t(\eta))=q(p(\eta),t)$ and $\sigma_t(\eta\cdot\gamma)=\sigma_t(\eta)\cdot \gamma$ hold for any $t\in I$, $\eta\in \cQ$, $\gamma\in \Gamma$. 
At $t=1$, we get $p(\sigma_1(\eta))=q(p(\eta),1)=q(T_\xi^{-1}(p(\eta)),0)$, and we may regard $\sigma_1$ as a homeomorphism from $\cQ$ onto itself 
satisfying $p(\sigma_1(\eta))=T_\xi^{-1}(p(\eta))$. 
Thus $\sigma_1\circ \hat{T}_\xi$ is an isomorphism from $T_\xi^*\cQ$ onto $\cQ$, and we denote by $f$ its inverse. 
Since $\Sigma(\sigma_1(\eta),0)=(T_\xi^{-1}(p(\eta)),0,\sigma_1(\eta))$ and 
$\Sigma(\eta,1)=(p(\eta),1,\sigma_1(\eta))$, we get
$$\cP\cong \frac{\cQ\times I}{(\sigma_1(\eta),0)\sim (\eta,1)}=\frac{\cQ\times I}{(\eta,0)\sim (\hat{T}_\xi\circ f(\eta),1)}.$$

(2) Thanks to (1), every extension is isomorphic to $\cQ_{f_0\circ \varphi}$ with $\varphi\in \Aut(\cQ)$. 
Since two homotopic principal $\Gamma$-bundles are isomorphic, 
the isomorphism class of $\cQ_{f_0\circ \varphi}$ depends only on the homotopy class of $\varphi$.     
Moreover, every isomorphism from $\cQ_{f_0\circ \varphi}$ to $\cQ_{f_0\circ \varphi'}$ is of the form 
$$[(\eta,t)]_{Q_{f_0\circ \varphi}}\mapsto [(\psi_t(\eta),t)]_{\cQ_{f_0\circ \varphi'}},$$
where $\{\psi_t\}_{t\in I}$ is a continuous family in $\Aut(\cQ)$ satisfying 
$$\psi_1\circ l_\xi\circ \varphi=l_\xi\circ \varphi'\circ \psi_0.$$
Therefore we get the statement by taking its homotopy class. 
\end{proof}

\begin{example}\label{Qf0}
Let $A$ be a $C^*$-algebra, let $\alpha$ be a  $G$-action on $A$, and 
let $\cQ$ be the restriction of $\cP_{\salpha}$ to $BN$, that is $\cQ=\cP_{\salpha|N}$. 
Let $\alpha'$ be an $N$-action on $A$ given by $\alpha'_{h}=\alpha_{h^{\xi^{-1}}}$. 
Then $T_\xi^*\cQ$ and $\hat{T}_\xi$ are identified with $\cP_{{\alpha'}^s}$ and $\hat{T}_\xi([(x,\gamma)])=[(\xi\cdot x,\gamma)]$ 
respectively. 
We can define an isomorphism $f_0:\cQ\to T_\xi^*\cQ$ by $f_0([(x,\gamma)])=[(x,\salpha_\xi\circ \gamma)]$. 
Then $\hat{T}_\xi\circ f_0([(x,\gamma)])=[(\xi\cdot x,\salpha_\xi\circ \gamma)]$, and 
we get $\cP_{\salpha}=\cQ_{f_0}$.  

Let $\beta$ be another outer action of $G$ on $A$ with $\beta_h=\alpha_h$ 
for any $h\in N$. 
Then $\cP_{\sbeta}$ is isomorphic to $\cQ_{f_0\circ \varphi}$ 
with $\varphi([(x,\gamma)])=[(x,{\salpha_\xi}^{-1}\circ \sbeta_\xi\circ \gamma)]$. 
\end{example}

The following statement holds in a much more general situation (see \cite[Proposition 4.17]{Sz18CMP}). 
We include an elementary proof in the specific case we need. 

\begin{lemma}\label{homotopy of mu} Let $G$ be a poly-$\Z$ group. 
Then there exists an outer $G$-action $\mu$ on $\cO_\infty$ that is homotopic to the trivial action 
within $G$-actions. 
\end{lemma}

\begin{proof} We first claim that there exists a faithful unitary representation 
$(\pi,\cH)$ of $G$ that is homotopic to the trivial representation $(1,\cH)$ within unitary representations 
in $\cH$. 
We show the claim by induction of the Hirsch length $h(G)$. 
If $h(G)=0$, there is nothing to prove. 
Assume that the claim is the case for $N$, and $(\pi_t,\cH)$ is a continuous family of unitary representations 
of $N$ such that $\pi_0$ is faithful and $\pi_1(h)=1$ for any $h\in N$. 
Let $\tilde{\cH}=\cH^{\oplus \Z}$. 
For $0\leq t\leq 1/2$, we set 
$$\tilde{\pi}_t(h)=\bigoplus_{n\in \Z}\pi_{2t}(h^{\xi^n}),$$ 
for $h\in N$, and $\tilde{\pi}(\xi)$ to be the shift of the index set of $\tilde{\cH}$. 
Then $\tilde{\pi}_t$ is a unitary representation of $G$ for each $0\leq t\leq 1/2$, and $\tilde{\pi}_0$ is faithful. 
Moreover $\tilde{\pi}_{1/2}(h)=1$ for any $h\in N$. 
Now by using functional calculus, we can construct a homotopy $\{\tilde{\pi}_t\}_{t\in [1/2,1]}$ with $\tilde{\pi}_1(g)=1$ for any $g\in G$, 
and the claim is shown. 

The quasi-free action of $G$ corresponding to the representation $(\tilde{\pi}_0,\tilde{\cH})$ has the desired property.  
\end{proof}

Now we prove the surjectivity of $\cB$ in Theorem \ref{main theorem s}. 
We fix a $G$-action $\mu$ on $\cO_\infty$ as above. 

\begin{theorem}[Dynamical Realization Theorem]\label{DRT} Let $G$ be a poly-$\Z$ group and let $A$ be a stable Kirchberg algebra. 
Then for every principal $\Aut(A)$-bundle $\cP$ over $BG$, there exists an outer $G$-action $\alpha$ on 
$A$ such that $\cP$ is isomorphic to $\cP_{\alpha}$. 
\end{theorem}

\begin{proof} 
Throughout the proof, we make identification of the $C^*$-algebras of compact operators on various separable infinite dimensional 
Hilbert spaces without mentioning it to avoid heavy notation.  
We also use Lemma \ref{tensoring} and Lemma \ref{de-stabilization} freely without mention. 
We prove the statement by induction of the Hirsch length $h(G)$ of $G$. 
The statement is trivially true for $h(G)=0$. 
Assume that the statement is true for $N$. 

Let $\cQ$ be the restriction of $\cP$ to $BN$. 
By induction hypothesis, there exists an outer $N$-action $\nu$ on $A$ such that 
$\cQ$ is isomorphic to $\cP_{\nu}$.  
Let $\nu'_h=\nu_{h^{\xi^{-1}}}$. 
Then we can identify $T_\xi^*\cQ$ with $\cP_{\nu'}$, and Lemma \ref{mappingtorus} shows that there exists an isomorphism  
$\varphi:\cP_{\nu}\to\cP_{\nu'}$ satisfying $\cQ_\varphi\cong \cP$. 
We denote by $\Phi:EN\to \Aut(A)$ the corresponding homotopy fixed point of $\Aut(A)$ with an $N$-action 
$h\cdot \gamma=\nu'_h \circ\gamma\circ {\nu_h}^{-1}$. 
Thanks to the injectivity of $\cB$ in Theorem \ref{main conjecture s}, there exist $\theta\in \Aut(A)$ and a $\nu$-cocycle 
$\{b_h\}_{h\in N}$ satisfying $\theta^{-1}\circ \nu'_h\circ \theta=\Ad b_h\circ \nu_{h}$. 
Furthermore, the proof of the injectivity of $\cB$  shows that we may assume $KK(\theta)=KK(\Phi(x))$ for any $x\in EN$.  
If we set $\tilde{\nu}_h=\nu_h$ for $h\in N$ and $\tilde{\nu}_\xi=\theta$, and extend $\tilde{\nu}$ in an appropriate way, 
we get a cocycle action $(\tilde{\nu},w)$ of $G$ on $A$ with $w_{h,k}=1$ for all $h,k\in N$. 
By replacing $\tilde{\nu}_\xi$ with $\kappa\otimes \tilde{\nu}_\xi$ if necessary, 
where $\kappa$ is an aperiodic automorphism of $\cO_\infty$, 
we may further assume that the cocycle action $(\tilde{\nu},w)$ is outer. 
Since $A$ is stable, the cocycle action $(\tilde{\nu},w)$ is equivalent to a genuine action, say $\rho$, and $\rho$ satisfies 
$\cP_{\rho}|_{BN}\cong \cQ$ by construction. 
Let $f_0$ be the isomorphism from $\cP_{\rho}|_{BN}$ to $T_\xi^*\cP_{\rho}|_{BN}$ as in Example \ref{Qf0}. 
Then $\cP$ is isomorphic to $(\cP_{\rho}|_{BN})_{f_0\circ \psi}$ with $\psi\in \Aut(\cQ)$, and by construction the homotopy fixed point 
$\Psi$ corresponding to $\psi$ satisfies $\Psi(EN)\subset \Aut(A)_0$. 

Remark \ref{perturbation} implies that we may assume that there exists an outer $G$-action 
$\alpha$ on $B$ in the Cuntz standard form with $\rho=\salpha$ by cocycle perturbation. 
For our purpose, we may replace $(B,\alpha,\Psi)$ with 
$(\cO_\infty\otimes  B,\id\otimes \alpha,\id\otimes \Psi)$, 
and we may further replace $\id\otimes \alpha$ with $\mu\otimes \alpha$ thanks to Lemma \ref{homotopy of mu} 
because homotopic principal bundles are isomorphic. 
Thus our task is to realize $(\cP_{{\mu\otimes \alpha}^{\mathrm{s}}}|_{BN})_{f'_0\circ(\id\otimes \psi)}$ by a group action, 
where $f_0'([(x,\gamma)])=[(x,(\mu_\xi\otimes \salpha_\xi)\circ \gamma)]$.

Theorem \ref{deformation} implies that there exist a continuous family $\{\gamma_t\}_{t\in I}$ of automorphisms 
of $\cO_\infty\otimes B$, a $\mu\otimes \alpha|_N$-cocycle $\{c_h\}_{h\in N}$, and a continuous map 
$$\tPsi:EN\times I\to 
\Aut(\cO_\infty\otimes B \otimes \K(\ell^2(G))\otimes \K)$$ satisfying 
\begin{itemize}
\item $\gamma_0=\id$, $\gamma_1\circ (\mu_h\otimes \alpha_h)\circ \gamma_1^{-1}=\Ad c_h\circ (\mu_h\otimes \alpha_h)$ for any $h\in N$. 
\item $\tPsi(h\cdot x,t)=(\gamma_t\otimes\id)\circ (\mu_h\otimes \salpha_h)\circ (\gamma_t\otimes \id)^{-1}\circ 
\tPsi(x,t)\circ (\mu_h\otimes \salpha_h)^{-1}$, 
for any $x\in EN$, $t\in I$, $h\in N$. 
\item $\tPsi(x,0)=\Psi(x)$, $\tPsi(x,1)=\Ad (W\otimes 1)$ for any $x\in EN$, where $W$ is a unitary in 
$M(\cO_\infty\otimes B\otimes \K(\ell^2(G)))$ satisfying $(W\otimes 1)\salpha_h(W^{-1}\otimes 1)=c_h\otimes 1$. 
\end{itemize}

Let $\Psi_t(x)=(\gamma_t\otimes \id)^{-1}\circ \tPsi(x,t)$. 
Then $\{\Psi_t\}_{t\in [0,1]}$ gives rise to a homotopy $\{\psi_t\}_{t\in [0,1]}$ in 
$\Aut(\cP_{{\mu\otimes \alpha}^{\mathrm{s}}}|_{BN})$ with $\psi_0=\psi$.
Thus $(\cP_{{\mu\otimes \alpha}^{\mathrm{s}}}|_{BN})_{f'_0\circ(\id\otimes \psi)}$ is isomorphic to 
$(\cP_{{\mu\otimes \alpha}^{\mathrm{s}}}|_{BN})_{f'_0\circ\psi_1}$, 
and the homotopy fixed point corresponding to $\psi_1$ is a constant map $(\gamma_1\otimes \id)^{-1}\circ (\Ad W\otimes \id)$. 
Note that it commutes with $\mu_h\otimes \salpha_h$ for any $h\in N$. 

Let $\beta_h=\mu_h\otimes \salpha_h$ for $h\in N$, and 
$$\beta_\xi=\salpha_\xi\circ (\gamma_1\otimes \id_{\K(\ell^2(G))\otimes \K})^{-1}\circ \Ad (W\otimes 1).$$
Then $\beta$ extends to an action of $G$ on $\cO_\infty\otimes B\otimes \K(\ell^2(G))\otimes \K$, and by construction, 
the principal bundle $\cP_{\beta}$ is isomorphic to $(\cP_{{\mu\otimes \alpha}^{\mathrm{s}}}|_{BN})_{f'_0\circ\psi_1}$. 
Taking tensor product with $(\cO_\infty,\mu)$ if necessary, we may assume that $\beta$ is outer, and we get the statement. 
\end{proof}

Let $A$ be a Kirchberg algebra and let $\alpha$ be a $G$-action with $KK(\alpha_g)=KK(\id)$ for any $g\in G$. 
We define a principal $\Aut(\sA)_0$-bundle $\cP_{\salpha}^0$ by 
$EG\times_G \Aut(\sA)_0$ with a $G$-action on $\Aut(\sA)_0$ given by $g\cdot \gamma=\salpha_g\circ \gamma$. 
We can show the following.

\begin{cor}\label{DRT0} Let $G$ be a poly-$\Z$ group and let $A$ be a stable Kirchberg algebra such that 
$\Aut(A)_0$ is open in $\Aut(A)$. 
Then for any principal $\Aut(A)_0$-bundle $p:\cP^0\to BG$, there exists an outer $G$-action on 
$A$ such that $KK(\alpha_g)=KK(\id)$ for any $g\in G$ and $\cP^0$ is isomorphic to $\cP_{\alpha}^0$. 
\end{cor}

\begin{proof} Let $\cP=\cP^0\times_{\Aut(A)_0}\Aut(A)$, which is a principal $\Aut(A)$-bundle over $BG$. 
We identify $\cP^0$ with the subset $\{[(x,\id)]\in \cP;\;x\in \cP^0\}$. 
Let $\Aut(A)=\sqcup_{j\in \pi_0(\Aut(A))}\Aut(A)_0\gamma_j$ be the coset decomposition, which is a topologically disjoint union 
as $\Aut(A)_0$ is an open subgroup of $\Aut(A)$. 
Then $\cP$ is a topologically disjoint union $\sqcup_{j\in \pi_0(\Aut(A))} \cP^0\cdot\gamma_j$. 
In particular $\cP^0$ is a (path) connected component in $\cP$, and 
$\pi_0(\Aut(A))$ acts on the set of connected components of $\cP$ freely and transitively. 

Thanks to Theorem \ref{DRT}, there exists an outer $G$-action $\alpha$ on $A$ with $\cP_{\alpha}$ isomorphic to $\cP$. 
Since we have 
$$[(x,\id)]=[(g\cdot x,\alpha_g)]=[(g\cdot x, \id)]\cdot\alpha_g$$ 
in $\cP_\alpha$ and $\pi_0(\Aut(A))$ acts on the set of connected components of $\cP_\alpha$ freely, 
we see that $\alpha_g\in \Aut(A)_0$ for all $g\in G$. 
Thus $\cP_\alpha^0$ makes sense and we have a decomposition 
$\cP_{\alpha}=\sqcup_{j\in \pi_0(\Aut(A))} \cP^0_{\alpha}\cdot\gamma_j$. 
Let $f:\cP\to \cP_{\alpha}$ be an isomorphism. 
Then there exists $j\in \pi_0(\Aut(A))$ with $f(\cP^0)=\cP_{\alpha}^0\cdot\gamma_j$,  
which means that there exists a homeomorphism $f_0:\cP^0\to \cP_\alpha^0$ satisfying $p(f_0(x))=p(x)$ and 
$f([(x,\id)])=[(f_0(x),\gamma_j)]$. 
For $x\in \cP^0$ and $\gamma\in \Aut(A)_0$, we get 
$$[(f_0(x\cdot \gamma),\gamma_j)]=f([(x,\gamma)])=f([(x,\id)])\cdot\gamma=[(f_0(x),\gamma_j\gamma)]=[(f_0(x)\cdot \gamma_j\gamma\gamma_j^{-1},\gamma_j)],$$
and $f_0(x\cdot \gamma)=f_0(x)\cdot (\gamma_j\gamma\gamma_j^{-1})$. 
Let $\beta_g=\gamma_j\circ \alpha_g\circ \gamma_j^{-1}$. 
Then $\cP^0$ is isomorphic to $\cP_{\beta}^0$. 
\end{proof}

\begin{remark} Let $A$ be a Kirchberg algebra and let $\alpha$ and $\beta$ be outer actions of a poly-$\Z$ group $G$ 
on $A$ with $KK(\alpha_g)=KK(\beta_g)=KK(\id)$ for any $g\in G$. 
Then $\alpha$ and $\beta$ are $KK$-trivially cocycle conjugate if and only if $\cP_{\salpha}^0$ and $\cP_{\sbeta}^0$ are isomorphic. 
Indeed, by Remark \ref{perturbation} the proof of the statement is  reduced to that in the unital case, which in turn follows from 
Theorem \ref{deformation}.   
  
\end{remark}

\section{Primary obstructions}
Let $G$ be a poly-$\Z$ group, let $A$ be a unital Kirchberg algebra, and let $\alpha$ and $\beta$ be $G$-actions on $A$ with 
$KK(\alpha_g)=KK(\beta_g)$. 
We introduced in \cite[Section 7]{IMI} an obstruction class $\fo^2(\alpha,\beta)\in H^2(G;KK^1(A,A))$ 
for $\alpha$ and $\beta$ to be $KK$-trivially cocycle conjugate. 
Let us recall the definition of $\fo^2(\alpha,\beta)$ first. 
Since $KK(\alpha_g)=KK(\beta_g)$, there exist unitaries $v_g\in A^\flat$ satisfying 
$$\lim_{t\to 1-0}\Ad v_g(t)\circ \alpha_g=\beta_g.$$
Let $w_{g,h}=v_g\alpha_g(v_h)v_{gh}^*\in A_\flat$, and let $\theta_g=\Ad v_g\circ \alpha_g|_{A_\flat}$. 
Then $(\theta,w)$ is a cocycle action of $G$ on $A_\flat$. 
Theorem \ref{KKtrivcc} shows that $\alpha$ and $\beta$ are $KK$-trivially cocycle conjugate if and only if $(\theta,w)$ 
is equivalent to a genuine $G$-action. 
Taking the $K_1$-classes of $w_{g,h}$, we get a 2-cocycle $([w_{g,h}]_1)_{g,h}$ of $G$ with values in $K_1(A_\flat)$, 
giving rise to a cohomology class 
$\fo^2(\alpha,\beta)\in H^2(G,K_1(A_\flat))$. 
Thanks to Theorem \ref{weak homotopy} and Dadarlat's result \cite{D07}, the coefficient group $K_1(A_\flat)$ is identified with 
$\pi_1(\Aut(\sA))\cong KK^1(A,A)$,  
The purpose of this section is to prove the following theorem. 
The reader is referred to \cite[Chapter 7]{DK} and \cite[Chapter IV]{W78} for obstruction theory. 

\begin{theorem}\label{obstruction} The obstruction class $\fo^2(\alpha,\beta)$ is identified with the primary obstruction in 
$H^2(BG;\pi_1(\Aut(\sA)))$ for the existence of a continuous section of the fiber bundle 
$\cI_{\salpha,\sbeta}^0\to BG$ introduced in Definition \ref{isobundle}. 
\end{theorem}

\begin{proof} 
We adopt Milgram's geometric bar construction \cite{Mi67}  as a model 
of $EG$, instead of $\R^n$, which is the only exception in this paper. 
Let $\Delta^n$ be the geometric $n$-simplex 
$$\Delta^n=\{(t_0,t_1,\cdots,t_n)\in \R^{n+1};\; \sum_{i=0}^nt_i=1,\;t_i\geq 0\}.$$
We define $d^i:\Delta^{n-1}\to \Delta^n$ for $0\leq i\leq n$, and $s_i:\Delta^{n+1}\to \Delta^n$ 
for $1\leq i\leq n$ by 
$$d^i(t_0,\cdots,t_{n-1})=(t_0,\cdots,t_{i-1},0,t_{i+1},\cdots,t_{n-1}),$$
$$s^i(t_0,\cdots,t_{n+1})=(t_0,\cdots,t_{i-1},t_i+t_{i+1},t_{i+2},\cdots,t_{n+1}).$$
Then 
$$EG=(\coprod_{k=0}^\infty G\times \Delta^k\times G^k)/\sim,$$
where the equivalence relation $\sim$ is generated by 
$$(g_0;d^i(t);g_1,\cdots,g_n) 
 \sim \left\{
\begin{array}{ll}
(g_0g_1;t;g_2,\cdots,g_n) , &\quad i=0 \\
(g_0;t;g_1,\cdots,g_ig_{i+1},\cdots,g_n) , &\quad 1\leq i\leq n-1 \\
(g_0;t;g_1,\cdots,g_{n-1}) , &\quad i=n
\end{array}
\right.$$
$$(g_0;t;g_1,\cdots,g_{i-1},e,g_{i+1},\cdots,g_n)\sim (g_0;s^i(t);g_1,\cdots,g_{i-1},g_{i+1},\cdots,g_n),$$
and a $G$-action is given by $g\cdot(g_0;t;g_1,\cdots,g_n)=(gg_0;t;g_1,\cdots,g_n)$. 
The $n$-skeleton on $EG$ is 
$$E_nG=(\coprod_{k=0}^n G\times \Delta^k\times G^k)/\sim,$$
and we set $B_nG=E_nG/G$. 

We introduce a $G$-action on $\Aut(\sA)$ by $g\cdot \gamma=\sbeta_g\circ \gamma \circ {\salpha_g}^{-1}$. 
We choose a norm continuous path $\{V_g(t)\}_{t\in [0,1)}$ in $U(M(\sA))$ satisfying $V_g(0)=1$ and $V_g(t)=v_g(t)\otimes 1$ for $t\in [1/2,1)$. 
For $g\in G$, we set 
$$\varphi_g(t)=\left\{
\begin{array}{ll}
\sbeta_g\circ {\salpha_g}^{-1} , &\quad t=1  \\
\Ad V_g(t) , &\quad 0\leq t<1
\end{array}
\right.
$$
Then $\{\varphi_g(t)\}_{t\in I}$ is a continuous path from $\id$ to $g\cdot \id$. 
We introduce a continuous $G$-equivariant map $\Phi:E_1G\to \Aut(\sA)_0$ by 
$\Phi([g_0;(1-t,t);g_1])=\sbeta_{g_0}\circ\varphi_{g_1}(t)\circ {\salpha_{g_0}}^{-1}$, 
which gives a partial section of $p:\cI_{\salpha,\sbeta}^0\to BG$ defined on $B_1G$. 
Note that we have 
$$\partial [e;\Delta^2;g_1,g_2]=[e;\Delta^1,g_1]\cup g_1\cdot [e;\Delta^1;g_2]\cup [e;\Delta^1;g_1g_2],$$
and the element $c_{g_1,g_2}\in \pi_1(\Aut(\sA))$ given by the restriction of $\Phi$ to $\partial  [e;\Delta^2;g_1,g_2]$ is 
the homotopy class of the concatenation of the following paths: $\{\varphi_{g_1}(t)\}_{t\in I}$, $\{g_1\cdot\varphi_{g_2}(t) \}_{t\in I}$, 
and $\{\varphi_{g_1g_2}(1-t)\}_{t\in I}$. 
Let 
$$f_1(t)=\left\{
\begin{array}{ll}
2t , &\quad 0\leq t\leq 1/2  \\
1 , &\quad 1/2\leq t\leq 1
\end{array}
\right.
$$
$$f_2(t)=\left\{
\begin{array}{ll}
0 , &\quad 0\leq t\leq 1/2  \\
2t-1 , &\quad 1/2\leq t\leq 1
\end{array}
\right.
$$
The concatenation of the first two paths is 
$$\{\varphi_{g_1}(f_1(t))\circ\salpha_{g_1}\circ\varphi_{g_2}(f_2(t))\circ{\salpha_{g_1}}^{-1}\}_{t\in I},$$ 
and it is homotopic to 
$\{\varphi_{g_1}(t)\circ\salpha_{g_1}\circ\varphi_{g_2}(t)\circ {\salpha_{g_1}}^{-1}\}_{t\in I}$, 
which we denote by $\{\psi(t)\}_{t\in I}$. 
On the other hand, the third path is homotopic to $\{\varphi_{g_1g_2}(1)\circ\varphi_{g_1g_2}(t)^{-1}\}_{t\in I}$.   
Thus the concatenation of $\{\psi(t)\}_{t\in I}$ and $\{\varphi_{g_1g_2}(1)\circ\varphi_{g_1g_2}(t)^{-1}\}_{t\in I}$ 
is  
$$\{\psi(f_1(t))\circ\varphi_{g_1g_2}(f_2(t))^{-1}\}_{t\in I},$$
which is homotopic to $\{\psi(t)\circ \varphi_{g_1g_2}(t)^{-1}\}_{t\in I}$. 
For $0\leq t<1$, we have 
$$\psi(t)\circ \varphi_{g_1g_2}(t)^{-1}=\Ad (V_{g_1}(t)\salpha_{g_1}(V_{g_2}(t))V_{g_1g_2}(t)^*),$$
and for $1/2\leq t<1$ we have  $V_{g_1}(t)\salpha_{g_1}(V_{g_2}(t))V_{g_1g_2}(t)^*=w_{g_1,g_2}(t)\otimes 1$. 
Thus from the construction of the map $\Pi_{A,0}:\pi_0(U(A_\flat))\to \pi_1(\Aut(\sA))$, we get 
$c_{g_1,g_2}=\Pi_{A,0}([w_{g_1,g_2}]_1)$, and we can identify $\fo^2(\alpha,\beta)$ with the primary obstruction for 
the existence of a section of $p:\cI_{\salpha,\sbeta}^0\to BG$. 
\end{proof}

When $\fo^2(\alpha,\beta)$ vanishes, we can define another obstruction class $\fo^3(\alpha,\beta)$, 
which is in a certain quotient of 
$$H^3(G,K_0(A_\flat))\cong H^3(G,KK(A,A))\cong H^3(BG,\pi_2(\Aut(A\otimes \K))),$$
(see \cite[Section 8]{IMI}). 
Although we believe that $\fo^3(\alpha,\beta)$ can be identified with the secondary obstruction, 
it would be too tedious an argument to present here. 
We will discuss it elsewhere. 

\section{The Cuntz algebra case}

We first recall Dadarlat's classification result \cite[Theorem 1.6]{D12} of continuous fields of Cuntz algebras, 
with a little improvement by Izumi-Sogabe \cite[Theorem 4.12]{IS}. 
We denote by $\tilde{H}^*(X)$ and $\tilde{K}^*(X)$ the reduced cohomology and reduced $K$-theory respectively. 

\begin{theorem} 
Let $m$ be a natural number, and let $\cO_{m+1}$ be the Cuntz algebra. 
Let $X$ be a finite CW-complex with $\Tor(H^*(X),\Z_m)=\{0\}$. 
Then 
$$|[X,B\Aut(\cO_{m+1})]|=|\tilde{H}^{\mathrm{even}}(X)\otimes \Z_m|=|\tilde{K}^0(X)\otimes \Z_m|.$$ 
\end{theorem}

Let $B$ be a locally trivial continuous field of $\cO_{m+1}$ over $X$ as above, and let $\iota_B:C(X)\hookrightarrow B$ 
be the inclusion map (in fact, local triviality is automatic thanks to \cite{D09}). 
We denote by $\M_{(m)}$ the UHF algebra 
$$\M_{(m)}=\bigotimes_{p\in P_m}\M_{p^\infty},$$
where $P_m$ denotes the set of all primes $p$ with $\gcd(p,m)=1$, and $\M_{p^\infty}$ denotes the UHF algebra of type $p^\infty$. 
Then we can identify $K_0(\M_{(m)})$ with $\Z_{(m)}$, the localization of $\Z$ at $m$, and 
$K_0(C(X)\otimes \M_{(m)})$ with $K^0(X)\otimes \Z_{(m)}$.  

Note that 
$$\ker((\iota_B\otimes \id_{\M_{(m)}})_*:K^0(X)\otimes \Z_{(m)}\to K_0(B\otimes \M_{(m)}))$$
is a complete invariant of $B$ as a $C(X)$-algebra. 
Indeed, it is shown in \cite[Theorem 4.12]{IS} that $B$ is always of the form $\cO_E\otimes M_{(m)}$, 
where $E$ is a vector bundle of rank $mr+1$ with $(m,r)=1$ and $\cO_E$ is the Cuntz-Pimsner algebra arising from $E$ 
as a $C(X)$-$C(X)$ bimodule. 
By \cite[Lemma 4.7]{IS}, the principal ideal $(1-[E])K^0(X)\otimes \Z_{(m)}$ in $K^0(X)\otimes \Z_{(m)}$ 
is a complete invariant of $\cO_E$ as $C(X)$-algebra. 
Let $\cT_E$ be the Toeplitz extension of $\cO_E$.  
Then since the inclusion map $C(X)\hookrightarrow \cT_E$ is a $KK$-equivalence, the 6-term exact sequence for 
$$0\to \K_E\otimes M_{(m)}\to \cT_E\otimes M_{(m)}\to \cO_E\otimes M_{(m)}\to 0,$$
shows that $(1-[E])K^0(X)\otimes \Z_{(m)}$ coincides the above kernel. 

Since $[E]\in mr+1+\tilde{K}^0(X)$ and $r$ is invertible in $K^0(X)\otimes \Z_{(m)}$, 
the above 6-term exact sequence takes the following form now:  
$$
{\begin{CD}
K^0(X)\otimes \Z_{(m)}@>(m+b)\times>>K^0(X)\otimes \Z_{(m)}  @>(\iota_B\otimes \id_{\M_{(m)}})_*>>K_0(B\otimes \M_{(m)})\\
@AAA@. @VVV\\
K_1(B\otimes \M_{(m)})@<(\iota_B\otimes \id_{\M_{(m)}})_*<< K^1(X)\otimes \Z_{(m)} @<(m+b)\times<<K^1(X)\otimes \Z_{(m)}
\end{CD}},
$$
where $b\in \tilde{K}^0(X)\otimes \Z_{(m)}$. 
The assumption $\Tor(H^*(X),\Z_m)=0$ with the Atiyah-Hirzebruch spectral sequence implies $\Tor(K^*(X),\Z_m)=\{0\}$ 
(see for example the proof of \cite[Theorem 5.3]{D12}), 
and the natural map from $K^*(X)\otimes \Z_{(m)}$ to $K^*(X)\otimes \Q$ is injective. 
Since $m+b$ is invertible in $K^0(X)\otimes \Q$, we get the following short exact sequence. 
\begin{equation}\label{exact}
0\to(m+b)(K^0(X)\otimes \Z_{(m)})\to K^0(X)\otimes \Z_{(m)}\xrightarrow{(\iota_B\otimes \id_{\M_{(m)}})_*}K_0(B\otimes \M_{(m)})\to 0.
\end{equation}

\begin{lemma} Let $m$ and $X$ be as above, and let $B_1$ and $B_2$ be locally trivial continuous fields of $\cO_{m+1}$ over $X$. 
If $B_1\otimes \K$ and $B_2\otimes \K$ are isomorphic as continuous fields over $X$, so are $B_1$ and $B_2$. 
\end{lemma}

\begin{proof} Let $e\in \K$ a minimal projection, and let $\varphi:B_2\otimes \K\to B_1\otimes \K$ be an isomorphism of $C(X)$-algebras. 
Then $B_2$ is isomorphic to $\varphi(1\otimes e)B_1\varphi(1\otimes e)$ as a $C(X)$-algebra. 
Since $B_1\cong B_1\otimes \cO_\infty$, the projection $\varphi(1\otimes e)$ is equivalent to the one of the form $p\otimes e$. 
Thus from the beginning we may and do assume $B_2=pB_1p$, and we identify $K_0(B_2)$ with $K_0(B_1)$. 
We simply denote $\iota_i=\iota_{B_i}$ for $i=1,2$. 

From (\ref{exact}), the map $(\iota_1\otimes \id_{\M_{(m)}})_*$ is surjective, and there exists $k\in K^0(X)\otimes \Z_{(m)}$ with $(\iota_1\otimes \id_{\M_{(m)}})_*(k)
=[p\otimes 1_{\M_{(m)}}]_0$. 
We can choose projections $p',p''\in C(X)\otimes \M_{(m)}\otimes \M_N$ satisfying $k=[p']_0-[p'']_0$, 
where $\M_N$ is the $N$-by-$N$ matrix algebra. 
Let $\varphi:C(X)\otimes \M_{(m)}\to C(X)\otimes\M_{(m)}\otimes \M_{(m)}$ be the embedding map given by 
$\varphi(f\otimes x)=f\otimes 1\otimes x$. 
Then $\varphi_*:K_0(C(X)\otimes \M_{(m)})\to K_0(C(X)\otimes \M_{(m)}\otimes \M_{(m)})$ is an isomorphism. 
Let $m:B_1\otimes C(X)\to B_1$ be the multiplication map. 
Then for any projection $q$ in $C(X)\otimes \M_{(m)}$, we get 
$$(\iota_2\otimes \id_{\M_{(m)}\otimes \M_{(m)}})_*\circ \varphi_*([q]_0)=(m\otimes \id_{\M_{(m)}\otimes \M_{(m)}})_*([p_{13}q_{24}]_0),$$
where $p_{13}$ and $q_{24}$ in $B_1\otimes C(X)\otimes \M_{(m)}\otimes \M_{(m)}$ are understood in an obvious sense.  
Let $d:C(X)\otimes C(X)\to C(X)$ be the homomorphism induced by restriction  to the diagonal.  
Then we have a commutative diagram 
$$\begin{CD}
C(X)\otimes C(X)@>d>>C(X)\\
@V \iota_1\otimes\id_{C(X)} VV @VV\iota_1 V\\
B_1\otimes C(X)@>m>> B_1
\end{CD}, 
$$ 
and 
\begin{align*}
\lefteqn{(m\otimes \id_{\M_{(m)}\otimes \M_{(m)}})_*([p_{13}q_{24}]_0)} \\
 &=(m\otimes \id_{\M_{(m)}\otimes \M_{(m)}})_*\circ(\iota_1\otimes \id_{C(X)}\otimes \id_{\M_{(m)}\otimes \M_{(m)}}) ([p'_{135}q_{24}]_0)\\
 &\quad -(m\otimes \id_{\M_{(m)}\otimes \M_{(m)}})_*\circ(\iota_1\otimes \id_{C(X)}\otimes \id_{\M_{(m)}\otimes \M_{(m)}}) ([p''_{135}q_{24}]_0)\\
 &=(\iota_1\otimes \id_{\M_{(m)}\otimes \M_{(m)}})_*\circ (d\otimes \id_{ \M_{(m)}\otimes \M_{(m)}  })_*([p'_{135}q_{24}]_0)\\
 &\quad -(\iota_1\otimes \id_{\M_{(m)}\otimes \M_{(m)}})_*\circ (d\otimes \id_{ \M_{(m)}\otimes \M_{(m)}  })_*([p''_{135}q_{24}]_0).
\end{align*}
Identifying $K_0(C(X)\otimes \M_{(m)})$ and $K_0(C(X)\otimes \M_{(m)}\otimes \M_{(m)})$ 
with $K^0(X)\otimes \Z_{(m)}$, we get 
$$(d\otimes \id_{\M_{(m)}\otimes \M_{(m)}})_*([p'_{135}q_{24}])
 -(d\otimes \id_{\M_{(m)}\otimes \M_{(m)}})_*([p''_{135}q_{24}])=k[q]_0,$$
and 
$$(\iota_2\otimes \id_{\M_{(m)}})_*(x)=(\iota_1\otimes \id_{\M_{(m)}})(kx),$$
for all $x\in K^0(X)\otimes \Z_{(m)}$. 
From this, we get the inclusion $\ker (\iota_1\otimes \id_{\M_{(m)}})_*\subset \ker (\iota_2\otimes \id_{\M_{(m)}})_*$. 
By switching the roles of $B_1$ and $B_2$, we get equality of the two sets, and $B_1$ and $B_2$ are isomorphic 
as continuous fields. 
\end{proof}

\begin{proof}[Proof of Theorem \ref{Cunzt algebra case}]
Let $m$ and $G$ be as in Theorem \ref{Cunzt algebra case}, and let $\alpha$ and $\beta$ be outer $G$ actions on $\cO_{m+1}$. 
If $\alpha$ and $\beta$ are cocycle conjugate, we have $\cP_{\salpha}\cong \cP_{\sbeta}$, and Lemma \ref{de-stabilization} implies 
$\cP_{\alpha\otimes \id_\K}\cong \cP_{\beta\otimes \id_\K}$. 
Thus the above lemma shows $\cP_\alpha\cong \cP_\beta$, and the map $\cB$ is well-defined. 

Assume $\cB([\alpha])=\cB([\beta])$, that is $\cP_\alpha\cong \cP_\beta$. 
Then $\cP_{\salpha}^0\cong \cP_{\sbeta}^0$, and Theorem \ref{deformation} implies that $\alpha$ and $\beta$ are cocycle conjugate. 
Thus $\cB$ is injective. 

It remains to show that $\cB$ is surjective. 
Let $\cP$ be a principal $\Aut(\cO_{m+1})$-bundle over $BG$, and let $B=\Gamma(\cP\times_{\Aut(\cO_{m+1})}\cO_{m+1})$ be the corresponding 
continuous field of $\cO_{m+1}$ over $BG$. 
Let $\cP'=\cP\times_{\Aut(\cO_{m+1})}\Aut(\cO_{m+1}\otimes \K)_0$.   
Then Corollary \ref{DRT0} shows that there exists an outer $G$-action $\gamma$ on $\cO_{m+1}\otimes \K$ with 
$\gamma_G\subset \Aut(\cO_{m+1}\otimes \K)_0$ satisfying 
$\cP'\cong \cP_\gamma^0$, which induces an isomorphism 
$\varphi:B\otimes \K\to M_\gamma$. 
Here we emphasize that $\varphi$ comes from an isomorphism of principal $\Aut(\cO_{m+1}\otimes \K)_0$-bundles, and the group 
$\Aut(\cO_{m+1}\otimes \K)_0$ acts on $K_0(\cO_{m+1})$ trivially. 

Choosing partial isometries $v_g\in \cO_{m+1}\otimes \K$ satisfying $v_g\gamma_g(1\otimes e)v_g^*=1\otimes e$, 
we get a cocycle action $(\theta,w)$ of $G$ on $\cO_{m+1}$ with $\theta_g(x)\otimes e=v_g\gamma_g(x\otimes e)v_g^*$ and 
$w_{g,h}\otimes e=v_g\alpha_g(v_h)v_{gh}^*$. 
Since $\hat{\hat{\theta}}$ as defined in Subsection 3.2 is a cocycle perturbation of $\gamma$, we may assume 
$\hat{\hat{\theta}}=\gamma$. 
Since the generalized mapping torus $M_{\hat{\hat{\theta}}}$ has a projection $\varphi(1\otimes e)$, 
Theorem \ref{cocycle actions} and Remark \ref{promap} imply that $(\theta,w)$ is equivalent to a genuine action, say $\nu$, and 
$\nu\otimes \id_\K$ is a cocycle perturbation of $\gamma$. 
Thus the previous lemma implies $\cP\cong \cP_\nu$, and $\cB$ is surjective. 
\end{proof}

\begin{remark}
When the Hirsch length $h(G)$ is less than or equal to 3, the map $\cB$ is always bijective regardless of 
the condition $\Tor(H^*(BG),\Z_m)=\{0\}$ because the map $\alpha\mapsto \fo^2_{\alpha,\id}$ induces a bijection 
between the cocycle conjugacy classes of outer $G$-actions on $\cO_{m+1}$ and the space of primary obstructions 
$H^2(BG,\pi_1(\Aut(\cO_{m+1})))$ for the existence of a section of $\cP_\alpha$  
(see Theorem \ref{obstruction} and \cite[Corollary 8.17]{IMI}).  
\end{remark}
\section{Appendix}\label{Appendix}
\subsection{Proof of Theorem \ref{model}}
We denote by $\{e_i\}_{i=1}^n$ the canonical basis of $\R^n$. 
Since matrices naturally act on row vectors on the right, it is more convenient to describe our computation 
using the right action convention, and we adopt it switching from the left action convention in the main body of the paper.  

We first recall that if $n:=h(G)\geq 5$, free and cocompact actions of $G$ on $\R^n$ are unique 
up to topological conjugacy. 
Indeed, assume that $\mu_1$ and $\mu_2$ are such actions. 
\cite[Theorem 15]{De02} shows that there exists a free and cocompact polynomial action $\nu$ of $G$ on $\R^n$ 
of bounded degree. 
Since $\R^n/\mu_i$ and $\R^n/\nu$ are $K(G,1)$ spaces, they are homotopy equivalent, and  
\cite[Theorem 15B1,(c)]{W99} shows that they are homeomorphic.  
Thus there exist two group automorphisms $\tau_1, \tau_2\in \Aut(G)$ such that 
$\mu_i$ is topologically conjugate to $\nu_{\tau_i(\cdot)}$ for $i=1,2$.   
Since any two free and cocompact polynomial actions of bounded degree are polynomially conjugate 
by \cite[Theorem 1]{BD2002}, we see that $\nu_{\tau_1(\cdot)}$ and $\nu_{\tau_1(\cdot)}$ 
are topologically conjugate, and so are $\mu_1$ and $\mu_2$. 

To prove Theorem \ref{model}, it suffices to show the statement for $h(G)\leq 5$. 
Indeed, assume that the statement is true for $h(G)\leq 5$.  
Assume that $G$ is a poly-$\Z$ group with $n=h(G)>5$ and we have constructed an action of $N=G_{n-1}$ on $\R^{n-1}$ 
satisfying the condition in Theorem \ref{model} for $N$ in place of $G$. 
We choose $\xi\in G$ such that $G$ is generated by $N$ and $\xi$. 
Then two actions of $N$ on $\R^{n-1}$ given by 
$(x,h)\mapsto x\cdot h$ and $(x,h)\mapsto x\cdot h^{\xi}$ for $x\in \R^{n-1}$ and $h\in N$ are topologically conjugate, 
and there exists a homeomorphism $\varphi$ of $\R^{n-1}$ satisfying 
$\varphi(x\cdot h)=\varphi(x)\cdot h^{\xi}$. 
Setting $(x,y)\cdot h=(x\cdot h,y)$ and $(x,y)\cdot \xi=(\varphi(x),y+1)$ 
for $x\in \R^{n-1}$ and $y\in \R$, we get a desired action of $G$ on $\R^n$. 
Thus by induction we get the statement for $G$. 

Furthermore, the same argument as above with the uniqueness result \cite[Theorem 1]{BD2002} implies that 
it suffices to show Theorem \ref{model} for $h(G)\leq 4$ with an extra assumption 
that the action is given by a polynomial map with bounded degree. 
Thus we prove Theorem \ref{model} by listing all poly-$\Z$ groups with $h(G)\leq 4$ and verifying the condition. 

For $G_1$ there is only one possibility $G_1=\langle\xi_1\rangle$.
In this case, we can put $x\cdot \xi_1=x+1$. 

For $G_2$ we have two possibilities, $G_2\cong \Z^2$ and $\xi_1^{\xi_2}=\xi_1^{-1}$. 
In the first case, we can put $(x,y)\cdot \xi_2=(x,y)+e_2$, and in the second case, we can put 
$(x,y)\cdot \xi_2=(-x,y+1)$. 

To go further, the following easy statement is useful. 

\begin{lemma}\label{commutator}
The commutator subgroup $[G_i,G_i]$ is generated by $[G_{i-1},G_{i-1}]$ and 
$\{h^{\xi_i}h^{-1};\;h\in G_{i-1}\}$ as a normal subgroup of $G_{i-1}$. 
\end{lemma}

Now we list all poly-$\Z$ groups of Hirsch length 3. 

\begin{lemma} \label{n=3}  Let $G$ be a poly-$\Z$ group with $h(G)=3$. 
Then $G$ is isomorphic to either one of the following isomorphism classes of groups: 
\begin{itemize}
\item[$(1)$] $G_2=\Z^2$ and $G=\Z^2\rtimes_A\Z$ with $A\in GL(2,\Z)$, 
that is $\xi_1^{\xi_3}=\xi_1^a\xi_2^b$, $\xi_2^{\xi_3}=\xi_1^c\xi_2^d$ with $A=\left(
\begin{array}{cc}
a &b  \\
c &d 
\end{array}
\right)
.$
For the condition of Theorem \ref{model} to be satisfied, we can define $(x,y,z)\cdot \xi_i=(x,y,z)+e_i$ for $i=1,2$ and 
$(x,y,z)\cdot \xi_3=((x,y)A,z)+e_3$. 
\item[$(2)$] $G_2=\langle \xi_1,\xi_2|\; \xi_1^{\xi_2}=\xi_1^{-1}\rangle$, and 
$\xi_1\xi_3=\xi_3\xi_1$, $\xi_2^{\xi_3}=\xi_2^{-1}$. 
For the condition of Theorem \ref{model} to be satisfied, we can define $(x,y,z)\cdot \xi_1=(x,y,z)+e_1$, 
$(x,y,z)\cdot \xi_2=(-x,y,z)+e_2$, and $(x,y,z)\cdot \xi_3=(x,-y,z)+e_3$.  
\item[$(3)$] $G_2=\langle \xi_1,\xi_2|\; \xi_1^{\xi_2}=\xi_1^{-1}\rangle$, and 
$\xi_1\xi_3=\xi_3\xi_1$, $\xi_2^{\xi_3}=\xi_2^{-1}\xi_1$. 
For the condition of Theorem \ref{model} to be satisfied, we can define $(x,y,z)\cdot \xi_1=(x,y,z)+e_1$, 
$(x,y,z)\cdot \xi_2=(-x,y,z)+e_2$, and $(x,y,z)\cdot \xi_3=(x,-y,z)+\frac{1}{2}e_1+e_3$.  
\end{itemize}
\end{lemma}

\begin{proof} If $G_2\cong \Z^2$, $G$ is isomorphic to a group in (1). 

Assume $\xi_1^{\xi_2}=\xi_1^{-1}$ now. 
Then $[G_2,G_2]=\langle \xi_1^2\rangle$ and $Z(G_2)=\langle \xi_2^2\rangle$. 
Since $\xi_1$ is a unique square root of $\xi_1^2$ in $G_2$, the subgroup $\langle \xi_1\rangle$ 
is characteristic in $G_2$. 
On the other hand, the set of square roots of $\xi_2^2$ in $G_2$ is $\{\xi_2\xi_1^a\}_{a\in \Z}$. 
Thus $\xi_1^{\xi_3}=\xi_1^{\epsilon_1}$ and $\xi_2^{\xi_3}=\xi_2^{\epsilon_2}\xi_1^a$ with 
$\epsilon_1,\epsilon_2\in \{1,-1\}$ and $a\in \Z$. 

Note that for fixed $G_2$, the splitting $G=G_2\rtimes \langle\xi_3\rangle$ is not unique and we have freedom 
to replace $\xi_3$ with $h\xi_3^{\pm 1}$, where $h\in G_2$. 
Thus by replacing $\xi_3$ with $\xi_2\xi_3$ if necessary, we may assume $\epsilon_1=1$. 
If $\epsilon_2=1$, the group $\langle \xi_1,\xi_3\rangle$ is a normal subgroup of $G$ isomorphic to $\Z^2$, 
and $G$ is of the form $\langle \xi_1,\xi_3\rangle\rtimes_{\xi_2}\Z\cong \Z^2\rtimes_A\Z$. 
This case is reduced to the case (1). 
Thus we may assume that $\epsilon_2=-1$. 
By replacing $\xi_3$ with $\xi_1^b\xi_3$ with an appropriate $b\in \Z$, we may assume that 
$a=0$ or $a=1$, and we are left with the two cases (2) and (3).  
\end{proof}

To prove Theorem \ref{model}, it suffices to show that for any group $G$ with its action on $\R^3$ in the 
above lemma and for any $\theta\in \Aut(G)$, 
there exists a polynomial homeomorphism $\varphi:\R^3\to\R^3$ such that the degree of $\varphi^n$ 
is at most two for any $n\in \Z$ and 
\begin{equation}\label{relation}
\varphi(x\cdot \xi_i)=\varphi(x)\cdot \xi_i^\theta
\end{equation} 
for any $x\in \R^3$ and $i=1,2,3$. 
Indeed, any poly-$\Z$ group with Hirsch length 4 is of the form $\Gamma=G\rtimes_\theta\Z$, 
and we can define a polynomial action of $\Gamma$ on $\R^4$ with degree at most 2 satisfying the 
condition of Theorem \ref{model} by $(x,y)\cdot g=(x\cdot g,y)$ for $g\in G$ and 
$(x,y)\cdot \xi_4=(\varphi(x),y+1)$ for $x\in \R^3$ and $y\in \R$. 

Before starting case-by-case analysis, we state an easy lemma. 
\begin{lemma} Let $A\in GL(2,\Z)$. 
\begin{itemize}
\item[$(1)$] If $\det A=-1$ and $\Tr A=0$, then $A$ is conjugate to 
either $\left(
\begin{array}{cc}
-1 &0  \\
 0& 1
\end{array}
\right)
$ or 
$\left(
\begin{array}{cc}
0 &1  \\
1 &0 
\end{array}
\right)
$ in $GL(2,\Z)$.
\item[$(2)$] If $\det A=1$ and $\Tr A=2$, then $A$ is conjugate to $\left(
\begin{array}{cc}
1 &m  \\
0 &1 
\end{array}
\right)
$ with $m\in \N\cup \{0\}$ in $SL(2,\Z)$. 
\end{itemize}  
\end{lemma}

Now we will finish the proof of Theorem \ref{model} by verifying Eq.(\ref{relation}) for each case. 

Case I.
Assume $G=\Z^2\rtimes_A\Z$ with $A\in GL(2,\Z)$ satisfying $(\Tr A,\det A)\neq (2,1),(0,-1)$. 
Then $\det(I-A)\neq 0$. 
Identifying $G_2$ with $\Z^2$, by Lemma \ref{commutator} we have $[G,G]=\{(x,y)(I-A)\in \Z^2;\; (x,y)\in G_2\}$, which has a finite index in $G_2$. 
Therefore 
$$G_2=\{g\in G;\; \exists k\in \Z\setminus\{0\},\; g^k\in [G,G]\},$$ 
and it is a characteristic subgroup of $G$. 
Thus there exist $B=\left(
\begin{array}{cc}
p &q  \\
r &s 
\end{array}
\right)
\in GL(2,\Z)$ and $h=\xi_1^{h_1}\xi_2^{h_2}\in G_2$ satisfying $\xi_1^\theta=\xi_1^p\xi_2^q$, 
$\xi_2^\theta=\xi_1^r\xi_2^s$, and $\xi_3^\theta=\xi_3^\epsilon h$ with $\epsilon\in \{1,-1\}$. 
The relations among $\xi_1,\xi_2,\xi_3$ imply $BA^\epsilon=AB$. 

On the other hand, any 
$$(B,\epsilon,h)\in GL(2,\Z)\times \{1,-1\}\times G_2$$ with $BA^\epsilon=AB$ 
gives rise to $\theta\in \Aut(G)$ as above. 
 
Letting $\varphi:\R^3\to\R^3$ be an invertible affine map given by  
$$\varphi((x,y,z))=((x,y)B,\epsilon z)+((h_1,h_2)(I-A^\epsilon)^{-1},0),$$ 
we can verify (\ref{relation}).

Case II. Assume $G=\Z^2\rtimes_A\Z$ with $A=\left(
\begin{array}{cc}
-1 &0  \\
0 &1 
\end{array}
\right)$. 
Then we have $[G,G]=\langle \xi_1^2\rangle$ and $Z(G)=\langle \xi_2,\xi_3^2\rangle$.  
Since $\xi_1$ is the only square root of $\xi_1^2$ in $G$, the subgroup $\langle \xi_1\rangle $ 
is characteristic in $G$. 
Thus there exists $\epsilon\in \{1,-1\}$ and $\left(
\begin{array}{cc}
a &b  \\
c &d 
\end{array}
\right)\in GL(2,\Z)$ satisfying $\xi_1^\theta=\xi_1^\epsilon$, $\xi_2^\theta=\xi_2^a\xi_3^{2b}$, and 
$(\xi_3^\theta)^2=\xi_2^c\xi_3^{2d}$. 
This implies that $c$ is even and $\xi_3^\theta=\xi_3^d\xi_2^{c/2}\xi_1^\alpha$ with $\alpha\in \Z$. 

On the other hand, for any 
$$(\left(
\begin{array}{cc}
a &b  \\
c &d 
\end{array}
\right),\epsilon,\alpha)\in GL(2,\Z)\times \{1,-1\}\times \Z,$$
with even $c$, there exists $\theta\in \Aut(G)$ as above. 
Letting $\varphi:\R^3\to\R^3$ be an invertible affine map given by 
$$\varphi(x,y,z)=(x,y,z)\left(
\begin{array}{ccc}
\epsilon &0 &0  \\
0 &a &2b  \\
0 &\frac{c}{2} &d 
\end{array}
\right)
+(\frac{\alpha}{2},0,0),$$
we can verify (\ref{relation}).

Case III. Assume $G=\Z^2\rtimes_A\Z$ with $A=\left(
\begin{array}{cc}
0 &1  \\
1 &0 
\end{array}
\right)$. 
Then $[G,G]=\langle \xi_1\xi_2^{-1}\rangle$ and $Z(G)=\langle \xi_1\xi_2,\xi_3^2\rangle$. 
Thus there exists $\left(
\begin{array}{cc}
a &b  \\
c &d 
\end{array}
\right)\in GL(2,\Z)$ and $\epsilon\in \{1,-1\}$ satisfying 
$(\xi_1\xi_2^{-1})^\theta=(\xi_1\xi_2^{-1})^\epsilon$, 
$(\xi_1\xi_2)^\theta=(\xi_1\xi_2)^a\xi_3^{2b}$, and $(\xi_3^{\theta})^2=(\xi_1\xi_2)^c\xi_3^{2d}$, 
which implies $(\xi_1^\theta)^2=\xi_1^{a+\epsilon}\xi_2^{a-\epsilon}\xi_3^{2b}$ and 
$(\xi_2^\theta)^2=\xi_1^{a-\epsilon}\xi_2^{a+\epsilon}\xi_3^{2b}$. 
This shows that $\xi_1^\theta$ is of the form $\xi_3^b\xi_2^p\xi_1^q$. 
If $b$ were odd, we would have $(\xi_1^\theta)^2=\xi_3^{2b}\xi_2^{p+q}\xi_1^{p+q}$, which is 
contradiction. 
Thus $b$ is even and $\xi_1^\theta=\xi_3^b\xi_2^{\frac{a-\epsilon}{2}}\xi_1^{\frac{a+\epsilon}{2}}$, 
$\xi_2^\theta=\xi_3^b\xi_2^{\frac{a+\epsilon}{2}}\xi_1^{\frac{a-\epsilon}{2}}$. 
Since $ad-bc=\pm 1$ and $b$ is even, $d$ is odd and we get 
$\xi_3^\theta=\xi_3^d\xi_2^{\frac{c-\alpha}{2}}\xi_1^{\frac{c+\alpha}{2}}$
with $c+\alpha\in 2\Z$. 

On the other hand, for any 
$$(\left(
\begin{array}{cc}
a &b  \\
c &d 
\end{array}
\right),\epsilon,\alpha)\in GL(2,\Z)\times \{1,-1\}\times \Z,$$
with even $b$ and $c+\alpha$, we
can get $\theta\in \Aut(G)$ as above. 
Letting $\varphi:\R^3\to\R^3$ be an invertible affine map given by 
$$\varphi(x,y,z)=(x,y,z)\left(
\begin{array}{ccc}
\frac{a+\epsilon}{2} &\frac{a-\epsilon}{2} &b  \\
\frac{a-\epsilon}{2} &\frac{a+\epsilon}{2} &b  \\
\frac{c}{2} &\frac{c}{2} &d 
\end{array}
\right)
+(\frac{\alpha}{2},0,0),$$
we can verify (\ref{relation}).

Case IV. Let $G=\Z^2\rtimes_A\Z$ with $A=\left(
\begin{array}{cc}
1 &m  \\
0 &1 
\end{array}
\right)$ and $m\in \N\cup\{0\}$. 
If $m=0$, we have $G=\Z^3$ and $\theta$ is given by $B\in GL(3,\Z)$. 
Hence $\varphi(x,y,z)=(x,y,z)B$ works. 
Assume $m\neq 0$. 
Then $G$ is a cocompact lattice of the Heisenberg group 
$$H_3=\{\left(
\begin{array}{ccc}
1 &r &s  \\
0 &1 &t  \\
0 &0 &1 
\end{array}
\right)
;\; r,s,t\in \R^3\},$$
with embedding 
$$\xi_3^c\xi_2^b\xi_1^a=\left(
\begin{array}{ccc}
1 &a &\frac{b}{m}  \\
0 &1 &c  \\
0 &0 &1 
\end{array}
\right).$$
We identify $(x,y,z)\in \R^3$ with $\left(
\begin{array}{ccc}
1 &x &\frac{y}{m}  \\
0 &1 &z  \\
0 &0 &1 
\end{array}
\right)\in H_3$, and let $G$ act on $H_3$ by right multiplication. 
This action is identified with the action given in Lemma \ref{n=3},(1). 
Thus we get $\varphi(h\cdot g)=\varphi(h)\cdot g^\theta$ for any $h\in H_3$ and $g\in G$.
Thanks to the Malcev rigidity theorem (see \cite[Theorem 2.11]{Rag72}), every $\theta\in \Aut(G)$ uniquely lifts to $\varphi\in \Aut(H_3)$. 
Since $H_3$ is a central extension 
$$0\to \R\to H_3\to \R^2\to0,$$
with a 2-cocycle given by $\omega((a_1,c_1),(a_2,c_2))=a_1c_2$, we can see that 
every automorphism of $H_3$ is a polynomial map of degree at most 2. 

Case V. Assume 
$G=\langle \xi_1,\xi_2,\xi_3|\;\xi_1^{\xi_2}=\xi_1^{-1},\;[\xi_1,\xi_3]=e,\; \xi_2^{\xi_3}=\xi_2^{-1} \rangle$. 
Then $[G,G]=\langle \xi_1^2,\xi_2^2\rangle$ and $Z(G)=\langle \xi_3^2\rangle$. 
Note that $G_2$ is characteristic as we have 
$$G_2=\{g\in G;\; \exists k\in \Z\setminus \{e\},\; g^k\in [G,G]\}.$$ 
Thus there exist $\epsilon_1,\epsilon_2,\epsilon_3\in \{1,-1\}$ and $a,b,c\in \Z$ satisfying 
$\xi_1^\theta=\xi_1^{\epsilon_1}$, $\xi_2^\theta=\xi_2^{\epsilon_2}\xi_1^a$, and $\xi_3^\theta=\xi_3^{\epsilon_3}\xi_2^b\xi_1^c$. 
Since $\xi_1^\theta$ commutes with $\xi_3^\theta$, the number $b$ is even, and so $(\xi_3^\theta)^2=\xi_3^{2\epsilon_3}\xi_1^{2c}$ 
implies $c=0$ as $(\xi_3^\theta)^2\in Z(G)$. 

On the other hand, for any 
$$(\epsilon_1,\epsilon_2,\epsilon_3,a,b)\in \{1,-1\}^3\times \Z^2$$
with even $b$, we get $\theta\in \Aut(G)$ as above. 
Letting $\varphi:\R^3\to\R^3$ be an invertible affine map given by 
$$\varphi(x,y,z)=(x,y,z)\left(
\begin{array}{ccc}
\epsilon_1 &0 &0  \\
0 &\epsilon_2 &0  \\
0&0 &\epsilon_3 
\end{array}
\right)
+(\frac{a}{2},\frac{b}{2},0),$$
we can verify (\ref{relation}). 

Case VI. Assume 
$G=\langle \xi_1,\xi_2,\xi_3|\; \xi_1^{\xi_2}=\xi_1^{-1},\; [\xi_1,\xi_3]=e,\;\xi_2^{\xi_3}=\xi_2^{-1}\xi_1 \rangle$. 
Then we have $[G,G]=\langle \xi_1^2, \xi_2^2\xi_1\rangle$ and $Z(G)=\langle \xi_1^{-1}\xi_3^2\rangle$. 
Since $G/[G,G]$ is abelian, the set 
$$\{g\in G;\; \exists k\in \Z\setminus \{e\},\; g^k\in [G,G]\},$$
is a characteristic subgroup of $G$, which coincides with $G_2$ as $\xi_1,\xi_2^4\in [G,G]$. 
Thus there exists $\epsilon_1,\epsilon_2,\epsilon_3\in \{1,-1\}$ and $a\in \Z$ satisfying
$\xi_1^\theta=\xi_1^{\epsilon_1}$, $\xi_2^\theta=\xi_2^{\epsilon_2}\xi_1^a$, and 
$\xi_3^\theta\in \xi_3^{\epsilon_3}G_2$. 
Since $(\xi_1^{-1}\xi_3^2)^\theta=(\xi_1^{-1}\xi_3^2)^{\epsilon_3}$, 
we see that there exists $b\in 2\Z$ satisfying 
$\xi_3^\theta=\xi_3^{\epsilon_3}\xi_2^b\xi_1^{\frac{\epsilon_1-\epsilon_3}{2}}$. 

On the other hand, for any 
$$(\epsilon_1,\epsilon_2,\epsilon_3,a,b)\in \{1,-1\}^3\times \Z^2$$
with even $b$, we get $\theta\in \Aut(G)$ as above. 
Letting $\varphi:\R^3\to\R^3$ be an invertible affine map given by 
$$\varphi(x,y,z)=(x,y,z)\left(
\begin{array}{ccc}
\epsilon_1 &0 &0  \\
0 &\epsilon_2 &0  \\
0&0 &\epsilon_3 
\end{array}
\right)
+(\frac{a}{2},\frac{b}{2},0),$$
we can verify (\ref{relation}).

\subsection{The map $\Pi_{A,X}$ comes from a weak homotopy equivalence}
Ulrich Pennig kindly informed us that the map $\Pi_{A,X}$ in Theorem \ref{weak homotopy} actually comes from 
a continuous map from $U(A_\flat)$ to $\Omega\Aut(A\otimes \K)$. 
We would like to thank him for his courtesy to allow us to include his argument here. 

Recall that $U(A_\flat)$ is homotopy equivalent to a CW complex as it is homotopic to an open subset of a Banach space (see \cite[Corollary IV.5.5]{LW}). 
Thus we choose a CW complex $Y$ and homotopy equivalence $f:Y\to U(A_\flat)$. 
Let $\{Y_i\}_i$ be the set of finite subcomplexes of $Y$, which is a direct set by inclusion. 
Thanks to the homotopy extension property of any subcomplexes $Y'\subset Y''$ of $Y$, the map 
$$[Y,Z]\to \lim_{i}[Y_i,Z]$$
is surjective for any space $Z$. 
Since $\{\Pi_{A,Y_i}([f|_{Y_i}])\}_i$ is an element in $\lim_i[Y_i,\Omega \Aut(\sA)]$ thanks to the naturality of 
$\Pi_{A,X}$ in $X$, there exists a map $g:Y\to \Omega\Aut(\sA)$ satisfying 
$[g|_{Y_i}]=\Pi_{A,Y_i}([f|_{Y_i}])$ for any $Y_i$. 
We show that $g$ is a weak homotopy equivalence. 

Assume we have $[h]\in \pi_n(Y)$ with  $g_*[h]=0$. 
Then we may assume that $h$ factors through a finite subcomplex, say $Y_i$, and 
$h=\iota_i\circ h_1$, where $h_i:S^n\to Y_i$, and $\iota_i$ is the inclusion map of $Y_i$ into $Y$. 
Thus 
$$0=[g\circ h]=[g\circ \iota_i\circ h_1]=h_1^*[g|_{Y_i}]=h_1^*(\Pi_{A,Y_i}([f|_{Y_i}]))=\Pi_{A,n}([f|_{Y_i}\circ h_1])=\Pi_{A,n}([f\circ h]),$$
and $[f\circ h]=0$. 
Since $f$ is a homotopy equivalence, we get $[h]=0$. 
This shows that the map $g_*:\pi_n(Y)\to \pi_n(\Omega\Aut(\sA))$ is injective. 

Let $[k]\in \pi_n(\Omega\Aut(\sA))$. 
Since $\Pi_{A,n}$ is an isomorphism, and $f$ is a homotopy equivalence, there exists $[k_1]\in \pi_n(Y)$ with 
$\Pi_{A,n}(f_*[k_1])=[k]$. 
As before, we may assume that the map $k_1$ factors through a finite subcomplex, and we get $g_*[k_1]=[k]$, showing that 
the map $g_*:\pi_n(Y)\to \pi_n(\Omega\Aut(\sA))$ is surjective.  
This shows that $g$ is a weak homotopy equivalence. 


\end{document}